%% file: Caldararu-hom_fin.tex
\documentclass{amsart}
\numberwithin{equation}{section}

\usepackage{latexsym,amsmath,amstext,mathrsfs,amsthm,amscd,amsopn,verbatim,amsrefs}

\usepackage{graphicx}
\usepackage[dvips]{color}
\usepackage{epsfig}
\usepackage{psfrag}
\usepackage[all]{xy}
\usepackage{amscd}
\usepackage{amssymb}

\newtheorem{Thm}{Theorem}[section]
\newtheorem{Prop}[Thm]{Proposition}
\newtheorem{Lem}[Thm]{Lemma}
\newtheorem{Cor}[Thm]{Corollary}

\theoremstyle{remark}

\newtheorem*{Ack}{Acknowledgment}

\theoremstyle{definition}

\newtheorem{Def}[Thm]{Definition}
\newtheorem{Exa}[Thm]{Example}

\renewcommand{\bul}{{\bullet}}

\newcommand{\C}{{\rm C}}

\newcommand{\coord}{{\rm coord}}
\newcommand{\aff}{{\rm aff}}
\newcommand{\poly}{{\rm poly}}
\newcommand{\Hom}{\mathrm{Hom}}


\title[C{\u{a}}ld{\u{a}}raru's conjecture and Tsygan's formality]
{C{\u{a}}ld{\u{a}}raru's conjecture and Tsygan's formality}
\author{Damien Calaque}
\address[Damien Calaque]{Department of Mathematics, ETH Z\"urich, 8092 Z\"urich, Switzerland \\ 
on leave from Institut Camille Jordan, CNRS \& Universit\'e Lyon 1, France}
\email{damien.calaque@math.ethz.ch,calaque@math.univ-lyon1.fr}
\author{Carlo A. Rossi}
\address[Carlo A. Rossi]{MPIM Bonn, Vivatsgasse 7, 53111 Bonn, Germany}
\email{crossi@mpim-bonn.mpg.de}
\author{Michel van den Bergh}
\address[Michel Van den Bergh]{Departement WNI, Hasselt University, Agoralaan, 3590 Diepenbeek, Belgium}
\email{michel.vandenbergh@uhasselt.be} \thanks{This paper was partly written when the first
  author was a visitor at the Department of Mathematics of ETH (on leave of absence from Universit\'e Lyon 1) 
  thanks to the support of the European Union through a Marie Curie
  Intra-European Fellowship (contract number MEIF-CT-2007-042212). }
\thanks{The second author has been partially supported by the Funda\c{c}\~{a}o para a Ci\^{e}ncia e a Tecnologia 
(FCT / Portugal) while writing the present paper.}  
\thanks{The third author is a director of research at the FWO}

\begin{document}

\maketitle

\begin{abstract}
  In this paper we complete the proof of C{\u{a}}ld{\u{a}}raru's conjecture on the
  compatibility between the module structures on differential forms
  over poly-vector fields and on Hochschild homology over Hochschild
  cohomology. In fact we show that twisting with the square root of
  the Todd class gives an isomorphism of precalculi between these
  pairs of objects.

  Our methods use formal geometry to globalize the local formality
  quasi-isomorphisms introduced by Kontsevich and Shoikhet (the
  existence of the latter was conjectured by Tsygan). We also
rely on the fact - recently proved by the first two authors - that Shoikhet's
quasi-isomorphism is compatible with cap product after twisting with a
Maurer-Cartan element.
\end{abstract}
\tableofcontents
\section{Introduction and statement of the main results}\label{ref-1-0}
Throughout $k$ is a ground field of characteristic $0$. 
In this
introduction $(X,\mathcal O)$ is a ringed site\footnote{We work over
  sites instead of spaces to cover some additional cases which are
  important for algebraic geometry (like algebraic spaces and Deligne--Mumford stacks). Readers not interested in such generality may
  assume that $(X,\mathcal O)$ is just a ringed space.}  such
that~$\mathcal O$ is a sheaf of commutative $k$-algebras. We fix in
addition a Lie algebroid $\mathcal L$ over $(X,\mathcal O)$.

Roughly speaking a Lie algebroid is a sheaf of $\mathcal
O$-modules which is also a sheaf of Lie algebras which acts on
$\mathcal O$ by derivations. See \S\ref{ref-3.1-17}.
Standard examples of Lie algebroids are the tangent bundle on
a smooth manifold and the  holomorphic tangent bundle on a complex manifold.  
Readers not familiar with Lie algebroids are advised to think of
  $\mathcal L$ as a tangent bundle (holomorphic or not) for the rest
  of this introduction. Concepts like ``connection'' take their
  familiar meaning in this context. In fact: our main reason for working
in the setting of Lie algebroids is that these allow us to treat the
algebraic, holomorphic and $C^\infty$-cases in a uniform way. 
%
\subsection{The Atiyah and Todd class of a Lie algebroid}\label{ref-1.1-1}
\emph{From now on we make the
additional assumption that the Lie algebroid $\mathcal L$ is locally
free of rank $d$ as an $\mathcal O$-module.}

The Atiyah class $A(\mathcal L)\in
\operatorname{Ext}^1(\mathcal L,\mathcal L^\ast\otimes \mathcal L)=\mathrm H^1(X,\mathcal{L}^\ast \otimes \operatorname{End}_{\mathcal O}(\mathcal L))$ of  
$\mathcal L$ may for example be defined as the obstruction against the existence
of a global $\mathcal{L}$-connection on $\mathcal{L}$. See \S\ref{ref-6-104} for more
details. 

The
$i$-th scalar Atiyah class $a_i(\mathcal L)$ of $\mathcal L$ is
defined as
\begin{equation*}
a_i(\mathcal L)=\mathrm{tr}(\bigwedge^i A(\mathcal L))\in \mathrm H^i(X,\bigwedge^i \mathcal L^*),
\end{equation*} 
where  $\bigwedge^i$ is the map
\[
\bigwedge^i:\left( \mathcal L^\ast\otimes \operatorname{End}(\mathcal L)\right)^{\otimes i} \rightarrow  \bigwedge^i \mathcal L^* \otimes\operatorname{End}(\mathcal L)
\]
given by composition on $\operatorname{End}(\mathcal L)^{\otimes i}$ and
the exterior product on $(\mathcal L^\ast)^{\otimes i}$ and where
$\mathrm{tr}$ is the usual trace on $\mathrm{End}(\mathcal L)$,
extended linearly to  a map $ \bigwedge^i \mathcal L^*   \otimes  \operatorname{End}(\mathcal L)  \rightarrow   \bigwedge^i \mathcal L^*
$.

The Todd class $\mathrm{td}(\mathcal L)$ of $\mathcal L$ is 
derived from the Atiyah class $A(\mathcal L)$ by the following familiar formula:
\begin{equation}\label{ref-1.2-3}
\mathrm{td}(\mathcal L)=\det\!\left(\frac{A(\mathcal L)}{1-\exp\left(- A(\mathcal L)\right)}\right)\in \bigoplus_{i\geq 0}\mathrm H^i(X,\bigwedge^i\mathcal L^*),
\end{equation}
where the function
\begin{equation}
\label{ref-1.3-4}
q(x)=\frac{x}{1-\exp(-x)} 
\end{equation}
is extended to $\bigwedge \mathcal L^*\otimes_{\mathcal
  O}\mathrm{End}(\mathcal L)$ {\em via} its formal Taylor expansion. In
this way the Todd class $\mathrm{td}(\mathcal L)$ of $\mathcal L$ can
be expressed in terms of the scalar Atiyah classes of $\mathcal L$.

\subsection{Gerstenhaber algebras and precalculi}
\label{ref-1.2-5}
By definition a Gerstenhaber algebra is a graded vector space equipped with a
Lie bracket $[-,-]$ of degree zero and a commutative, associative cup product $\cup$
of degree one\footnote{Note that our grading conventions are shifted
with respect to the usual ones.} such that the Leibniz
rule is satisfied
\[
[a,b\cup c]=[a,b]\cup c+(-1)^{|a|(|b|+1)}b\cup [a,c]
\]
If $A$ is a Gerstenhaber algebra then a precalculus \cite{DTT2}  over $A$ is a quadruple $(A,M,\imath,\mathrm{L})$ where
$M$ is a graded vector space and $\imath:A\otimes M\rightarrow M$ and
$\mathrm{L}:A\otimes M\rightarrow M$ are linear maps of degree $1$
and $0$ respectively such that $\imath$ makes $M$ into an
$(A[-1],\cup)$-module and $\mathrm{L}$ makes $M$ into an
$(A,[-,-])$-Lie module and such that the following compatibilities
hold for $a,b\in A$
\begin{align}
\imath_a \mathrm{L}_{{b}} -(-1)^{(|a|+1)|b|}\mathrm{L}_b \imath_a&=\imath_{[a,b]}\\
\mathrm{L}_a\imath_b+(-1)^{|a|+1}\imath_a\mathrm{L}_b&=\mathrm{L}_{a\cup b}
\label{ref-1.5-6}
\end{align}
A precalculus is \emph{not} the same as a Gerstenhaber module. The
second equation in the previous display is not correct for a
Gerstenhaber module.

Below $\imath$ will be referred to as ``contraction'' and $\mathrm L$ as the 
``Lie derivative''. 
Furthermore we will often write $a\cap m$ for $\imath_a(m)$ and as such
refer to it as the ``cap product''.

\subsection{Poly-vector fields, poly-differential operators, differential forms and Hochschild chains in the Lie algebroid framework}\label{ref-1.3-7}
For a Lie algebroid $\mathcal L$ the sheaves
of $\mathcal L$-poly-vector fields and 
$\mathcal L$-differential forms are defined as
\[
T_\mathrm{poly}^\mathcal L(X)=\bigoplus_{n\geq -1} \bigwedge^{n+1}\mathcal L,\qquad\Omega^\mathcal L(X)=\bigoplus_{n\leq 0} \bigwedge^{-n}\mathcal L^*
\]
where the wedge products are taken over $\mathcal O_X$.

The sheaf $T_\mathrm{poly}^\mathcal L(X)$ becomes a sheaf of
Gerstenhaber algebras when endowed with the trivial differential, the Lie
algebroid version of the Schouten--Nijenhuis Lie bracket and the exterior
product.  Our grading convention is such that the Lie bracket and
wedge product are of degree $0$ and $1$ respectively. 

We equip $\Omega^\mathcal L(X)$ with the trivial(!)
differential\footnote{The De Rham differential $\mathrm d_L$ on
  $\Omega^\mathcal L(X)$ is not part of the precalculus structure.  In
  the operadic setting of \cite{DTT2}, $\mathrm d_L$ appears as a unary
  operation and not as a differential.}, and also with the contraction
operator and  Lie derivative 
with respect to\ $\mathcal L$-poly-vector fields.  In this way the pair
($T_\mathrm{poly}^\mathcal L(X)$, $\Omega^\mathcal L(X)$) becomes a
sheaf of precalculi. 
In our conventions the contraction operator and Lie derivative have degrees $1$ and
$0$ respectively.  

\medskip

The Lie algebroid generalization of the sheaf of $\mathcal
L$-poly-differential operators is denoted by $D_\mathrm{poly}^\mathcal
L(X)$~\cite{Xu,Cal}. It is the tensor algebra over $\mathcal{O}$ of the
universal enveloping algebra of $\mathcal L$ (see \S\ref{ref-3.3-33} below). 

The sheaf $D_\mathrm{poly}^\mathcal L(X)$ has similar properties as
the standard sheaf of poly-differential operators on $X$ (see e.g.\ \cite{K}). In particular
it is a differential graded Lie algebra (shortly, from now on, a DG-Lie algebra) 
and also a Gerstenhaber algebra \emph{up to homotopy}. For the 
definition of the differential, the Lie bracket (of degree $0$) and
the cup product (of degree $1$) see \S\ref{ref-3.3-33}.

\medskip

The sheaf of $\mathcal L$-Hochschild chains $C_\poly^\mathcal L(X)$
may be defined as the $\mathcal O$-dual of $D_\mathrm{poly}^\mathcal
L(X)$ (although we use a slightly different
approach). 
Furthermore there is a differential $\mathrm b_H$ as well as actions $\cap$,
$\mathrm L$ of $D_\mathrm{poly}^\mathcal
L(X)$ 
on
$C_\poly^\mathcal L(X)$
which make
the pair ($D_\mathrm{poly}^\mathcal L(X),C_\poly^\mathcal L(X)$)  into a precalculus up 
to homotopy.
We refer to \S\ref{ref-3.4-38}. 

\medskip

Finally, we recall that there is a Hochschild--Kostant--Rosenberg (HKR
for short) quasi-isomor\-phism from $T_\mathrm{poly}^\mathcal L(X)$ to
$D_\mathrm{poly}^\mathcal L(X)$; dually, there is a HKR
quasi-isomor\-phism from $C_\mathrm{poly}^\mathcal L(X)$ to
$\Omega^\mathcal L(X)$. As in the classical case where $\mathcal L$ is the
tangent bundle neither of these  HKR quasi-isomorphisms
is compatible with the Gerstenhaber and precalculus structures up to homotopy.

\subsection{Main results}\label{ref-1.4-8}
Now we consider the derived category $D(X)$
of sheaves of $k$-vector spaces over $X$. When equipped with the
derived tensor product this becomes a symmetric monoidal category.  Furthermore, viewed as objects in $D(X)$,
both $T_\mathrm{poly}^\mathcal L(X)$ and $D_\mathrm{poly}^\mathcal
L(X)$ are honest Gerstenhaber algebras
and their combination with
$\Omega^\mathcal
L(X)$ and $C_\mathrm{poly}^\mathcal L(X)$ yields precalculi. 

Our first main result relates the 
Todd class of a Lie
algebroid (as discussed in~\S\ref{ref-1.1-1}) to the failure of
the HKR isomorphisms to preserve these precalculi structures. 
\begin{Thm}\label{ref-1.1-9}
  Let $\mathcal L$ be a locally free Lie algebroid of rank $d$ over
  the ringed site $(X,\mathcal O_X)$.
  Then we have the following commutative diagram of precalculi in the category $D(X)$:
\begin{equation}\label{ref-1.6-10}
  \xymatrix{T_\mathrm{poly}^\mathcal L(X)
\ar[rr]^{\mathrm{HKR}\circ \iota_{\sqrt{\mathrm{td}(\mathcal L)}}}\ar@{~>}[d] 
&& 
D_\mathrm{poly}^\mathcal L(X)\ar@{~>}[d]\\
\Omega^\mathcal L(X) 
&& 
C_\mathrm{poly}^\mathcal L(X)
\ar[ll]^{(\sqrt{\mathrm{td}(\mathcal L)}\wedge-)\circ\mathrm{HKR}},} 
\end{equation}
where the vertical arrows indicate  \emph{actions}
and the horizontal arrows
are isomorphisms.
Here $\wedge$ denotes the left multiplication in $\Omega^\mathcal
L(X)$ and $\iota$ denotes the contraction action of $\Omega^\mathcal
L(X)$ on $T^\mathcal L(X)$.\footnote{Note that normally we
  view $\Omega^\mathcal L(X)$ as a module over $T^\mathcal L(X)$. In
  the definition of the horizontal arrows in the diagram
  \eqref{ref-1.6-10} the opposite actions appear for
  reasons that are mysterious to the authors.}
\end{Thm}
\emph{The convention that wavy arrows indicate actions will be used
throughout below.}

\medskip

\noindent The following corollary will be applied to C{\u{a}}ld{\u{a}}raru's
conjecture below.
\begin{Cor}
There is a
commutative diagram of precalculi:
\begin{equation}\label{ref-1.8-13}
  \xymatrix{\bigoplus_{m,n\geq 0} \mathrm H^m(X,\bigwedge^n \mathcal L)\ar@{~>}[d]\ar[rr]^-{\mathrm{HKR}\circ\iota_{\sqrt{\mathrm{td}(\mathcal L)}}} && \mathbb H^\bul(X,D_{\poly}^{\mathcal L}(X))\ar@{~>}[d]\\
    \bigoplus_{m,n\geq 0}\mathrm H^m(X,\bigwedge^n \mathcal L^*) && \mathbb H^\bul(X,C_\mathrm{poly}^\mathcal L(X))\ar[ll]^{(\sqrt{\mathrm{td}(\mathcal L)}\wedge-)\circ \mathrm{HKR}}.}
\end{equation}
with $\mathbb H^\bul(X,-)$ denoting the hypercohomology functor.
\end{Cor}
\begin{proof}
This follows by applying the functor $\mathbb H^\bul(X,-)$ to the
commutative diagram (\ref{ref-1.6-10}).
\end{proof}
If we consider only the Lie brackets and the Lie algebra actions then the
horizontal isomorphisms in the commutative diagram (\ref{ref-1.6-10})
are obtained from the horizontal arrows in 
diagram \eqref{ref-1.7-12} below, which is part of our second main result:
\begin{Thm}\label{ref-1.2-11} {}{Assume that $\mathbb{R}\subset k$.}
  Let $\mathcal L$ be a locally free
  Lie algebroid of rank $d$ over the ringed site $(X,\mathcal O)$.
  There exist sheaves of differential graded Lie algebras $(\mathfrak
  g_i^\mathcal L,\mathrm d_i,[\ ,\ ]_i)$ and sheaves of DG-Lie 
  modules $(\mathfrak m_i^\mathcal L,\mathrm b_i,\mathrm L_i)$ over
  them as well as $L_\infty$-quasi-isomorphisms $\mathfrak U_\mathcal
  L$ from $\mathfrak g_1^\mathcal L$ to $\mathfrak g_2^\mathcal L$ and
  $\mathfrak S_\mathcal L$ from $\mathfrak m_2^\mathcal L$ to
  $\mathfrak m_1^\mathcal L$, which fit into the following commutative
  diagram:
\begin{equation}\label{ref-1.7-12}
\xymatrix{T_\mathrm{poly}^\mathcal L(X) \ar@{^{(}->}[r]\ar@{~>}[d]_{\mathrm L} & \mathfrak g_1^\mathcal L \ar[r]^{\mathfrak U_\mathcal L}\ar@{~>}[d]_{\mathrm L_1}& \mathfrak g_2^\mathcal L\ar@{~>}[d]_{\mathrm L_2} & D_\mathrm{poly}(X)\ar@{_{(}->}[l]\ar@{~>}[d]_{\mathrm L}\\
\Omega^\mathcal L(X)\ar@{^{(}->}[r] & \mathfrak m_1^\mathcal L & \mathfrak m_2^\mathcal L\ar[l]_{\mathfrak S_\mathcal L} & C_\mathrm{poly}^\mathcal L(X)\ar@{_{(}->}[l]}.
\end{equation}
where the hooked arrows are strict (i.e.~DG-Lie) quasi-isomorphisms.
\end{Thm}
\subsubsection{Comments on the results and the proofs}

The proof of Theorems~\ref{ref-1.1-9} and \ref{ref-1.2-11} depends on
the simultaneous globalization of a number of local formality results
due to Kontsevich \cite{K} (see also \cite{MT}), Tsygan \cite{T},
Shoikhet \cite{Sh} and the first two authors \cite{CR, CR1,CR2}. This
globalization is performed by a functorial version of formal geometry
\cite{CVdB} (see also \cite{Y}).

The proof of Theorem~\ref{ref-1.1-9} roughly speaking involves the
construction of a morphism up to homotopy between the precalculus
structures up to homotopy on ($T_\mathrm{poly}^\mathcal L(X)$,
$\Omega^\mathcal L(X)$) and ($D_\mathrm{poly}^\mathcal
L(X),C_\poly^\mathcal L(X)$). In this paper we do not construct a
full ``precalculus$_\infty$''-quasi-isomorphism between these structures (in the case
that $\mathcal{L}$ is a tangent bundle this has been done in \cite{DTT2} using operadic
methods; actually, in {\em loc.\ cit.}, the authors work in the ``calculus$_\infty$'' setting, encoding also the De Rham differential, which is not part of the precalculus structure as observed before).  
On the other hand, in contrast to {\em loc.\ cit.}, the results we
prove are explicit and this fact is essential to recover C{\u{a}}ld{\u{a}}raru's
conjecture as formulated in \cite{Cald} (see Theorem \ref{ref-1.3-14} below).  

We are able to obtain such explicit results by starting with
the local quasi-isomorphisms of Kontsevich and Shoikhet which are given by explicit formul\ae\ (in contrast
to, say, Tamarkin's local $G_\infty$-quasi-isomorphism \cite{Tam}). 
While these are {\em a priori} only $L_\infty$-quasi-isomorphisms they
are nonetheless compatible with products up to homotopy \cite{K,CR2} in a
strong explicit sense and this turns out to be enough for our
purposes.

\medskip

  {}
  As the local quasi-isomorphisms of Kontsevich and Shoikhet are
  defined over~$\mathbb{R}$ (see \S\ref{ref-5.4-84}) we have to assume
  $\mathbb{R}\subset k$ in the statement of Theorem \ref{ref-1.2-11}.  However enough
  coefficients are rational (and computable), which in turn allows us to prove
  Theorem \ref{ref-1.1-9} over an arbitrary field of characteristic
  zero. This idea was already used in \cite{CVdB}. See
  \S\ref{sec-arbitrary}.  For Theorem \ref{ref-1.2-11} we could likely
  have started with a Tamarkin-style local quasi-isomorphism 
  \cite{Tam} defined over $\mathbb{Q}$, but
  since the coefficients of such a local quasi-isomorphism are not
  explicit, the result would not be immediately applicable to Theorem
  \ref{ref-1.1-9}.

\medskip

The existence of the upper horizontal isomorphism in
(\ref{ref-1.6-10}) has been proved independently in~\cite{DTT1,CVdB},
while its explicit form has been computed in~\cite{CVdB}.  The
existence of the lower horizontal isomorphism has been shown
in~\cite{DTT2}. As observed above, our approach {\em via} Kontsevich's and
Shoikhet's local formality formul\ae\ allows us to compute it explicitly.


\subsection{C{\u{a}}ld{\u{a}}raru's conjecture}
Assume now that $X$ is a smooth algebraic or complex variety.
C{\u{a}}ld{\u{a}}raru's conjecture (stated originally in the algebraic case)
asserts the existence of various compatibilities between the
Hochschild (co)homology and tangent (co)homology of $X$ (see
below). For the full statement we refer to \cite{Cald}. The results
in this paper complete the proof of C{\u{a}}ld{\u{a}}raru's conjecture.

\medskip

We now explain this in more detail.  The Hochschild (co)homology
\cite{Swan} of $X$ is defined as
\begin{align*}
\mathrm{HH}^n(X)&=\operatorname{Ext}^n_{\mathcal{O}_{X\times X}}(\mathcal{O}_\Delta,\mathcal{O}_\Delta) 
\qquad (n\geq0)\\
\mathrm{HH}_n(X)&=\operatorname{Tor}_{-n}^{\mathcal{O}_{X\times X}}(\mathcal{O}_\Delta,\mathcal{O}_\Delta)
\qquad (n\leq0)
\end{align*}
where $\Delta\subset X\times X$ is the diagonal. From these definitions
it is clear that  $\mathrm{HH}^\bullet(X)$ has a canonical algebra structure
(by the Yoneda product) and $\mathrm{HH}_\bullet(X)$ is a module over it. 

Furthermore if we put $\mathcal L=T_X$ then it is proved in
\cite{CRVdB2} (and partially in \cite{Y1}) that there are
isomorphisms of algebras and modules
\[
\xymatrix{
\mathrm{HH}^\bullet(X) \ar@{~>}[d]\ar[rr]  &&\mathbb H^\bul(X,D_{\poly}^{\mathcal L}(X))\ar@{~>}[d]\\
\mathrm{HH}_\bullet(X)\ar[rr]&&\mathbb H^\bul(X,C_\mathrm{poly}^\mathcal L(X))
}
\]
where on the right-hand side we consider only the part of the precalculus given by the cup and cap product. 

We define the tangent (co)homology of $X$ by
\[
\mathrm{HT}^\bul(X)=\bigoplus \mathrm H^\bul(X,\bigwedge^\bul T_X),\ \mathrm{H\Omega}_\bul(X)=\bigoplus \mathrm H^\bul(X,\Omega_X^{-\bul}).
\]
where now $\Omega_X^\bul$ denotes the
graded sheaf of differential forms on $X$.

The commutative diagram (\ref{ref-1.8-13}) then yields the following
\begin{Thm}[``C{\u{a}}ld{\u{a}}raru's conjecture'']\label{ref-1.3-14}
  For a smooth algebraic or complex variety $X$ over $k$ there is a
  commutative diagram of 
$k$-algebras and 
  modules
\begin{equation}\label{ref-1.9-15}
  \xymatrix{\mathrm{HT}^\bul(X)
\ar@{~>}[d]\ar[rrr]^-{\mathrm{HKR}\circ\iota_{\sqrt{\mathrm{td}(X)}}} &&& \mathrm{HH}^\bul(X)\ar@{~>}[d]\\
\mathrm{H\Omega}_\bul(X)     &&& \mathrm{HH}_\bul(X)\ar[lll]^{(\sqrt{\mathrm{td}(X)}\wedge -)\circ\mathrm{HKR}},}
\end{equation}
where $\mathrm{td}(X)$ is the Todd class for $\mathcal L=T_X$.
\end{Thm}
Theorem \ref{ref-1.3-14} completes the proof of the
parts of C{\u{a}}ld{\u{a}}raru's conjecture \cite{Cald} which do not depend on $X$
being proper. The cohomological part (the upper row in the above
diagram) had already been proved in \cite{CVdB} and is also an
unpublished result of Kontsevich. 

In the proper case there is an additional assertion in C{\u{a}}ld{\u{a}}raru's
conjecture which involves the natural bilinear form on
$\mathrm{HH}_{\bullet}(X)$. We do not consider this assertion in the
present paper as it has already been proved by Markarian~\cite{Mar} and
Ramadoss~\cite{Ram}.  If we combine Theorem \ref{ref-1.9-15} with the 
results of Markarian and Ramadoss we obtain a full proof of C{\u{a}}ld{\u{a}}raru's
conjecture. Let us also mention that in the compact Calabi-Yau case
C{\u{a}}ld{\u{a}}raru's conjecture has been proved in~\cite{H}.

\begin{Ack}
We express our gratitude to the anonymous referees for the careful reading of the manuscript.
Their comments have helped us improve the paper.
\end{Ack}

\section{Notation and conventions}
As stated already we always work over a ground field
$k$ of characteristic $0$; unadorned tensor products are over $k$.

Most objects we consider are equipped with a topology which will be
explicitly specified when needed. However if an object is introduced
without a specific topology, or if the topology is not clear from the
context, then it is assumed to be equipped with the discrete topology.

Many objects we will encounter are $\mathbb Z$-graded.
Koszul's sign rule is always assumed in this context. For a double or higher
complex we apply the sign rule with respect to total degree.

\section{Some recollections on Lie algebroids and related topics}\label{ref-3-16}

\subsection{Generalities on Lie algebroids}\label{ref-3.1-17}
In this section $R$ is a commutative $k$-algebra.
\begin{Def}\label{ref-3.1-18}
  A {\bf Lie algebroid} $L$ over $R$ is a Lie algebra over $k$ which is in addition an $R$-module and is
  endowed with an anchor map $\rho:L\to \mathrm{Der}_k(R)$
satisfying the compatibility
\begin{equation}\label{ref-3.1-19}
[l_1,r l_2]=\rho(l_1)(r)l_2+r[l_1,l_2],\ r\in R,\ l_i\in L,\ i=1,2.
\end{equation}
\end{Def}
The basic example of a Lie algebroid over $R$ is
$L=\mathrm{Der}_k(R)$ with the identity anchor map and the commutator Lie
bracket.

If $L$ is a Lie algebroid then $R\oplus L$ is a Lie algebra with Lie
bracket $[(r,l),(r',l')]=(\rho(l)(r')-\rho(l')(r),[l,l'])$.  We define
the \textbf{universal enveloping algebra} $\mathrm U_R(L)$ of~$L$ to
be the quotient of the enveloping algebra associated to the Lie
algebra $R\oplus L$ by the relation $r\otimes l=rl$ ($r\in R$, $l\in
R\oplus L$).

For the sake of simplicity, below we will usually omit the anchor map $\rho$
from the notation, unless it is necessary for the sake of clarity.

The universal enveloping algebra of a Lie algebroid satisfies a
universal property similar to that of an ordinary enveloping algebra.
This implies for example that the anchor map $\rho$ uniquely extends
to an algebra morphism from $\mathrm U_R(L)$ to $\mathrm{End}_k(R)$,
or equivalently: it yields a left $\mathrm U_R(L)$-module structure on
$R$.
\emph{For reasons which will become clear later we assume that our Lie
algebroids  are free of rank $d$ over $R$.}

\subsubsection{$L$-poly-vector fields and $L$-differential forms over $R$}\label{ref-3.1.1-20}
To a Lie algebroid $L$ over $R$ we associate 
\begin{align}
\label{ref-3.2-21} T_\mathrm{poly}^L(R)&=\bigoplus_{n\geq -1}\wedge^{n+1}_R L,\\
\label{ref-3.3-22} \Omega^L(R)&=\bigoplus_{n\leq 0} \wedge_R^{-n}L^*,\ L^*=\mathrm{Hom}_R(L,R).
\end{align}
We refer to (\ref{ref-3.2-21}) and (\ref{ref-3.3-22}) as the
spaces of {\bf $L$-poly-vector fields} and {\bf $L$-forms} on~$R$.

As an exterior algebra $T_\mathrm{poly}^L(R)$ has a wedge product
which we denote by $\cup$ (``the cup product'').  The extension of the
Lie bracket on $L$ to a bi-derivation on $T_\mathrm{poly}^L(R)$ defines
a Lie bracket which is called the Schouten--Nijenhuis bracket and is
denoted by $[-,-]$. Note that with our grading conventions the
cup product has degree one and the Lie bracket has degree zero. The
cup product and the Lie bracket make $T_\mathrm{poly}^L(R)$ into a
(shifted) Gerstenhaber algebra with trivial differential.

On the other hand, $\Omega^L(R)$ is obviously a graded algebra with respect to
the wedge product. In addition there is an analogue $\mathrm d_L$ of the De Rham
differential on $\Omega^L(R)$, which is given on generators by
\begin{align*}
\mathrm d_L(r)(l)&=l(r)\\
\mathrm d_L(l^*)(l_1,l_2)&=l_1(l^*(l_2))-l_2(l^*(l_1))-l^*([l_1,l_2]),
\end{align*}
for $r\in R$, $l,l_i\in L$, $i=1,2$, $l^*\in L^*$, and is extended uniquely by Leibniz's rule.

The natural contraction operation of $L$-forms on $R$ with respect to
$L$-poly-vector fields is denoted by $\cap$ (the ``cap product'').  The
Lie derivative $\mathrm L$ of $L$-forms on $R$ with respect to $L$-poly-vector
fields is specified in the usual way {\em via} Cartan's homotopy formula as the
commutator of $\mathrm d_L$ and the contraction.  
The pair 
$$
((T_\mathrm{poly}^L(R),[-,-],\cup),(\Omega^L(R),\cap,\mathrm L))
$$
forms a precalculus (see \S\ref{ref-1.2-5}). 
\subsubsection{ $L$-connections}\label{ref-3.1.2-23}
As usual  $L$ is a Lie algebroid over $R$.
\begin{Def}\label{ref-3.2-24}
  Let $M$ be an $R$-module $M$. An {\bf $L$-connection} on $M$ is a
  $k$-linear map $\nabla$ from $M$ to $L^*\otimes_R M$, which satisfies
  Leibniz's rule
\begin{equation}\label{ref-3.4-25}
\nabla (rm)=\mathrm d_L(r)\otimes m+r\nabla m,\ r\in R,\ m\in M.
\end{equation}
\end{Def}
The $L$-connection $\nabla$  is said to be {\bf flat}, if
$\nabla^2=0$. Equivalently, the assignment $l\mapsto
\nabla_l$, where $\nabla_l$ denotes the action of $\nabla$
followed by contraction with respect to $l$, defines a Lie algebra morphism
from $L$ to $\mathrm{End}_k(M)$.

If we let $l\in L$ act as $\nabla_l$ then a flat $L$-connection on $M$
extends to a left $\mathrm
U_R(L)$-module structure on $M$.

Furthermore a flat $L$-connection $\nabla$ on $M$ can be extended to a
differential (denoted by the same symbol) on the graded $R$-module
$\Omega^L(R)\otimes_R M$ {\em via} Leibniz's rule
\[
\nabla (\omega\otimes_R m)=\mathrm d_L\omega\otimes_R m+(-1)^{|\omega|}\omega\wedge \nabla m,\ \omega\in \Omega^L(R),\ m\in M.
\]
\subsubsection{$L$-differential operators over $R$}\label{ref-3.1.3-26}
In this section we define the algebra of poly-differential operators
of a Lie algebroid and we list some of its properties. We give
some explicit formul\ae\ along the lines of~\cite{CR2}.

As in the case of ordinary Lie algebras $\mathrm U_R(L)$ (see
\S\ref{ref-3.1-17}) may be naturally filtered by giving $R$ filtered degree 0 and
$L$ filtered degree 1. In particular
\[
F_0\mathrm U_R(L)=R,\ F_1\mathrm U_R(L)=R\oplus L,
\]
We view $\mathrm U_R(L)$ as an $R$-\emph{central} bimodule {\em via} the natural embedding of $R$ into $\mathrm
U_R(L)$.  Explicitly, if we denote this embedding by $i$
then\footnote{Note that there is an at first sight more natural right $R$-module
  structure on $U_R(L)$ given by the formula $Dr=Di(r)$. This
  alternative right module structure will not be used in this paper.}
\begin{equation}
\label{ref-3.5-27}
rD\overset{!}{=}Dr\overset{\text{def}}{=}i(r)D,\ r\in R, D\in \mathrm U_R(L)
\end{equation}
Moreover 
$\mathrm U_R(L)$ is an $R$-coalgebra \cite{Xu}, i.e.\ $\mathrm U_R(L)$
possesses an $R$-linear coproduct $\Delta:\mathrm U_R(L)\to \mathrm
U_R(L)\otimes_R\mathrm U_R(L)$ and an $R$-linear counit, satisfying
the usual axioms. The comultiplication actually takes values in
\begin{multline*}
(\mathrm U_R(L)\otimes_R\mathrm U_R(L))'\\
=\left\{\sum_j D_j\otimes E_j\in U_R(L)\otimes_R\mathrm U_R(L)\mid \forall r\in R:
\sum_j D_ji(r)\otimes E_j=\sum_j D_j\otimes E_ji(r)\right\}
\end{multline*}
which is an $R$-algebra even though $\mathrm U_R(L)\otimes_R\mathrm \mathrm U_R(L)$  is not. 

 The comultiplication $\Delta$ and counit $\epsilon$
are given by similar formul\ae\ as in the Lie algebra case
\begin{equation}\label{ref-3.6-28}
\begin{aligned}
\Delta(r)&=r\otimes_R 1=1\otimes_R r&& r\in R\\
\Delta(l)&=l\otimes_R 1+1\otimes_R l&&l\in L\\
\Delta(DE)&=D_{(1)}E_{(1)}\otimes D_{(2)}E_{(2)} && D,E\in \mathrm U_R(L)\\
\epsilon(D)&=D(1)
\end{aligned}
\end{equation}
In the third formula we have used Sweedler's convention. 
 The expression on the right-hand side is well defined because it is the product inside the algebra
$(\mathrm U_R(L)\otimes_R\mathrm U_R(L))'$.
 In the fourth formula we have used the natural
action of $\mathrm U_R(L)$ on $R$ (see \S\ref{ref-3.1-17}).

The algebra (better: in the terminology of~\cite{Xu,Cal} ``the Hopf
algebroid'') $\mathrm U_R(L)$ may be thought of as an algebra of
$L$-differential operators on $R$: in the case $L=\mathrm{Der}_k(R)$ and
$R$ smooth over $k$ then
$\mathrm U_R(L)$ coincides with the algebra of differential operators
on $R$.
\subsubsection{$L$-jets}\label{ref-3.1.4-29}
Let $(\mathrm U_R(L))_{\le n}$ be the elements of degree $\le n$ with respect to the canonical filtration
on $\mathrm U_R(L)$ introduced in \S\ref{ref-3.1.3-26}. The $L$-$n$-jets
are defined as 
\[
J^nL=\Hom_R(\mathrm U_R(L)_{\le n},R)
\]
(this is unambiguous, as the left and right $R$-modules structures on
$\mathrm U_R(L)$ are the same, see \eqref{ref-3.5-27}). 
We also put 
\begin{equation}\label{ref-3.7-30}
J\!L=\Hom_R(\mathrm U_R(L),R)=\projlim_n J^nL\qquad(\textrm{as }\mathrm U_R(L)=\injlim_n (\mathrm U_R(L))_{\le n}).
\end{equation}
$J\!L$ has a natural commutative
algebra structure obtained from the comultiplication on $\mathrm U_R(L)$. Thus for
$\phi_1,\phi_2\in J\!L$, $D\in \mathrm U_R(L)$ we have
\[
(\phi_1\phi_2)(D)=\phi_1(D_{(1)})\phi_2(D_{(2)})\,,
\]
and the unit in $J\!L$ is given by the counit on $\mathrm U_R(L)$. 

In addition $J\!L$ has two commuting left $\mathrm U_R(L)$-module structures
which we now elucidate. 
First of all there are two distinct monomorphisms of
$k$-algebras
\begin{align*}
\alpha_1:&R\rightarrow J\!L:r\mapsto (D\mapsto r\epsilon(D))\,,\\
\alpha_2:&R\rightarrow J\!L:r\mapsto (D\mapsto D(r))\,.
\end{align*}
It will be convenient to write $R_i=\alpha_i(R)$
and to view $J\!L$ as an $R_1-R_2$-bimodule. 

There are also two distinct commuting actions by derivations
of $L$ on $J\!L$. Let $l\in L$, $\phi\in J\!L$, $D\in \mathrm U_R(L)$.
\begin{align*}
{}^1\nabla_l(\phi)(D)&=l(\phi(D))-\phi(lD)\\
{}^2\nabla_l(\phi)(D)&=\phi(Dl)
\end{align*}
Again it will be convenient to write $L_i$ for $L$ acting by
${}^i\nabla$. Then ${}^i\nabla$ defines a flat $L_i$-connection on
$J\!L$, considered as an $R_i$-module. The connection ${}^1\nabla$ is
the well-known {\bf Grothendieck connection}. It follows that $J\!L$ is a
$\mathrm U_R(L)_1-\mathrm U_R(L)_2$-bimodule (with both $\mathrm
U_R(L)_1$ and $\mathrm U_R(L)_2$ acting on the left).

The $\mathrm U_R(L)_2$ action on $J\!L$ takes the very simple form
\[
(D\cdot \phi)(E)= \phi(ED)
\]
(for $D,E\in \mathrm U_R(L)_2$, $\phi \in J\!L$). 

Define $\epsilon:J\!L\rightarrow R$ by $\epsilon(\phi)=\phi(1)$ and put
$J^cL=\ker \epsilon$. Then $J\!L$ is complete for the $J^cL$-adic topology
and the filtration on $J\!L$ induced by \eqref{ref-3.7-30} coincides with the
$J^cL$-adic filtration.  If we filter $J\!L$ with the $J^cL$-adic
filtration then we obtain
\begin{equation}
\label{ref-3.8-31}
\operatorname{gr} J\!L=\textsc{S}_{R} L^\ast
\end{equation}
and the $R_1$ and $R_2$-action on the r.h.s. of this equation coincide (here
 and below the letter $\textsc{S}$ stands for ``symmetric algebra'').

The induced actions on $\operatorname{gr} J\!L=\textsc{S}_R L^\ast$ of $l\in L$, considered as an
element of $L_1$ and $L_2$, are given by the contractions $i_{-l}$ and
$i_l$, respectively.

\medskip

In case $R$ is the coordinate ring of a smooth affine algebraic
variety and $L=\operatorname{Der}_k(R)$ then we may identify $J\!L$ with
the completion $R\widehat\otimes R$ of $R\otimes R$ at the kernel of the multiplication map
$R\otimes R\rightarrow R$. The two actions of $R$ on $J\!L$ are
respectively $R\widehat\otimes 1$ and $1\widehat\otimes R$.

Similarly a derivation on $R$ can be extended to $R\widehat\otimes R$ in two
ways by letting it act respectively on the first and second factor. Since
derivations are continuous they act on adic completions and hence
in particular on $J\!L$.  This provides the two actions of $L$ on $J\!L$. 

\medskip

\emph{In the sequel we will view the action labelled by ``\/$1$\!'' as the
default action. I.e.\ we will usually not write the $1$\! explicitly.}

\subsection{Relative poly-vector fields, poly-differential operators}\label{ref-3.2-32}
We need relative poly-differential operators and poly-vector fields. So
assume that $A\rightarrow B$ is a morphism of commutative $k$-algebras. Then
\begin{align*}
T_{\poly,A}(B)&=\bigoplus_{n\ge -1} T^n_{\poly,A}(B)\\
D_{\poly,A}(B)&=\bigoplus_{n\ge -1} D^n_{\poly,A}(B)
\end{align*}
where $T^n_{\poly,A}(B)=\bigwedge^{n+1}_B\operatorname{Der}_A(B)$.  
Similarly, $D^n_{\poly,A}(B)\subseteq \mathrm{Hom}_A(B^{\otimes_A(n+1)},B)$ consists of those $A$-linear maps from $B^{\otimes_A(n+1)}$ to $B$, which are $A$-linear differential operators on $B$ in each argument.

It is easy to see that $T_{\poly,A}(B)$ is a Gerstenhaber algebra
when equipped with the Schouten bracket and the exterior product. 
Similarly $D_{\poly,A}(B)$ is a graded subspace of the relative Hochschild
complex $C^\bullet_A(B)$ and since differential operators are closed
under composition one easily sees that it is in fact a sub-$B_\infty$-algebra, see Appendix~\ref{ref-A-127} for more details on $B_\infty$-algebras.

If $A$ and $B$ are DG-algebras then we equip
$T_{\poly,A}(B)$, $D_{\poly,A}(B)$ with the total
differentials $[\mathrm d_B,-]$ and $[\mathrm d_B,-]+\mathrm d_H$ where $\mathrm d_H$ denotes the Hochschild
differential. Similar results now apply.

\subsection{The sheaf of $L$-poly-differential operators}

\begin{Def}\label{ref-3.3-33}
  For a Lie algebroid $L$ over $R$ we define the graded vector space
  $D_\mathrm{poly}^L(R)$ of {\bf $L$-poly-differential operators on
    $R$} as the tensor algebra over $R$ of $\mathrm U_R(L)$ with
  shifted degree, i.e.\
\[
D_\mathrm{poly}^L(R)=\bigoplus_{n\geq -1}\mathrm U_R(L)^{\otimes_R
  (n+1)}.
\]
\end{Def}
The action of $\mathrm U_R(L)$ on $R$ extends to a map 
\begin{equation}\label{ref-3.9-34}
D_\mathrm{poly}^{L,n}(R)\rightarrow \operatorname{Hom}_k(R^{\otimes n+1}, R)
\end{equation}
defined by
\[
(D_1\otimes\cdots \otimes D_{n+1})(r_1\otimes\cdots\otimes r_{n+1})
\mapsto D_1(r_1)\cdots D_{n+1}(r_{n+1})
\]
whose image lies in the space $D_\poly(R)$ of poly-differential
operators on $R$.

$D_{\poly}^L(R)$ is a $B_\infty$-algebra. In particular it is 
a DG-Lie algebra and furthermore it is a Gerstenhaber
algebra up to homotopy. In Appendix \ref{ref-A-127} we give the formul\ae\
for the full $B_\infty$-structure. Here we content ourselves by reminding
the reader of the basic operations.

\medskip

The \textbf{Gerstenhaber bracket} on
$D_\mathrm{poly}^L(R)$ is defined by
\begin{equation}\label{ref-3.10-35}
[D_1,D_2]=D_1\{D_2\}-(-1)^{|D_1||D_2|}D_2\{D_1\},\ D_i\in D_\mathrm{poly}^L(R),\ i=1,2.
\end{equation}
where 
\[
D_1\{D_2\}=\sum_{i=0}^{|D_1|} 
(-1)^{i|D_2|}(\operatorname{id}^{\otimes i}\otimes \Delta^{|D_2|}
\otimes \operatorname{id}^{\otimes |D_1|-i})(D_1)
\cdot (1^{\otimes i}\otimes D_2\otimes 1^{\otimes |D_1|-i})\,.
\]
It is a Lie bracket  of degree $0$.
The special element $\mu=1\otimes_R 1\in
D^{L,1}_\mathrm{poly}(R)=\mathrm U_R(L) \otimes_R \mathrm U_R(L)$ satisfies
$[\mu,\mu]=0$. The {\bf Hochschild differential} is defined as
 the operator $\mathrm d_\mathrm H=[\mu,-]$.

\medskip

The \textbf{cup product} on  $D_\mathrm{poly}^L(R)$ is defined by
\begin{equation}
\label{ref-3.11-36}
D_1\cup D_2=(-1)^{(|D_1|-1)(|D_2|-1)}D_1\otimes_R D_2
\end{equation}
 (See also Appendix~\ref{ref-A-127} for an explicit derivation of the previous Formula). 
 
\medskip

One may now show that these operations make the $4$-tuple
$(D_\mathrm{poly}^L(R),\mathrm d_\mathrm H,[\ ,\ ],\cup)$ into a
Gerstenhaber algebra up to homotopy (see Lemma~\ref{ref-A.1-130}).
Indeed if $R$ is smooth over $k$ and $L=\operatorname{Der}_k(R)$ is the tangent bundle then
the operations we have defined are the same  as those one obtains
from the identification $D^L_{\poly}(R)=D_{\poly}(R)$ where we view the
right-hand side as a sub-$B_\infty$-algebra of the Hochschild complex $C^\bullet(R)$ of $R$
(cfr. \S\ref{ref-3.2-32}). 

\medskip

It is in fact, as we explain now, not necessary to verify that
we have defined a homotopy Gerstenhaber structure on
$D^L_{\poly}(R)$. Indeed the results can be obtained directly from the
known results for the Hochschild complex (see
\cite{GJ,GV}). Similarly it is not
necessary to write explicit formul\ae\ for $[-,-]$ and $\cup$ (or for
the whole $B_\infty$-structure for that matter). This point of view
will be useful when we consider Hochschild chains as in that case the formul\ae\ become more complicated.

\medskip 

The $L_2$-action on $J\!L$ commutes with the $R_1$-action (see \S\ref{ref-3.1.4-29})
so we obtain a ring homomorphism
\[
\mathrm U_{R_2}(L_2)\rightarrow D_{R_1}(J\!L):D\mapsto (\theta\mapsto D(\theta))\,.
\]
and hence a map
\begin{equation}
\label{ref-3.12-37}
D^{L_2}_{\poly}(R_2)\rightarrow D_{\poly, R_1}(J\!L)\,.
\end{equation}
of Gerstenhaber algebras up to homotopy. 
The right-hand side has an $R_1$-connection given by $[{}^1\nabla,-]$ and it
follows from \cite[Prop.\ 4.2.4, Lemma 4.3.4]{CVdB} that the left-hand side
of \eqref{ref-3.12-37} is given by the horizontal sections for this connection. 

Now as discussed in \S\ref{ref-3.2-32}, we know that $D_{\poly,
  R_1}(J\!L)$ is a $B_\infty$-algebra and it is an easy verification
that the braces and the differential, which make up the $B_\infty$-structure, are
horizontal for $[{}^1\nabla,-]$. Hence the $B_\infty$-structure on
$D_{\poly, R_1}(J\!L)$ descends to $D^{L}_{\poly}(R)$ and one
verifies that its basic operations are indeed given by the formul\ae\ we gave earlier. 
\subsection{The Hochschild complex of $L$-chains over $R$}\label{ref-3.4-38}
We start with the following definition.
\begin{Def}\label{ref-3.4-39}
For a Lie algebroid $L$ over $R$, the graded $R$-module 
\begin{equation}\label{ref-3.13-40}
C_{\mathrm{poly},p}^L(R)=\begin{cases}
J\!L^{\widehat{\otimes}_R -p},& p<0\\
R,& p=0,
\end{cases}
\end{equation}
is called the space of {\bf Hochschild  $L$-chains over $R$}.
\end{Def}
Our aim in this section will be to show that the pair
\[
(D_{\poly}^L(R),C_{\mathrm{poly}}^L(R))
\]
is a precalculus up to homotopy. We will do this without relying on
explicit formul\ae\ (as they are quite complicated). Instead we will
reduce to a relative version of \cite{CR2} which discusses Hochschild
(co)homology. Explicit formul\ae\ are given in Appendix~\ref{ref-B-136}.

\medskip

Let us first remind the reader that if $A$ is a $k$-algebra then the
pair $(C^{\bullet}(A),C_{\bullet}(A))$ consisting of the spaces of
Hochschild cochains and chains is a precalculus up to homotopy. For
$C^{\bullet}(A)$ this is just the (shifted) homotopy Gerstenhaber structure
which we have already mentioned in \S\ref{ref-3.3-33} and which was introduced in
\cite{GJ,GV}.

The full precalculus structure up to homotopy on
$(C^{\bullet}(A),C_{\bullet}(A))$ is a more intricate object. A
complete treatment in a very general setting has been given in
\cite{CR2}.  It is shown that the precalculus structure can
be obtained from two interacting $B_\infty$-module structures on
$C_{\bullet}(A)$. These $B_\infty$-module structures are obtained from
brace-type operations.  For more operadic approaches see
\cite{DTT2}.

Although we do not really use them, for the benefit of the reader we state the well-known formul\ae\ for the
contraction, the Lie derivative and the differential. If $P\in C^{m-1}(A)=\Hom(A^{\otimes m},A)$ and $(a_0|\cdots|a_t)\in C_{-t}(A)=A^{\otimes t+1}$ then we
have
\begin{align*}
  \imath_P (a_0|\cdots|a_t)&=(a_0P(a_1,\ldots,a_m)|a_{m+1}|\cdots|a_t)\\
\mathrm  L_P (a_0|\ldots|a_t)&=\sum_{i=0}^{t-m+1} (-1)^{(m-1)i} (a_0|\cdots|
  a_{i-1} |P(a_i,\ldots, a_{i+m-1})|a_{i+m}|\cdots| a_t)
  \\
  &\qquad +\sum_{l=t-m+2}^{t+1}
  (-1)^{lt}(P(a_l,\ldots,a_t,a_0,\ldots,a_{m-t+l-2})|a_{m-t+l-1}|\cdots| a_{l-1})
\end{align*}
The differential $\mathrm b_H$ is defined as $\mathrm L_\mu$ where $\mu$ is the multiplication,
considered as an element of $\Hom(A^{\otimes 2}, A)$. 

\medskip

To construct the precalculus structure up to homotopy on
$
(D_{\poly}^L(R),C_{\mathrm{poly}}^L(R))
$
we proceed as in \S\ref{ref-3.3-33}. We first define an object that is
larger than $C_{\mathrm{poly}}^L(R)$. 
\begin{Def}\label{ref-3.5-41}
  The space of {\bf $L$-poly-jets over $R$} is the \emph{completed}
  space of relative Hochschild chains
  $\widehat{C}_{R_1,\bullet}(J\!L)$. Explicitly
\begin{equation}\label{ref-3.14-42}
\widehat{C}_{R_1,\bullet}(J\!L)=\bigoplus_{p\leq 0} J\!L^{\widehat{\otimes}_{R_1} -p-1}
\end{equation}
\end{Def}

\medskip

The Grothendieck connection ${}^1\nabla$ on $J\!L$ (see \S\ref{ref-3.1.4-29})
yields a connection on $\widehat{C}_{R_1,\bullet}(J\!L)$
 by Leibniz's rule to
which we also refer to as the Grothendieck connection. The following
result was proved in \cite{CDH}.
\begin{Prop}\label{ref-3.6-43} 
  For a Lie algebroid $L$ over a commutative ring $R$ as above, there is
an isomorphism of graded vector spaces
\begin{equation}
\label{ref-3.15-44}
\widehat{C}_{R_1,\bullet}(J\!L)^{{}^1\nabla}\rightarrow C_{\mathrm{poly}}^L(R)
\end{equation}
which sends
\[
\phi_1\otimes\phi_2\otimes \cdots\otimes \phi_p\in \widehat{C}_{R_1,-p}(J\!L)^{{}^1\nabla}
\]
to
\[
\epsilon(\phi_1)\phi_2\otimes \cdots\otimes \phi_p \in C_{\mathrm{poly},1-p}^L(R)
\]
\end{Prop}
\begin{proof}
The arguments of the proof of \cite[Prop.\ 1.11]{CDH}, can be repeated almost {\em verbatim}.
\end{proof} 
The formul\ae\ from \cite{CR2} for the Hochschild complexes now yield
that
\[
(D_{\poly,R_1}(J\!L),\widehat{C}_{R_1,\bullet}(J\!L))
\subset (C^{\text{cont},\bullet}_{R_1}(J\!L),\widehat{C}_{R_1,\bullet}(J\!L))
\]
is a precalculus up to homotopy.  Furthermore one verifies that the
formul\ae\ in \cite{CR2} are compatible with the Grothendieck connection
${}^1\nabla$. Hence the precalculus descends to one on 
\begin{equation}
\label{ref-3.16-45}
(D_{\poly,R_1}(J\!L)^{{}^1\nabla},\widehat{C}_{R_1,\bullet}(J\!L)^{{}^1\nabla})
=(D_{\poly}^L(R),C^L_{\poly}(R))
\end{equation}
where we use \eqref{ref-3.12-37} as well as Proposition \ref{ref-3.6-43}. 

\medskip

It remains to check that this construction coincides with the standard one for a smooth commutative algebra. Namely, if $R/k$ is smooth and $L$
is the tangent bundle then we have
\[
D^L_{\poly}(R)=D_{\poly}(R)
\]
We also have $J\!L=R\widehat{\otimes} R$ (see \S\ref{ref-3.1.4-29}) and in this
way we obtain an isomorphism
\begin{equation}
\label{ref-3.17-46}
C^L_{\poly,-p}(R)=(R\widehat\otimes R)^{\widehat{\otimes}_{R_1} p}\rightarrow
R^{\widehat{\otimes} p+1}:(r_1\widehat\otimes s_1)\widehat\otimes\cdots \widehat\otimes (r_p\widehat\otimes s_p)
\mapsto (r_1\cdots r_p)\widehat\otimes s_1\widehat\otimes\cdots\widehat\otimes s_p
\end{equation}
which yields an isomorphism of graded vector spaces
\[
C^L_{\poly}(R)=\widehat{C}_{\bullet}(R)
\]
Thus we have an isomorphism of pairs of graded vector spaces
\begin{equation}
\label{ref-3.18-47}
(D^L_{\poly}(R),C^L_{\poly}(R))=(D_{\poly}(R), \widehat{C}_{\bullet}(R))
\end{equation}
The right-hand side is a precalculus up to homotopy (as it is basically
a pair of spaces of Hochschild chains/cochains). 
\begin{Lem}
The  precalculus up to homotopy on the right-hand side of \eqref{ref-3.18-47}
is the same one as the one we have constructed on the left-hand side. 
\end{Lem}
\begin{proof}
Note that going from the pair $(k,R)$ to $(R,J\!L)$ is a base
extension by $R$ (since $J\!L=R\widehat{\otimes} R$). Since the formul\ae\ 
in \cite{CR2} are clearly compatible with base extension we have
that the precalculus structure on
\begin{equation}
\label{ref-3.19-48}
(D_{\poly,R}(J\!L),\widehat{C}_{R,\bullet}(J\!L))=(R\widehat{\otimes} D_{\poly}(R),
R\widehat{\otimes} \widehat C_{\bullet}(R))
\end{equation}
is obtained by base extension from the one on 
\[
( D_{\poly}(R), \widehat C_{\bullet}(R))
\]
Furthermore one checks that the Grothendieck connections on
$D_{\poly,R}(J\!L)$ and $\widehat{C}_{R,\bullet}(J\!L)$  under the isomorphism
\eqref{ref-3.19-48} act by the standard Grothendieck connection on the copy of $R$ appearing on the left
of $\widehat{\otimes}$ and trivially on $D_{\poly}(R)$, $ \widehat C_{\bullet}(R)$. Hence its
invariants are precisely $D_{\poly}(R)$, $\widehat{C}_{\bullet}(R)$. This finishes
the proof.
\end{proof}

\subsection{The Hochschild--Kostant--Rosenberg Theorem in the Lie algebroid framework}\label{ref-3.5-49}

We recall the Lie algebroid version of the famous cohomological Hochschild--Kostant--Rosenberg (shortly, HKR) quasi-isomorphism; for a proof, we refer to~\cite{Cal}.
\begin{Thm}\label{ref-3.8-50}
We consider a Lie algebroid $L$ over $R$ in the sense of Definition~\ref{ref-3.1-18}, which is assumed to be free of rank $d$ over $R$.

Then, the map 
\begin{equation}\label{ref-3.20-51}
\mathrm{HKR}(l_1\wedge\cdots\wedge l_p)=(-1)^{\frac{p(p-1)}2}\frac{1}{p!}\sum_{\sigma\in\mathfrak S_p}(-1)^\sigma l_{\sigma(1)}\otimes_R\cdots\otimes_R l_{\sigma(p)}
\end{equation}
defines a quasi-isomorphism of complexes from $(T_\mathrm{poly}^L(R),0)$ to $(D_\mathrm{poly}^L(R),\mathrm d_\mathrm H)$.
\end{Thm}
There is a dual version of Theorem~\ref{ref-3.8-50}, which will also be needed.
\begin{Thm}\label{ref-3.9-52}
The quasi-isomorphism~\eqref{ref-3.20-51} induces the quasi-isomorphism
\[
\mathrm{HKR}(a)=a\circ \mathrm{HKR}
\]
of complexes from $(C_\mathrm{poly}^L(R),\mathrm b_\mathrm H)$ to $(\Omega^L(R),0)$.
\end{Thm}

\section{Fedosov resolutions in the Lie algebroid framework}\label{ref-4-53}

\subsection{Introduction}
\label{sec-intro-formal}
{}{The aim of this section is to discuss  Fedosov
resolutions~\cite{Dol} in the Lie algebroid framework. These are
needed to formulate and prove the globalization result, which in turn
leads to the main results.

To help the reader understand our algebraic setup (which was inspired
by \cite{Y}) we give some motivation for the definitions in the
subsequent sections.  
For the sake of exposition
we assume in this introduction that $X$ is some kind of $d$-dimensional smooth space
and $\mathcal{L}$ is an appropriate version of the tangent bundle of
$X$.   

One of the applications of formal geometry is the
globalization of local \emph{coordinate dependent} constructions. For
example using the Darboux Lemma it is trivial to quantize a symplectic
manifold locally but such local quantizations are coordinate dependent
and they do not globalize easily.
The same is true for the local
formality morphisms (see \S\ref{ref-5.4-84} below for more details) which we
use in this paper.

The idea is then to replace $X$ by a much larger infinite dimensional space
$X^{\coord}\rightarrow X$ that parametri\-zes  formal
local coordinate systems on $X$.  For example if $X$ is an algebraic
variety then the fiber at $x\in X$ in $X^{\coord}$ is given by the
$k$-algebra isomorphisms $
\widehat{\mathcal{O}}_{X,x}\rightarrow k[[t_1,\ldots,t_d]] 
$.  An equivalent way of saying this is that
$X^{\coord}$ universally trivializes the jet bundle
$(\widehat{\mathcal{O}}_{X,x})_{x\in X}$ over $X$.

Local constructions can be tautologically globalized to $X^{\coord}$
and this should be followed by some type of descent for
$X^{\coord}/X$. A general procedure to do this is to resolve $\mathcal{O}_X$ by a De Rham-type complex over $\mathcal{O}_{X^{\coord}}$ but
this does not really work as the fibers of $X^{\coord}\rightarrow X$
are not contractible.

However in the aforementioned examples the local constructions are all compatible with \emph{linear} coordinate changes. 
So if we define $X^{\aff}=X^{\coord}/\operatorname{Gl}_d$ then the constructions
descend to $X^{\aff}$ and as the fibers of $X^{\aff}/X$ are contractible
we can descend further to $X$. 

In this paper we work over a general locally free Lie algebroid $\mathcal{L}$
rather than~$T_X$. In this setting we define the analogue of  $X^{\coord}$
as the space which universally trivializes the space
of jet bundles for $\mathcal{L}$ (see \S\ref{ref-3.1.4-29}). 
}
\subsection{Setup}
As a general principle we work on the presheaf level in this paper, performing
 sheafification only as the very last step of the constructions.
 This
means that we may throughout replace all spaces by rings and locally free
sheaves may be treated as free modules.



As before we consider a Lie algebroid $L$ over a ring $R$ in the sense
of Definition~\ref{ref-3.1-18}, {\em i.e.} $L$ is free
of rank $d$ over $R$.

First we
 discuss Fedosov resolutions of $L$-poly-vector fields and
$L$-poly-differential operators as Gerstenhaber algebras up to
homotopy, referring to~\cite{CVdB} for details.  Finally, we
discuss Fedosov resolutions of $\Omega^L(R)$ (see~\eqref{ref-3.3-22}) and $C_\mathrm{poly}^L(R)$
(see~\eqref{ref-3.13-40}) which are compatible with the precalculus structure up to homotopy.

\subsection{The (affine) coordinate space of a Lie algebroid}\label{ref-4.1-54}

For a Lie algebroid $L$ over $R$ as above, its {\bf
  coordinate space $R^{\mathrm{coord},L}$} has been introduced and
discussed in details in~\cite{VdB,CVdB}, to which we refer for a more
extensive treatment.

{}{As explained in \S\ref{sec-intro-formal}, the main property of $R^{\mathrm{coord},L}$ is} the existence of an isomorphism of
$R^{\mathrm{coord},L}$-algebras
\begin{equation}\label{ref-4.1-55}
t:R^{\mathrm{coord},L}\widehat{\otimes}_{R_1} J\!L\to R^{\mathrm{coord},L}[\![x_1,\dots,x_d]\!]=R^{\mathrm{coord},L}\widehat{\otimes}F,\ F=k[\![x_1,\dots,x_d]\!],
\end{equation} 
and $R^{\mathrm{coord},L}$ is universal with respect to this property, {\em i.e.} if there is an $R$-algebra $W$, such that there is a $W$-linear isomorphism $W\widehat\otimes_{R_1}J\!L\to W[\![x_1,\dots,x_d]\!]$, then there exists a unique morphism $R^{\mathrm{coord},L}\to W$.

We note in particular that in contrast to $J\!L$ the ring $R^{\mathrm{coord},L}$ is \emph{not} an adic
topological ring: it is equipped with the discrete topology (like~$R$).
\begin{Exa}
\label{exa-poly}
Assume $R=k[x_1,\dots,x_d]$ and $L=\mathrm{Der}_k(R)$.
As explained in \cite[\S6.1.5]{VdB}\cite{Y}
we have 
\[
R^{\mathrm{coord},L}=R[y_{i,\alpha}:\ i=1,\dots,d,\ \alpha\in \mathbb N^d\smallsetminus\{0\}]_{\det(y_{i,e_j})},
\]
where $e_j$ is the $j$-th standard basis vector in $\mathbb{Z}^d$, and
the subscript $\det(y_{i,e_j})$ refers to the localization at the
indicated element.  
As
in this case $X=\operatorname{Spec} R$ has global coordinates
$x_1,\ldots,x_n$ the coordinate ring of the jet bundle $J\!L$
is equal to $R[\![y_1,\dots,y_d]\!]$ where $y_i$ is a local version of
the global coordinate $x_i$. The morphism~$t$ is the ``universal
Taylor expansion'' morphism
\[
t(y_i)=\sum_\alpha y_{i,\alpha}t^\alpha.
\]
\end{Exa}

As a consequence of the universal property of $R^{\mathrm{coord},L}$,
$R^{\mathrm{coord},L}$ admits an action of $\mathrm{GL}_d(k)$, such
that the following identity holds true on $R^{\mathrm{coord},L}\widehat{\otimes} F$
for $A\in \mathrm{GL}_d(k)$
\[
(A^{-1}\widehat{\otimes} A){\mid}_{J\!L}=\operatorname{Id}_{J\!L}
\]
where $J\!L$ is considered as a subalgebra of $R^{\mathrm{coord},L}\widehat{\otimes}F$ through \eqref{ref-4.1-55}.

By means of $R^{\mathrm{coord},L}$, we consider the graded algebra
$C^{\mathrm{coord},L}=\Omega_{R^{\mathrm{coord},L}}\otimes_{\Omega_{R_1}}\Omega^{L_1}(R_1)$.
It has the structure of a DG-algebra with
differential $\mathrm d_{C^{\mathrm{coord},L}}=\mathrm
d_{\Omega_{R^{\mathrm{coord},L}}}\otimes_{\Omega_{R_1}}
1+1\otimes_{\Omega_{R_1}} \mathrm d_{L_1}$, and inherits from
$R^{\mathrm{coord},L}$ a rational $\mathrm{GL}_d(k)$-action.

The universal isomorphism~\eqref{ref-4.1-55} extends to an isomorphism
\begin{equation}\label{ref-4.2-56}
t:C^{\mathrm{coord},L}\widehat{\otimes}_{R_1} J\!L\to C^{\mathrm{coord},L}[\![x_1,\dots,x_d]\!],
\end{equation}
where we used the respective obvious identifications
\[
\begin{aligned}
C^{\mathrm{coord},L}\widehat{\otimes}_{R^{\mathrm{coord},L}}\left(R^{\mathrm{coord},L}\widehat{\otimes}_{R_1} J\!L\right)&\cong \Omega_{R^{\mathrm{coord},L}}\widehat{\otimes}_{\Omega_{R_1}}\left(\Omega^{L_1}(R_1)\otimes_{R_1} J\!L\right)\cong C^{\mathrm{coord},L}\widehat{\otimes}_{R_1} J\!L,\\
C^{\mathrm{coord},L}\widehat{\otimes}_{R^{\mathrm{coord},L}}\left(R^{\mathrm{coord},L}\widehat{\otimes} F\right)&\cong C^{\mathrm{coord},L}[\![x_1,\dots,x_d]\!].
\end{aligned}
\]
We endow the graded algebra on the left-hand, resp.\ right-hand, side
of~\eqref{ref-4.2-56} with the following natural differential
\begin{align}
\label{ref-4.3-57} {}^1\nabla^\mathrm{coord}&=\mathrm d_{\Omega_{R^{\mathrm{coord},L}}}\widehat\otimes_{\Omega_{R_1}} 1+1\widehat\otimes_{\Omega_{R_1}}{{}^1\nabla},\ \text{resp.}\\
\label{ref-4.4-58} \mathrm d&=\mathrm d_{C^{\mathrm{coord},L}}\widehat\otimes 1, 
\end{align}
where ${}^1\nabla$ has been introduced in Subsection~\ref{ref-3.4-38}.
Both~\eqref{ref-4.3-57} and~\eqref{ref-4.4-58} are, by construction, flat $C^{\mathrm{coord},L}$-connections on the respective spaces, and the obvious inclusions from $C^{\mathrm{coord},L}$ into $C^{\mathrm{coord},L}\widehat{\otimes}_{R_1} J\!L$ and $C^{\mathrm{coord},L}[\![x_1,\dots,x_d]\!]$ are morphisms of DG-algebras.

The main property of the connections~\eqref{ref-4.3-57}
and~\eqref{ref-4.4-58} lies in the existence of a {\bf canonical Maurer--Cartan
  {}{element} in $C^{\mathrm{coord},L}$}: namely, according
to~\cite[Subsection 1.6]{VdB} and~\cite[Subsection 5.2]{CVdB}, there
exists a unique element $\omega$ of
$C^{\mathrm{coord},L}\widehat{\otimes}\mathrm{Der}(F)$ of degree $1$,
satisfying
\[
t\circ {}^1\nabla^\mathrm{coord}\circ t^{-1}-\mathrm d=\omega,
\] 
where the expression on the left-hand side is naturally viewed as a
$C^{\mathrm{coord},L}$-linear derivation of $F$.  Furthermore,
$\omega$ satisfies the Maurer--Cartan equation in the DG-Lie algebra
$C^{\mathrm{coord},L}\widehat{\otimes}\mathrm{Der}(F)$, i.e.\
\[
\mathrm d\omega+\frac{1}2[\omega,\omega]=0,
\]
(which implies that $\mathrm d+[\omega,\bullet]$ is a flat connection on $C^{\mathrm{coord},L}[\![x_1,\dots,x_d]\!]$) and the verticality condition
\begin{equation}
\label{ref-4.5-59}
\iota_v\omega=1\otimes v,\ v\in \mathfrak{gl}_d(k)
\end{equation}
(here, $\iota_v$ on the left-hand side denotes the contraction operation on $C^{\mathrm{coord},L}$ with respect to $v$, coming from the infinitesimal action of $\mathfrak{gl}_d(k)$ on $R^{\mathrm{coord},L}$; $v$ on the right-hand side denotes the linear vector field associated to $v$, acting on $F$).

Finally, we consider the {\bf affine coordinate space $R^{\mathrm{aff},L}$} of a Lie algebroid $L$ over $R$: it is simply the $\mathrm{GL}_d(k)$-invariant ring 
\[
R^{\mathrm{aff},L}=\left(R^{\mathrm{coord},L}\right)^{\mathrm{GL}_d(k)}
\]
It is an $R$-algebra in an obvious way, and
enjoys a universal property similar to the one satisfied by
$R^{\mathrm{coord},L}$, for which we refer to~\cite[Subsection
5.4]{CVdB}.
\begin{Exa}
Continuing Example \ref{exa-poly}, assume $R=k[x_1,\dots,x_d]$ and $L=\mathrm{Der}_k(R)$. We now have
\[
R^{\mathrm{aff},L}=R[y_{i,\alpha}:\ i=1,\dots,d,\ |\alpha|\geq 2],
\]
where $|\bullet|$ denotes the norm of a multiindex in $\mathbb N^d$.
We observe that $R^{\mathrm{aff},L}$ is an (infinite) polynomial ring,
while $R^{\mathrm{coord},L}$ is not, due to the localization.
\end{Exa}
Similarly, we have the DG-algebra
$C^{\mathrm{aff},L}=\Omega_{R^{\mathrm{aff},L}}\otimes_{\Omega_{R_1}}\Omega^{L_1}(R_1)$,
with differential $\mathrm d_{C^{\mathrm{aff},L}}=\mathrm
d_{\Omega_{R^{\mathrm{aff},L}}}\otimes_{\Omega_{R_1}}1+1\otimes_{\Omega_{R_1}}\mathrm
d_{L_1}$.  We may further consider the graded algebra
\[
C^{\mathrm{aff},L}\widehat{\otimes}_{R^{\mathrm{aff},L}}\left(R^{\mathrm{aff},L}\widehat{\otimes}_{R_1}
  J\!L\right)\cong
\Omega_{R^{\mathrm{aff},L}}\widehat{\otimes}_{\Omega_{R_1}}\left(\Omega^{L_1}(R_1)\otimes_{R_1}
  J\!L\right)\cong C^{\mathrm{aff},L}\widehat{\otimes}_{R_1} J\!L,
\]
endowed with the natural differential 
\[
{}^1\nabla^{\mathrm{aff}}=\mathrm d_{\Omega_{R^{\mathrm{aff},L}}}\widehat\otimes_{\Omega_{R_1}} 1+1\widehat\otimes_{\Omega_{R_1}} {{}^1\nabla},
\]
making the natural inclusion $C^{\mathrm{aff},L}\hookrightarrow C^{\mathrm{aff},L}\widehat{\otimes}_{R_1}J\!L$ into a morphism of DG-algebras.
Obviously, ${}^1\nabla^\mathrm{coord}$ descends by its very construction to $C^{\mathrm{aff},L}\widehat\otimes_{R_1}J\!L$ and identifies with ${}^1\nabla^\mathrm{aff}$.
\begin{Lem} \label{ref-4.1-60} $R^{\aff,L}$ is of the form $S\otimes R$ where $S$ is an (infinitely generated) polynomial ring.
\end{Lem}
\begin{proof}
See \cite[\S5.3]{CVdB}.
\end{proof}
Note that {}{this depends on our standing assumption that $L$ is free and furthermore} the decomposition $R^{\aff,L}=S\otimes R$ is not canonical. 

\subsection{Fedosov resolutions of $L$-poly-vector fields and $L$-poly-differential operators on $R$}\label{ref-4.2-61}

In this section, we recall briefly the main results
of~\cite[\S4.3]{CVdB}, to which we refer for more details.
We consider relative poly-differential operators and poly-vector fields (see
\S\ref{ref-3.2-32}) 
in the following situation: $(A,\mathrm d_A)=(C^{\mathrm{aff},L},\mathrm d_{
  C^{\mathrm{aff},L}})$ and $(B,\mathrm
d_B)=(C^{\mathrm{aff},L}\widehat{\otimes}_{R_1}J\!L,{}^1\nabla^\mathrm{aff})$.
\begin{Thm}\label{ref-4.2-62}
For a Lie algebroid $L$ over $R$ as above, there exist quasi-isomorphisms of Gerstenhaber algebras up to homotopy
\begin{align}
\label{ref-4.6-63}(T_\mathrm{poly}^L(R),0,[\ ,\ ],\cup)=(T_\mathrm{poly}^{L_2}(R_2),0,[\ ,\ ],\cup)&\hookrightarrow \left(T_{\mathrm{poly},C^{\mathrm{aff},L}}(C^{\mathrm{aff},L}\widehat{\otimes}_{R_1}J\!L),{}^1\nabla^\mathrm{aff},[\ ,\ ],\cup\right),\\
\label{ref-4.7-64}(D_\mathrm{poly}^L(R),\mathrm d_\mathrm H,[\ ,\
],\cup)=(D_\mathrm{poly}^{L_2}(R_2),\mathrm d_\mathrm H,[\ ,\
],\cup)&\hookrightarrow
\left(D_{\mathrm{poly},C^{\mathrm{aff},L}}(C^{\mathrm{aff},L}\widehat{\otimes}_{R_1}J\!L),{}^1\nabla^\mathrm{aff}+\mathrm
  d_\mathrm H,[\ ,\ ],\cup\right).
\end{align} 
\end{Thm} 
\begin{proof}
We refer to \cite{CVdB} for details. For example the map \eqref{ref-4.7-64}
is derived by suitable base extension from \eqref{ref-3.12-37}. For the fact that
the maps are quasi-isomorphisms we refer to \cite[Prop.\ 7.3.1]{CVdB}.
\end{proof}

\subsection{The Fedosov resolution of $L$-forms on $R$}\label{ref-4.3-65}

We consider the precalculus $(\Omega^L(R),0,\mathrm L,\cap)$ of
$L$-forms over the Gerstenhaber algebra $(T_\mathrm{poly}^L(R),0,[\ ,\
],\cup)$, described in \S\ref{ref-3.1-17}: we describe now a well-suited
resolution of $(\Omega^L(R),0,\mathrm L,\cap)$ which is compatible
with the Fedosov resolution
$\left(T_{\mathrm{poly},C^{\mathrm{aff},L}}(C^{\mathrm{aff},L}\widehat{\otimes}_{R_1}J\!L),{}^1\nabla^\mathrm{aff},[\
  ,\ ],\cup\right)$ from Theorem~\ref{ref-4.2-62}.
\begin{Thm}\label{ref-4.3-66}
  For a Lie algebroid $L$ over $R$ as above, there exists a
  quasi-isomorphism of precalculi
  as in the following commutative diagram \footnote{$\Omega_A$, for a topological
  $k$-algebra $A$, denotes the continuous De Rham complex. A similar convention holds
for an extension of topological algebras $B/A$.}:
\[
\xymatrix{(T_\mathrm{poly}^L(R),0,[\ ,\ ],\cup)=(T_\mathrm{poly}^{L_2}(R_2),0,[\ ,\ ],\cup)\ar@{^{(}->}[r]\ar@{~>}[d] &\left(T_{\mathrm{poly},C^{\mathrm{aff},L}}(C^{\mathrm{aff},L}\widehat{\otimes}_{R_1}J\!L),{}^1\nabla^\mathrm{aff},[\ ,\ ],\cup\right)\ar@{~>}[d]\\
(\Omega^L(R),0,\mathrm L,\cap)=(\Omega^{L_2}(R_2),0,\mathrm L,\cap)\ar@{^{(}->}[r] &\left(\Omega_{C^{\mathrm{aff},L}\widehat{\otimes}_{R_1}J\!L/C^{\mathrm{aff},L}},{}^1\nabla^\mathrm{aff},\mathrm L,\cap\right)},
\]
the vertical arrows denoting the contraction and Lie derivative.
\end{Thm}
\begin{proof}
  We refer to~\cite[\S  4.3.3]{CVdB}: we observe that the
  construction of the quasi-isomorphism uses a dualization of the
  construction of the quasi-isomorphism~\eqref{ref-4.6-63}, and that
  contraction operations and differentials are preserved by the above
  quasi-isomorphism, whence all algebraic structures are preserved.
\end{proof}

\subsection{The Fedosov resolution of $L$-chains on $R$}\label{ref-4.4-67}

We consider the DG-algebra
$(C^{\mathrm{aff},L}\widehat{\otimes}_{R_1}J\!L,{}^1\nabla^\mathrm{aff})$,
and to it we associate the $C^{\mathrm{aff},L}$-relative Hochschild
chain complex, i.e.\
\[
\begin{aligned}
  \widehat{C}_{C^{\mathrm{aff},L},\bullet}(C^{\mathrm{aff},L}\widehat{\otimes}_{R_1}J\!L)&=
  \bigoplus_{p\leq 0}
  \left(C^{\mathrm{aff},L}\widehat{\otimes}_{R_1}J\!L\right)^{\widehat{\otimes}_{C^{\mathrm{aff},L}}(-p+1)}\\
  & \cong \bigoplus_{p\leq 0} \left(C^{\mathrm{aff},L}\widehat{\otimes}_{R_1}J\!L^{\widehat{\otimes}_{R_1}(-p+1)}\right)=\\
  &=
  C^{\mathrm{aff},L}\widehat\otimes_{R_1}\widehat{C}_{R_1,\bullet}(J\!L),
\end{aligned}
\]
Further, we have the identification
\[
C^{\mathrm{aff},L}\widehat{\otimes}_{R_1} \widehat{C}_{R_1,\bullet}(J\!L)\cong
\Omega_{R^{\mathrm{aff},L}}\widehat{\otimes}_{\Omega_{R_1}}\left(\Omega^{L_1}(R_1)\widehat{\otimes}_{R_1}
\widehat{C}_{R_1,\bullet}(J\!L)\right).
\]
and one checks that the differentials coming from the Grothendieck
connection on each side are the same. I.e.
\[
{}^1\nabla^\mathrm{aff}=\mathrm d_{\Omega_{R^{\mathrm{aff},L}}}
\,\widehat{\otimes}\,1+1\,\widehat{\otimes}\,{}^1\nabla,
\]
\begin{Prop}\label{ref-4.4-68}
  For a Lie algebroid $L$ over $R$ as above, the cohomology of
  $\left(\widehat{C}_{C^{\mathrm{aff},L},\bullet}(C^{\mathrm{aff},L}\widehat{\otimes}_{R_1}J\!L),{}^1\nabla^\mathrm{aff}\right)$
  is concentrated in degree $0$, where
\[
\mathrm H^0\!\left(\widehat{C}_{C^{\mathrm{aff},L},\bullet}(C^{\mathrm{aff},L}\widehat{\otimes}_{R_1}J\!L),{}^1\nabla^\mathrm{aff}\right)\cong C_\mathrm{poly}^L(R).
\]
\end{Prop}
\begin{proof}
Taking the inverse of \eqref{ref-3.15-44} we obtain a morphism
\[
 C_{\mathrm{poly}}^L(R)\cong \widehat{C}_{R_1,\bullet}(J\!L)^{{}^1\nabla}
\hookrightarrow \widehat{C}_{R_1,\bullet}(J\!L)
\]
which extends to a morphism 
\begin{equation}
(C_{\mathrm{poly}}^L(R),0)\rightarrow (\Omega_{R^{\mathrm{aff},L}}\widehat{\otimes}_{\Omega_{R_1}}\left(\Omega^{L_1}(R_1)\widehat{\otimes}_{R_1}\widehat{C}_{R_1,\bullet}(J\!L)\right)
,\mathrm d_{\Omega_{R^{\mathrm{aff},L}}}
\,\otimes\,1+1\,\otimes\,{}^1\nabla)
\end{equation}
We will show that it is a quasi-isomorphism.
To this end we make use of the identification $R^{\aff,L}=S\otimes R$ given in 
Lemma \ref{ref-4.1-60}. The right-hand side of the extended morphism becomes
\[
(\Omega_S\,\widehat{\otimes}\,\ \left(\Omega^{L_1}(R_1)\widehat{\otimes}_{R_1}\widehat{C}_{R_1,\bullet}(J\!L)\right)
,\mathrm d_S\otimes 1+1\otimes {}^1\nabla)
\]
Using a filtration argument together with a suitable version of Poincar\'e's
Lemma for $S$, the previous complex is quasi-isomorphic to
\[
(\Omega^{L_1}(R_1)\otimes_{R_1}\widehat{C}_{R_1,\bullet}(J\!L),{}^1\nabla)
\]
It remains to show that for each $p\leq0$
\[
(\Omega^{L_1}(R_1)\otimes_{R_1}\widehat C_{R_1,p}(J\!L),{}^1\nabla)
\]
has cohomology in degree $0$. 
Filtering this complex with respect to the $J$-adic
filtration and taking the associated graded complex one verifies that
one obtains 
\[
(\Omega^{L}(R)\otimes_{R_1}\mathrm S(L^\ast)^{\otimes -p-1},\mathrm d)
\]
where the differential $\mathrm d$ is obtained from the action of  
$L$ on $\mathrm S(L^\ast)^{\otimes p+1}$ by contraction. Using again a suitable version
of Poincar\'e's Lemma one finds that the resulting complex is indeed exact in
degrees $<0$. 
\end{proof}
\begin{Thm}\label{ref-4.5-69}
For a Lie algebroid $L$ over $R$ as above, there is a quasi-isomorphism of 
precalculi up to homotopy as in the following commutative diagram:
\[
\xymatrix{(D_\mathrm{poly}^L(R),\mathrm d_\mathrm H,[\ ,\ ],\cup)=(D_\mathrm{poly}^{L_2}(R_2),\mathrm d_\mathrm H,[\ ,\ ],\cup)\ar@{^{(}->}[r]\ar@{~>}[d] &\left(D_{\mathrm{poly},C^{\mathrm{aff},L}}(C^{\mathrm{aff},L}\widehat{\otimes}_{R_1}J\!L),{}^1\nabla^\mathrm{aff}+\mathrm d_\mathrm H,[\ ,\ ],\cup\right)\ar@{~>}[d]\\
  (C_\mathrm{poly}^L(R),\mathrm b_\mathrm H,\mathrm
  L,\cap)\ar@{^{(}->}[r] &
  \left(\widehat{C}_{C^{\mathrm{aff},L},\bullet}(C^{\mathrm{aff},L}\widehat{\otimes}_{R_1}J\!L),{}^1\nabla^\mathrm{aff}+\mathrm
    b_\mathrm H,\mathrm L,\cap\right)},
\]
the vertical arrows denoting the contraction and Lie derivative.
\end{Thm}

\section{Globalization of Tsygan's formality in the Lie algebroid framework}\label{ref-5-70}

The present section is devoted to the proof of Theorem~\ref{ref-1.2-11}.
We first briefly review some basic facts on $L_\infty$-algebras,
$L_\infty$-modules and related morphisms. This is discussed in
\cite[\S6]{CVdB} for $L_\infty$-morphisms. Here we add a discussion on
the descent procedure for $L_\infty$-modules over $L_\infty$-algebras
and related morphisms.

Then, we add a short {\em excursus} on Kontsevich's and Shoikhet's formality theorems: we focus on the main properties of both formality morphisms, without delving into the technical details of their respective constructions.

Finally, we give the main lines, along which the globalization of Tsygan's
formality can be proved: the proof is a combination of the properties
of Kontsevich's and Shoikhet's $L_\infty$-morphisms with the Fedosov
resolutions from \S\ref{ref-4-53}.

\subsection{Descent for $L_\infty$-algebras and $L_\infty$-modules.}

We discuss a series of descent scenarios for $L_\infty$-algebras,
$L_\infty$-modules and related morphisms, which are modelled after the
formalism for descent of differential forms in differential geometry.
The verification of the results in this section are along the same
lines as~\cite[\S7.6,\ \S7.7]{VdB}.
To clearly separate all the various cases we have numbered them.

\medskip

\noindent
\textbf{(1)} To start it is convenient to work over an arbitrary DG operad
$\mathcal{O}$ with underlying graded operad $\widetilde{\mathcal{O}}$ (thus we
forget the differential on $\mathcal{O})$.  Assume
that $\frak{g}$ is an algebra over $\mathcal{\mathcal{O}}$ and consider a set of
$\widetilde{\mathcal O}$-derivations $(\iota_v)_{v\in \mathfrak{s}}$ of degree
$-1$ on $\frak{g}$ ($\mathfrak{s}$ is an index set, without any additionnal structure). 
Put $\mathrm L_v=\mathrm d_{\frak{g}}\iota_v+\iota_v\mathrm d_{\frak{g}}$.  This is a
derivation of $\frak{g}$ of degree zero which commutes with $\mathrm d_{\frak{g}}$.
Put 
\begin{equation}\label{ref-5.1-71}
\frak{g}^{\mathfrak{s}}=\{w\in \frak{g}\mid \forall v\in \mathfrak{s}:\iota_v w= \mathrm L_v w=0\}
\end{equation}
It is easy to see that $\frak{g}^{\mathfrak{s}}$ is an algebra over $\mathcal{O}$ as well. Informally
we will call such a set of derivations $(\iota_v)_{v\in \mathfrak{s}}$ an $\mathfrak{s}$-action.

\medskip

\noindent
\textbf{(2)} Assume that $M$ is a $\frak{g}$-module and assume that $\mathfrak s$
also acts on $M$, in a way compatible with the action of $\mathfrak{s}$ on
$\mathfrak g$, i.e.\ a general element $v$ of $\mathfrak s$ determines
an operator $\iota_v$ on $M$, such that Leibniz's rule holds true for
the operations $\widetilde{\mathcal{O}}(n)\otimes \left(\mathfrak{g}^{\otimes n-1}
\otimes M\right)\rightarrow M$. Again, we set $\mathrm L_v=\mathrm
d_M\iota_v+\iota_v\mathrm d_M$, which is a derivation of degree~$0$ on
$M$ compatible with the derivations $\mathrm L_v$ on $\mathfrak{g}$, $\mathrm d_M$ being the differential on $M$.

\medskip
\noindent
\textbf{(3)} The above constructions apply in particular if $\frak{g}$ is an $L_\infty$-algebra.  Assume that it has Taylor coefficients $Q_n$,
$n\geq 1$. Then $\mathrm L_v$ is defined by means of $\mathrm
d_{\mathfrak{g}}=Q_1$, and the derivation property of $\iota_{v}$ reads as
\begin{equation}
\label{ref-5.2-72}
\iota_v \left(Q_n(x_1,\dots,x_n)\right)=\sum_{i=1}^n(-1)^{\sum_{j=1}^{i-1}|x_j|+i}Q_n(x_1,\dots,\iota_v x_i,\dots,x_n),\ x_j\in\mathfrak g,\ j=1,\dots,n.
\end{equation}
Under these conditions the $L_\infty$-structure descends to $\mathfrak
g^\mathfrak s$. 

\medskip
\noindent
\textbf{(4)}
Similarly if $M$ is an $L_\infty$-module over
$\mathfrak{g}$ defined by Taylor coefficients $R_n$ then the compatibility
condition is
\begin{equation}
\label{ref-5.3-73}
\begin{aligned}
  \iota_v \left(R_n(x_1,\dots,x_n;m)\right)&=\sum_{i=1}^n(-1)^{\sum_{j=1}^{i-1}|x_j|+i}R_n(x_1,\dots,\iota_v x_i,\dots,x_n;m)+\\
  &\phantom{=}+(-1)^{\sum_{i=1}^n|x_i|+n-1} R_n(x_1,\dots,x_n;\iota_v
  m),\ m\in M,\ x_j\in\mathfrak g,\ j=1,\dots,n.
\end{aligned}
\end{equation}
If this holds true then $M^{\mathfrak{s}}$ becomes an $L_\infty$-module over
$\mathfrak{g}^{\mathfrak{s}}$.

\medskip
\noindent
\textbf{(5)} We also need descent for $L_\infty$-morphisms. This does not
immediately fall under the operadic framework given in (1), (2) but it is
easy enough to give explicit formul\ae\ like
\eqref{ref-5.2-72}, \eqref{ref-5.3-73}. Thus assume $\psi:\mathfrak{g}\rightarrow
\mathfrak{h}$ is an $L_\infty$-morphism between $L_\infty$-algebras
with $\mathfrak{s}$-action.  Under the following compatibility
condition
\begin{equation}
\label{ref-5.4-74}
\iota_v(\psi_n(x_1,\ldots,x_n))=\sum_{i=1}^n (-1)^{\sum_{j=1}^{i-1}|x_j|+(i-1)} \psi_n(x_1,\ldots,
\iota_vx_i,\ldots,x_n)
\end{equation}
$x_j\in\mathfrak g,\ j=1,\dots,n,\ n\geq 1$, $\psi$ descends to an
$L_\infty$-morphism $\psi^{\mathfrak{s}}:\frak{g}^{\mathfrak{s}}\rightarrow\frak{h}^{\mathfrak{s}}$. 

\medskip
\noindent
\textbf{(6)} Let $\psi:\mathfrak{g}\rightarrow\mathfrak{h}$ be a morphism between
$L_\infty$-algebras with $\mathfrak{s}$-action such that the descent 
condition \eqref{ref-5.4-74} holds, and let $N$ be an $L_\infty$-module over
$\mathfrak{h}$ equipped with a compatible $\mathfrak{s}$-action.
Let $N_\psi$ be the pullback of $N$ along $\psi$. Then
the $\mathfrak{s}$-action on $N_\psi$ is compatible with the $\mathfrak{s}$-action
on $\mathfrak{g}$. 

\medskip
\noindent
\textbf{(7)}
Now assume that $\mathfrak{g}$ is an $L_\infty$-algebra and $M$, $N$
are $L_\infty$-modules over $\mathfrak{g}$. Assume that all objects
are equipped with an $\mathfrak{s}$-action and that the descent
conditions are satisfied.

Assume that $\varphi:M\rightarrow N$ is an $L_\infty$-module
morphism. Then the condition for $\varphi$ to descend to an
$L_\infty$-morphism $M^{\mathfrak{s}}\rightarrow N^{\mathfrak{s}}$ is
\begin{equation}
\label{ref-5.5-75}
\begin{aligned}
\iota_v\left(\varphi_n(x_1,\dots,x_n;m)\right)&=\sum_{i=1}^n (-1)^{\sum_{j=1}^{i-1}|x_j|+(i-1)} \varphi_n(x_1,\dots,\iota_v x_i,\dots,x_n;m)+\\
&\phantom{=}+(-1)^{\sum_{i=1}^n|x_i|+n}\varphi_n(x_1,\dots,x_n;\iota_v m),
\end{aligned}
\end{equation}
for $m\in M,\ x_j\in\mathfrak g,\ j=1,\dots,n,\ n\geq 1$.
\subsection{Twisting of $L_\infty$-algebras and $L_\infty$-modules}\label{ref-5.2-76}
We refer to~\cite[\S2]{Dol}, for a very detailed exposition of
$L_\infty$-algebras, $L_\infty$-modules and the associated twisting
procedures.  See also \cite{Y2}. 

\medskip

\noindent \textbf{Convention. }We will work with infinite sums. We assume throughout that the occurring sums
are convergent and that standard series manipulations are allowed. This
will be the case in our applications. 

\medskip

If $(\mathfrak{g},Q)$ is an $L_\infty$-algebra then the Maurer-Cartan equation 
is defined as
\begin{equation}
\label{ref-5.6-77}
\sum_{j=1}^\infty \frac{1}{j!} Q_n(\underbrace{\omega,\cdots,\omega}_j)=0,
\end{equation}
and a solution $\omega\in \mathfrak{g}_1$ is called a Maurer--Cartan
element (MC element for short).   Below we will only use DG-Lie
algebras and in this case \eqref{ref-5.6-77} reduces to the finite sum
\[
d\omega+\frac{1}{2}[\omega,\omega]=0.
\]
A MC element defines a new ``twisted'' DG-Lie structure on
$\mathfrak{g}$ (denoted by $\mathfrak{g}_\omega$) with Taylor coefficients
\[
Q_{\omega,n}(x_1,\ldots,x_n)=\sum_{j} \frac{1}{j!} Q_{n+j}(
\underbrace{\omega,\ldots,\omega}_j, x_1,\ldots,x_n),\quad n\ge 1
\]
If $\mathfrak{g}$ is actually a DG-Lie algebra then twisting keeps the
bracket but changes the differential to
\[
\mathrm d_{\omega}=\mathrm d_{\mathfrak{g}}+[\omega,-].
\]
If $\mathfrak{h}$ is another $L_\infty$ algebra, $\psi$ is an
$L_\infty$-morphism from $\mathfrak g$ to $\mathfrak h$ and $\omega$
is a MC element in $\mathfrak g$ then
\begin{equation}
\label{ref-5.7-78}
\psi(\omega)=\sum_{n\geq 1}\frac{1}{n!}\psi_n(\underset{n}{\underbrace{\omega,\dots,\omega}})
\end{equation}
is a MC element in $\mathfrak h$.

We may also twist $\psi$ with respect to $\omega$, so as to
get an $L_\infty$-morphism $\psi_\omega$ from $\mathfrak g_\omega$ to
$\mathfrak g_{\psi(\omega)}$, where
\[
\psi_{\omega,n}(x_1,\ldots,x_n)=\sum_{j\geq 0}\frac{1}{j!}\psi_{n+j}(\underset{j}{\underbrace{\omega,\dots,\omega}},x_1,\ldots,x_n),\ n\geq 1.
\]
If $M$ is an $L_\infty$-module over a DG-Lie algebra with Taylor 
coefficients $R_n$ and $\omega\in
\mathfrak{g}_1$ is a MC element then we may define a twisted $L_\infty$
structure on $M_\omega$ over $\mathfrak{g}_\omega$ by the formula
\[
R_{\omega,n}(x_1,\ldots,x_n;m)=\sum_{j\ge 0}\frac{1}{j!}R_{n+j}(
\underbrace{\omega,\ldots,\omega}_j,x_1,\ldots,x_n;m),\quad n\ge 0
\]
If $\mathfrak{g}$ is a DG-Lie algebra and $M$ is a DG-Lie module over
$\mathfrak{g}$ then twisting keeps the $\mathfrak{g}$-action on $M$ but
changes the differential on $M$ to
\[
\mathrm d_\omega=d+\omega\bullet.
\]
Twisting of modules is compatible with pullback. More precisely if
$\psi:\mathfrak{g}\rightarrow\mathfrak{h}$ is an $L_\infty$-morphism,
$N$ is an $L_\infty$-module over $\mathfrak{h}$ and $\omega\in
\mathfrak{g}_1$ is a MC element then we have
\begin{equation}
\label{ref-5.8-79}
(N_{\psi(\omega)})_{\psi_\omega}=(N_\psi)_\omega
\end{equation}
If $\varphi:M\rightarrow N$ is an $L_\infty$-morphism of DG-Lie modules over the
DG-Lie algebra $\mathfrak{g}$ and $\omega$ is a MC element in $\mathfrak{g}_1$
then we obtain a twisted $L_\infty$-morphism $\varphi_\omega:M_\omega\rightarrow N_\omega$ which is defined by 
\begin{equation}
\label{ref-5.9-80}
\varphi_{\omega,n}(x_1,\ldots,x_n;m)= \sum_{j\geq
  0}\frac{1}{j!}\varphi_{n+j}(\underset{j}{\underbrace{\omega,\dots,\omega}},x_1,\ldots,x_n;m),\
n\geq 1.
\end{equation}
\subsection{Compatibility of twisting and descent}
\label{ref-5.3-81}
Assume now that $\mathfrak{g}$ is a DG-Lie algebra equipped with an
$\mathfrak{s}$-action and that $\omega\in \mathfrak{g}_1$ is a MC
element.  
Then $\mathfrak{s}$ still acts on $\mathfrak{g}_\omega$, where we forget here about the differential: in fact,
the concept of an $\mathfrak{s}$-action only refers to the underlying
Lie algebra structure on $\mathfrak{g}$. 
However $\mathfrak{g}^{\mathfrak{s}}$ and $\mathfrak{g}_\omega^{\mathfrak{s}}$
will be different (as the Lie derivative $\mathrm L_v$ for
$v\in\mathfrak{s}$ will be different).

If $(M,R)$ is an $L_\infty$-module over~$\mathfrak{g}$ which is also equipped with a compatible
$\mathfrak{s}$-action then the $\mathfrak{s}$-actions on
$\mathfrak{g}_\omega$ and $M_\omega$ are compatible provided the following
condition holds 
\begin{equation}
\label{ref-5.10-82}
R_n(\iota_v\omega,x_2,\dots,x_n;m)=0, x_i\in\mathfrak g,\ i=2,\dots,n,\ n\geq 2,\ m\in M.
\end{equation}
This condition is automatic if $M$ is a DG-Lie module. 

\medskip

If $\psi:\mathfrak{g}\rightarrow\mathfrak{h}$ is an
$L_\infty$-morphism of DG-Lie algebras equipped with an
$\mathfrak{s}$-action and the descent condition \eqref{ref-5.4-74} is
satisfied for $\psi$ then an easy computation (see e.g.\
\cite[\S7.7]{VdB}) shows that the same descent condition will be
satisfied for $\psi_\omega$ if the following condition holds
\begin{equation}
\label{ref-5.11-83}
\psi_n(\iota_v\omega,x_2,\dots,x_n)=0,\ x_i\in\mathfrak g,\ i=2,\dots,n\ s\in \mathfrak s,\ n\geq 2,
\end{equation}
Furthermore if in this setting $N$ is an $L_\infty$-module  over $\frak{h}$
with compatible $\mathfrak{s}$-action such that the compatibility condition
\eqref{ref-5.10-82} holds then the corresponding condition will
hold for $N_\psi$.

Similarly if we have an $L_\infty$-morphism $\varphi:M\rightarrow N$
between DG-Lie modules over a DG-Lie algebra $\frak{g}$ such that
$\frak{g}$, $M$, $N$ are equipped with compatible
$\mathfrak{s}$-actions in such a way that the descent condition \eqref{ref-5.5-75}
holds for $\varphi$ then the
same descent condition will be satisfied for $\varphi_\omega$ if 
the following condition holds 
\[
\varphi_n(\iota_v\omega,x_2,\dots,x_n;m)=0,\ m\in M,\ x_i\in\mathfrak g,\ i=2,\dots,n\ n\geq 2,\ s\in \mathfrak s.
\]

\subsection{Kontsevich's and Shoikhet's formality theorems}\label{ref-5.4-84}

In this brief section, we quote (without proofs) Kontsevich's and
Shoikhet's formality theorems, along with the relevant properties, which
we will need later in the proof of globalization results.

We consider the algebra $F=k[\![x_1,\dots,x_d]\!]$ of formal power series in $d$ variables over a field $k$ containing $\mathbb R$.

To $F$, we associate the DG-Lie algebras $(T_\mathrm{poly}(F),0,[\ ,\
])$, resp.\ $(D_\mathrm{poly}(F),\mathrm d_\mathrm H,[\ ,\ ])$, of
formal poly-vector fields, resp.\ formal poly-differential operators, on
$F$; further, we consider the DG-Lie modules $(\Omega_F,0,\mathrm L)$,
resp.\ $(\widehat{C}_\bullet(F),\mathrm b_\mathrm H,\mathrm L)$,
over $(T_\mathrm{poly}(F),0,[\ ,\ ])$, resp.\
$(D_\mathrm{poly}(F),\mathrm d_\mathrm H,[\ ,\ ])$, where $\Omega_F$
denotes the continuous De Rham complex of~$F$ with De Rham
differential $\mathrm d$, and $\widehat{C}_\bullet(F)$ is the continuous Hochschild
chain complex of $F$.

The following is Kontsevich's celebrated ``formality'' result. 
\begin{Thm}[\cite{K}]\label{ref-5.1-85}
There is an $L_\infty$-quasi-isomorphism
\[
\mathcal U:(T_\mathrm{poly}(F),0,[\ ,\ ])\to (D_\mathrm{poly}(F),\mathrm d_\mathrm H,[\ ,\ ]),
\]
enjoying the following properties:
\begin{enumerate}
\item[$i)$] the first Taylor coefficient of $\mathcal U$ coincides with the Hochschild--Kostant--Rosenberg quasi-isomorphism
\[
\mathrm{HKR}(\partial_{i_1}\wedge\cdots\wedge \partial_{i_p})=(-1)^{\frac{p(p-1)}2}\frac{1}{p!}\sum_{\sigma\in\mathfrak S_p}(-1)^\sigma \partial_{i_{\sigma(1)}}\otimes\cdots\otimes\partial_{i_{\sigma(p)}}
\]
from the DG-vector space $(T_\mathrm{poly}(F),0)$ to the DG-vector space $(D_\mathrm{poly}(F),\mathrm d_\mathrm H)$.
\item[$ii)$] If $n\geq 2$, and $\gamma_i$, $i=1,\dots,n$, are elements of $T_\mathrm{poly}^0(F)$, then
\begin{equation*}
\mathcal U_n(\gamma_1,\dots,\gamma_n)=0.
\end{equation*} 
\item[$iii)$] If $n\geq 2$, $\gamma_1$ is a linear vector field on $F$ (i.e.\ an element of $\mathfrak{gl}_d$), $\gamma_i$, $i=2,\dots,n$ are general elements of $T_\mathrm{poly}(F)$, then
\begin{equation*}
\mathcal U_n(\gamma_1,\gamma_2,\dots,\gamma_n)=0.
\end{equation*} 
\end{enumerate}
\end{Thm}
By composing the action $\mathrm L$ of $D_\mathrm{poly}(F)$ on
$\widehat{C}_\bullet(F)$ with the $L_\infty$-quasi-isomorphism~$\mathcal U$ from
Theorem~\ref{ref-5.1-85}, $\widehat{C}_\bullet(F)$ inherits an $L_\infty$-module structure
over the DG-Lie algebra $(T_\mathrm{poly}(F),0,[\ ,\ ])$.

The first part of the following theorem was a conjecture by Tsygan \cite{T}
which has been proved by Shoikhet in \cite{Sh}. The second part has
been proved in~\cite{Dol}.
\begin{Thm}\label{ref-5.2-86}
There is an $L_\infty$-quasi-isomorphism 
\[
\mathcal S:(\widehat{C}_\bullet(F),\mathrm b_\mathrm H,\mathrm L\circ \mathcal U)\to (\Omega_F,0,\mathrm L)
\]
of $L_\infty$-modules over the DG-Lie algebra $(T_\mathrm{poly}(F),0,[\ ,\ ])$, enjoying the following properties:
\begin{enumerate}
\item[$i)$] the $0$-th Taylor coefficient of $\mathcal S$ coincides with the Hochschild--Kostant--Rosenberg quasi-isomorphism
\[
\mathrm{HKR}((a_0|\cdots|a_p))=\frac{1}{p!}a_0\mathrm da_1\cdots\mathrm d a_p
\]
from the DG-vector space $(\widehat{C}_\bullet(F),\mathrm b_\mathrm H)$ to the DG-vector space $(\Omega_F,0)$.
\item[$ii)$] If $n\geq 1$, $\gamma_1$ is a linear vector field on $F$,
  $\gamma_i$, $i=2,\dots,n$ are general elements of
  $T_\mathrm{poly}(F)$, and $c$ is a general element of
  $\widehat{C}_\bullet(F)$, then
\begin{equation*}
\mathcal S_n(\gamma_1,\dots,\gamma_n;c)=0.
\end{equation*} 
\end{enumerate}
\end{Thm}

\subsection{Formality theorem in the ring case}\label{ref-5.5-87}

This section is devoted to the proof of a Tsygan-like formality theorem in
the case of a Lie algebroid $L$ over a $k$-algebra $R$, such that $L$
is free over $R$ of rank $d$: the proof combines Shoikhet's formality
theorem~\ref{ref-5.2-86} with the Fedosov resolutions from
\S\ref{ref-4-53}.  
\begin{Thm}\label{ref-5.3-88}
{}{Assume $\mathbb{R}\subset k$.}  For a Lie algebroid $L$ over $R$ as above, there exist DG-Lie
  algebras $(\mathfrak g_i^L,\mathrm d_i,[\ ,\ ]_i)$, DG-Lie modules
  $(\mathfrak m_i^L,\mathrm b_i,\mathrm L_i)$ over $\mathfrak g_i$,
  $i=1,2$, and $L_\infty$-quasi-isomorphisms $\mathfrak U_L$ from
  $\mathfrak g_1^L$ to $\mathfrak g_2^L$ and $\mathfrak S_L$ from
  $\mathfrak m_2^L$ to $\mathfrak m_1^L$, which fit into the following
  commutative diagram:
\begin{equation}\label{ref-5.12-89}
\xymatrix{T_\mathrm{poly}^L(R) \ar@{^{(}->}[r]\ar@{~>}[d]_{\mathrm L} & \mathfrak g_1^L \ar[r]^{\mathfrak U_L}\ar@{~>}[d]_{\mathrm L_1}& \mathfrak g_2^L\ar@{~>}[d]_{\mathrm L_2} & D_\mathrm{poly}^L(R)\ar@{_{(}->}[l]\ar@{~>}[d]_{\mathrm L}\\
\Omega^L(R)\ar@{^{(}->}[r] & \mathfrak m_1^L & \mathfrak m_2^L\ar[l]_{\mathfrak S_L} & C_\mathrm{poly}^L(R)\ar@{_{(}->}[l]},
\end{equation}
such that the induced maps 
\[
T_\mathrm{poly}^L(R)\to \mathrm H^\bullet(D_\mathrm{poly}^L(R),\mathrm d_\mathrm H),\ \mathrm H^\bullet(C_\mathrm{poly}^L(R),\mathrm b_\mathrm H)\to \Omega^L(R)
\]
on (co)homology coincide with the respective HKR-quasi-isomorphisms.
The morphisms indicated by hooked arrows are actual quasi-isomorphisms
of DG-Lie algebras and DG-Lie modules respectively.
\end{Thm}
\begin{proof}
  The first step in the proof of Theorem~\ref{ref-5.3-88} may be
  borrowed from~\cite[\S7.3]{CVdB}.
Namely, we consider the following graded vector spaces:
\begin{align*}
  C^{\mathrm{coord},L}\widehat{\otimes}T_\mathrm{poly}(F)&\cong T_{\mathrm{poly},C^{\mathrm{coord},L}}(C^{\mathrm{coord},L}\widehat{\otimes}F),\\ C^{\mathrm{coord},L}\widehat{\otimes}D_\mathrm{poly}(F)&\cong D_{\mathrm{poly},C^{\mathrm{coord},L}}(C^{\mathrm{coord},L}\widehat{\otimes}F),\\
  C^{\mathrm{coord},L}\widehat{\otimes}\Omega_F&\cong
  \Omega_{C^{\mathrm{coord},L}\widehat{\otimes}F/C^{\mathrm{coord},L}},\
  \\C^{\mathrm{coord},L}\widehat{\otimes}\widehat{C}_\bullet(F)&\cong
  \widehat{C}_{C^{\mathrm{coord},L},\bullet}(C^{\mathrm{coord},L}\widehat{\otimes}F),
\end{align*}
where the DG-algebra $C^{\mathrm{coord},L}$ has been introduced in
\S\ref{ref-4.1-54}, and where
$\widehat{C}_{C^{\mathrm{coord},L},\bullet}(C^{\mathrm{coord},L}\widehat{\otimes}F)$
denotes the $C^{\mathrm{coord},L}$-relative Hochschild chain complex
of the DG-algebra $C^{\mathrm{coord},L}\widehat{\otimes}F$.

The Maurer--Cartan form on $C^{\mathrm{coord},L}\widehat{\otimes}F$
introduced in \S\ref{ref-4.1-54} defines a twisted differential $\mathrm
d_\omega=\mathrm d+\omega$ on the listed graded vector spaces and as
explained in \S\ref{ref-5.2-76} $\mathrm d_\omega$ is compatible with  
the respective DG-Lie algebra and DG-Lie module structures.

\medskip

Thus, formal geometry provides us with the following DG-Lie algebras and 
respective DG-Lie modules:
\[
\xymatrix{\left(T_{\mathrm{poly},C^{\mathrm{coord},L}}(C^{\mathrm{coord},L}\widehat{\otimes}F),\mathrm d_\omega,[\ ,\ ]\right),\ar@{~>}[d] & \left(D_{\mathrm{poly},C^{\mathrm{coord},L}}(C^{\mathrm{coord},L}\widehat{\otimes}F),\mathrm d_\omega+\mathrm d_\mathrm H,[\ ,\ ]\right)\ar@{~>}[d]\\
  \left(\Omega_{C^{\mathrm{coord},L}\widehat{\otimes}F/C^{\mathrm{coord},L}},\mathrm
    d_\omega,\mathrm L\right) &
  \left(\widehat C_{C^{\mathrm{coord},L}}(C^{\mathrm{coord},L}\widehat{\otimes}F),\mathrm
    d_\omega+\mathrm b_\mathrm H,\mathrm L\right)}.
\]
We repeat that, viewing all DG-Lie algebra and DG-Lie module structures above as
$L_\infty$-structures, the differential $\mathrm d_\omega$ is the
twist of the standard structures with respect to the MC element
$\omega$ of
$C^{\mathrm{coord},L}\widehat{\otimes}\mathrm{Der}(F)=T_{\mathrm{poly},C^{\mathrm{coord},L}}^0(C^{\mathrm{coord},L}\widehat{\otimes}F)$.

The $L_\infty$-quasi-isomorphism $\mathcal U$ of Theorem~\ref{ref-5.1-85}
extends $C^{\mathrm{coord},L}$-linearly to an
$L_\infty$-quasi-isomorphism
\[
\mathcal U_L:\left(T_{\mathrm{poly},C^{\mathrm{coord},L}}(C^{\mathrm{coord},L}\widehat{\otimes}F),\mathrm d,[\ ,\ ]\right)\to \left(D_{\mathrm{poly},C^{\mathrm{coord},L}}(C^{\mathrm{coord},L}\widehat{\otimes}F),\mathrm d+\mathrm d_\mathrm H,[\ ,\ ]\right). 
\]
The composition of the DG-Lie action $\mathrm L$ of
$D_{\mathrm{poly},C^{\mathrm{coord},L}}(C^{\mathrm{coord},L}\widehat{\otimes}F)$
on
$\widehat{C}_{C^{\mathrm{coord},L},\bullet}(C^{\mathrm{coord},L}\widehat{\otimes}F)$
with the $L_\infty$-quasi-isomorphism $\mathcal U_L$ endows the latter
graded vector space with a structure of $L_\infty$-module over the
DG-Lie algebra
$T_{\mathrm{poly},C^{\mathrm{coord},L}}(C^{\mathrm{coord},L}\widehat{\otimes}F)$,
which is obtained by $C^{\mathrm{coord},L}$-base extension of the corresponding $L_\infty$-module
structure of $\widehat{C}_\bullet(F)$ over $T_{\poly}(F)$. 

Accordingly, the $L_\infty$-quasi-isomorphism $\mathcal S$ of
Theorem~\ref{ref-5.2-86} extends to an $L_\infty$-quasi-isomorphism of
$L_\infty$-modules
\[
\mathcal
S_L:\left(\widehat{C}_{C^{\mathrm{coord},L},\bullet}(C^{\mathrm{coord},L}\widehat{\otimes}F),
  \mathrm d+\mathrm b_\mathrm H,\mathrm L\circ \mathcal U_{L}\right)
\to
\left(\Omega_{C^{\mathrm{coord},L}\widehat{\otimes}F/C^{\mathrm{coord},L}},\mathrm d,L\right)
\]
both viewed as $L_\infty$-modules over
$T_{\mathrm{poly},C^{\mathrm{coord},L}}(C^{\mathrm{coord},L}\widehat{\otimes}F)$.

As outlined in \S\ref{ref-5.2-76} we may apply the twisting procedures for
$L_\infty$-algebras, $L_\infty$-modules and $L_\infty$-morphisms to
the present case, where the MC element is the Maurer--Cartan form
$\omega$: thus, we get an $L_\infty$-morphism $\mathcal U_{L,\omega}$
\[
\mathcal
U_{L,\omega}:\left(T_{\mathrm{poly},C^{\mathrm{coord},L}}(C^{\mathrm{coord},L}\widehat{\otimes}F),\mathrm
  d_\omega,[\ ,\ ]\right)\to
\left(D_{\mathrm{poly},C^{\mathrm{coord},L}}(C^{\mathrm{coord},L}\widehat{\otimes}F),\mathrm
  d_\omega+\mathrm d_\mathrm H,[\ ,\ ]\right),
\]
where here and below we used Property $ii)$ of Theorem~\ref{ref-5.1-85},
which yields that the MC element $\mathcal U(\omega)$ equals $\omega$.

The $L_\infty$-morphism $U_{L,\omega}$  yields an
$L_\infty$-module structure on
$\widehat{C}_{C^{\mathrm{coord},L},\bullet}(C^{\mathrm{coord},L}\widehat{\otimes}F)$
over the $\omega$-twisted DG-Lie algebra
$T_{\mathrm{poly},C^{\mathrm{coord},L}}(C^{\mathrm{coord},L}\widehat{\otimes}F)$.

Translating \eqref{ref-5.8-79} to the present case we have
\begin{multline*}
\left(\widehat{C}_{C^{\mathrm{coord},L},\bullet}(C^{\mathrm{coord},L}\widehat{\otimes}F),\mathrm
  d_\omega+\mathrm b_\mathrm H,\mathrm L\circ \mathcal
  U_{L,\omega}\right)\\=
\left(\widehat{C}_{C^{\mathrm{coord},L},\bullet}(C^{\mathrm{coord},L}\widehat{\otimes}F),\mathrm
  d+\mathrm b_\mathrm H,\mathrm L\circ \mathcal
  U_{L}\right)_\omega
\end{multline*}
from which we get an $L_\infty$-quasi-isomorphism
\[
\mathcal
S_{L,\omega}:\left(\widehat{C}_{C^{\mathrm{coord},L},\bullet}(C^{\mathrm{coord},L}\widehat{\otimes}F),\mathrm
  d_\omega+\mathrm b_\mathrm H,\mathrm L\circ \mathcal
  U_{L,\omega}\right)\to
\left(\Omega_{C^{\mathrm{coord},L}\widehat{\otimes}F/C^{\mathrm{coord},L}},\mathrm
  d_\omega,\mathrm L\right)
\]
of $L_\infty$-modules.

Using the isomorphism~\eqref{ref-4.2-56} we obtain isomorphisms of DG-Lie algebras and respective DG-Lie modules
\begin{align*}
  \left(T_{\mathrm{poly},C^{\mathrm{coord},L}}(C^{\mathrm{coord},L}\widehat{\otimes}F),\mathrm d_\omega,[\ ,\ ]\right)&\cong \left(T_{\mathrm{poly},C^{\mathrm{coord},L}}(C^{\mathrm{coord},L}\widehat{\otimes}_{R_1}J\!L),{}^1\nabla^\mathrm{coord},[\ ,\ ]\right),\\
  \left(D_{\mathrm{poly},C^{\mathrm{coord},L}}(C^{\mathrm{coord},L}\widehat{\otimes}F),\mathrm d_\omega+\mathrm d_\mathrm H,[\ ,\ ]\right)&\cong \left(D_{\mathrm{poly},C^{\mathrm{coord},L}}(C^{\mathrm{coord},L}\widehat{\otimes}_{R_1}J\!L),{}^1\nabla^\mathrm{coord}+\mathrm d_\mathrm H,[\ ,\ ]\right),\\
  \left(\Omega_{C^{\mathrm{coord},L}\widehat{\otimes}F/C^{\mathrm{coord},L}},\mathrm d_\omega,\mathrm L\right)&\cong \left(\Omega_{C^{\mathrm{coord},L}\widehat{\otimes}_{R_1}J\!L/C^{\mathrm{coord},L}},{}^1\nabla^\mathrm{coord},\mathrm L\right),\\
  \left(\widehat{C}_{C^{\mathrm{coord},L},\bullet}(C^{\mathrm{coord},L}\widehat{\otimes}F),\mathrm
    d_\omega+\mathrm b_\mathrm H,\mathrm L\right)&\cong
  \left(\widehat{C}_{C^{\mathrm{coord},L},\bullet}(C^{\mathrm{coord},L}\widehat{\otimes}_{R_1}J\!L),{}^1\nabla^\mathrm{coord}+\mathrm
    b_\mathrm H,\mathrm L\right),
\end{align*}
an $L_\infty$-morphism
\[
\mathcal U_L^\mathrm{coord}:\left(T_{\mathrm{poly},C^{\mathrm{coord},L}}(C^{\mathrm{coord},L}\widehat{\otimes}_{R_1}J\!L),{}^1\nabla^\mathrm{coord},[\ ,\ ]\right)\to \left(D_{\mathrm{poly},C^{\mathrm{coord},L}}(C^{\mathrm{coord},L}\widehat{\otimes}_{R_1}J\!L),{}^1\nabla^\mathrm{coord}+\mathrm d_\mathrm H,[\ ,\ ]\right),
\]
which yields an $L_\infty$-module structure on $\widehat{C}_{C^{\mathrm{coord},L},\bullet}(C^{\mathrm{coord},L}\widehat{\otimes}_{R_1}J\!L)$ over $T_{\mathrm{poly},C^{\mathrm{coord},L}}(C^{\mathrm{coord},L}\widehat{\otimes}_{R_1}J\!L)$, and finally an $L_\infty$-morphism 
\[
\mathcal
S_L^\mathrm{coord}:\left(\widehat{C}_{C^{\mathrm{coord},L},\bullet}(C^{\mathrm{coord},L}\widehat{\otimes}_{R_1}J\!L),{}^1\nabla^\mathrm{coord}+\mathrm
  b_\mathrm H,\mathrm L\circ \mathcal U_L^\mathrm{coord}\right)\to
\left(\Omega_{C^{\mathrm{coord},L}\widehat{\otimes}F/C^{\mathrm{coord},L}},{}^1\nabla^\mathrm{coord},\mathrm
  L\right).
\]
We recall from \S\ref{ref-4.1-54} that there is a rational action of
$\mathrm{GL}_d(k)$ on $C^{\mathrm{coord},L}$, which extends in a
natural way to a (topological) rational action on all DG-Lie algebras and
DG-Lie modules above.  The previous actions determine infinitesimally
actions of $\mathfrak{gl}_d(k)$ on all DG-Lie algebras and DG-Lie
modules considered so far in the sense of \S\ref{ref-5.2-76}.

The $L_\infty$-morphism $\mathcal U_{L,\omega}$ descends with respect
to the action of the set $\mathfrak s=\mathfrak{gl}_d(k)$ (using the
notation of \S\ref{ref-5.2-76}), because the descent condition
\eqref{ref-5.11-83} is satisfied as a consequence of Property $iii)$ of
Theorem~\ref{ref-5.1-85} and of the verticality property
\eqref{ref-4.5-59} of $\omega$.

Similarly, Property $ii)$ of Theorem~\ref{ref-5.2-86}, together with the
verticality property of $\omega$ implies that $\mathcal S_{L,\omega}$
descends with respect to the action of $\mathfrak{gl}_d(k)$ (see
\S\ref{ref-5.3-81}).  Summarizing all arguments so far, and because
of the compatibility of the $\mathrm{GL}_d(k)$-action with the
isomorphism~\eqref{ref-4.2-56}, we get $L_\infty$-morphisms
\begin{multline*}
\mathcal
(U_L^\mathrm{coord})^{\mathfrak{gl}_d(k)}:\left(T_{\mathrm{poly},C^{\mathrm{coord},L}}(C^{\mathrm{coord},L}\widehat{\otimes}_{R_1}J\!L),{}^1\nabla^\mathrm{coord},[\
  ,\ ]\right)^{\mathfrak{gl}_d(k)}\\\to
\left(D_{\mathrm{poly},C^{\mathrm{coord},L}}(C^{\mathrm{coord},L}\widehat{\otimes}_{R_1}J\!L),{}^1\nabla^\mathrm{coord}+\mathrm
  d_\mathrm H,[\ ,\ ]\right)^{\mathfrak{gl}_d(k)}
\end{multline*}
and
\begin{multline*}
\mathcal
(S_L^\mathrm{coord})^{\mathfrak{gl}_d(k)}
:\left(\widehat{C}_{C^{\mathrm{coord},L},\bullet}(C^{\mathrm{coord},L}\widehat{\otimes}_{R_1}J\!L),{}^1\nabla^\mathrm{coord}+\mathrm
  b_\mathrm H,\mathrm L\circ \mathcal U_L^\mathrm{coord}\right)^{\mathfrak{gl}_d(k)}\\\to
\left(\Omega_{C^{\mathrm{coord},L}\widehat{\otimes}F/C^{\mathrm{coord},L}},{}^1\nabla^\mathrm{coord},\mathrm
  L\right)^{\mathfrak{gl}_d(k)},
\end{multline*}
Repeating almost {\em verbatim} the arguments at the end of~\cite[\S7.3.3]{CVdB}, there are obvious isomorphisms of DG-Lie algebras  and DG-Lie modules
\[
\begin{aligned}
  \left(T_{\mathrm{poly},C^{\mathrm{aff},L}}(C^{\mathrm{aff},L}\widehat{\otimes}_{R_1}J\!L),{}^1\nabla^\mathrm{aff},[\ ,\ ]\right)&\cong 
\left(T_{\mathrm{poly},C^{\mathrm{coord},L}}(C^{\mathrm{coord},L}\widehat{\otimes}_{R_1}J\!L),{}^1\nabla^\mathrm{coord},[\ ,\ ]\right)^{\mathfrak{gl}_d(k)},\\
  \left(D_{\mathrm{poly},C^{\mathrm{aff},L}}(C^{\mathrm{aff},L}\widehat{\otimes}_{R_1}J\!L),{}^1\nabla^\mathrm{aff}+\mathrm d_\mathrm H,[\ ,\ ]\right)&\cong 
\left(D_{\mathrm{poly},C^{\mathrm{coord},L}}(C^{\mathrm{coord},L}\widehat{\otimes}_{R_1}J\!L),{}^1\nabla^\mathrm{coord}+\mathrm d_\mathrm H,[\ ,\ ]\right)^{\mathfrak{gl}_d(k)},\\
  \left(\Omega_{C^{\mathrm{aff},L}\widehat{\otimes}F/C^{\mathrm{aff},L}},{}^1\nabla^\mathrm{aff},\mathrm L\right)&\cong 
\left(\Omega_{C^{\mathrm{coord},L}\widehat{\otimes}F/C^{\mathrm{coord},L}},{}^1\nabla^\mathrm{coord},\mathrm L\right)^{\mathfrak{gl}_d(k)},\\
  \left(\widehat{C}_{C^{\mathrm{aff},L},\bullet}(C^{\mathrm{aff},L}\widehat{\otimes}_{R_1}J\!L),{}^1\nabla^\mathrm{aff}+\mathrm
    b_\mathrm H,\mathrm L\right)&\cong
  \left(\widehat{C}_{C^{\mathrm{coord},L},\bullet}(C^{\mathrm{coord},L}\widehat{\otimes}_{R_1}J\!L),{}^1\nabla^\mathrm{coord}+\mathrm
    b_\mathrm H,\mathrm L\right)^{\mathfrak{gl}_d(k)}.
\end{aligned}
\]
We now set
\[
\begin{aligned}
\mathfrak g_1^L&=T_{\mathrm{poly},C^{\mathrm{aff},L}}(C^{\mathrm{aff},L}\widehat{\otimes}_{R_1}J\!L),\ & \ \mathfrak g_2^L&=D_{\mathrm{poly},C^{\mathrm{aff},L}}(C^{\mathrm{aff},L}\widehat{\otimes}_{R_1}J\!L),\\
\mathfrak m_1^L&=\Omega_{C^{\mathrm{aff},L}\widehat{\otimes}_{R_1}J\!L/C^{\mathrm{aff},L}},\ & \ \mathfrak m_2^L&=\widehat{C}_{C^{\mathrm{aff},L},\bullet}(C^{\mathrm{aff},L}\widehat{\otimes}_{R_1}J\!L),
\end{aligned}
\]
and $\mathfrak U_L=\mathcal U_L^\mathrm{aff}$, $\mathfrak S_L=\mathcal
S_L^\mathrm{aff}$: combining all the results so far, we get the
commutative diagram~\eqref{ref-5.12-89}, and to prove the claim, it
remains to show that $\mathfrak U_L$ and $\mathfrak S_L$ are
$L_\infty$-quasi-isomorphisms.

The proof of the fact that $\mathfrak U_L$ is a quasi-isomorphism can
be found in~\cite[\S 7.3.4]{CVdB}; the
proof of the fact that $\mathfrak S_L$ is a quasi-isomorphism is dual. We
will now sketch it. 

The $L_\infty$-morphism $\mathfrak S_L$ is obtained from $\mathcal
S_{L,\omega}$ using the
isomorphism~\eqref{ref-4.2-56} and by~\eqref{ref-5.9-80}
the Taylor components of $\mathcal S_{L,\omega}$ are given by
\[
\begin{aligned}
  &\mathcal S_{L,\omega,n}(\gamma_1,\dots,\gamma_n;c)=\sum_{m\geq 0}\frac{1}{m!}\mathcal S_{L,n+m}(\underset{m}{\underbrace{\omega,\dots,\omega}},\gamma_1,\dots,\gamma_n;c),\\
  &\gamma_i\in
  T_{\mathrm{poly},C^{\mathrm{coord},L}}(C^{\mathrm{coord},L}\widehat{\otimes}F),\
  c\in \widehat{C}_{C^{\mathrm{coord},L},\bullet}(C^{\mathrm{coord},L}\widehat\otimes
  F).
\end{aligned}
\]
$T_{\mathrm{poly},C^{\mathrm{coord},L}}(C^{\mathrm{coord},L}\widehat{\otimes}F)$,
$D_{\mathrm{poly},C^{\mathrm{coord},L}}(C^{\mathrm{coord},L}\widehat{\otimes}F)$,
$\Omega_{C^{\mathrm{coord},L}\widehat{\otimes}F/C^{\mathrm{coord},L}}$
and
$\widehat{C}_{C^{\mathrm{coord},L},\bullet}(C^{\mathrm{coord},L}\widehat\otimes
F)$ are bi-graded complexes: the first degree is the natural degree
coming from $C^{\mathrm{coord},L}$, while the second degree is
associated to poly-vector degree, (shifted) Hochschild degree,
(negative) form degree and (negative) Hochschild degree respectively.

The component $\mathcal S_{L,\omega,0}$ can be written into a sum
\[
\mathcal S_{L,\omega,0}(c)=\sum_{n\geq 0}\frac{1}{n!}\mathcal S_n(\underset{n}{\underbrace{\omega,\dots,\omega}};c);
\]
the grading property of the $L_\infty$-quasi-isomorphism $\mathcal S$
of Theorem~\ref{ref-5.2-86} implies that the
component $\mathcal S_{L,\omega,0}^n$ of $\mathcal S_{L,\omega,0}$
indexed by $n$ has bi-degree $(n,-n)$.

Dualizing~\cite[Lemma 7.3.2]{CVdB}, and using Property $i)$ of
Theorem~\ref{ref-5.2-86}, we get the following
commutative diagram of graded vector spaces:
\begin{equation}\label{ref-5.13-90}
\xymatrix{C_\mathrm{poly}^L(R)\ar@{^{(}->}[r]\ar[d]_{\mathrm{HKR}} & C^{\mathrm{coord},L}\widehat{\otimes} \widehat{C}_\bullet(F)\ar[d]^{\mathcal S_{L,\omega,0}^0}\\
\Omega^L(R)\ar@{^{(}->}[r] & C^{\mathrm{coord},L}\widehat{\otimes} \Omega_F},
\end{equation}
where the morphism $\mathrm{HKR}$ on the left vertical arrow has been defined in Theorem~\ref{ref-3.9-52}.

The twisting procedure and the descent procedure by the isomorphism~\eqref{ref-4.2-56} produce the commutative diagram 
\[
\xymatrix{C_\mathrm{poly}^L(R)\ar@{^{(}->}[r]\ar[d]_{\mathrm{HKR}} & \widehat{C}_{C^{\mathrm{aff},L},\bullet}(C^{\mathrm{aff},L}\widehat\otimes_{R_1}J\!L)\ar[d]^{\mathfrak S_{L,0}^0}\\
\Omega^L(R)\ar@{^{(}->}[r]& \Omega_{C^{\mathrm{aff},L}\widehat\otimes_{R_1}J\!L/C^{\mathrm{aff},L}}}
\] 
out of the commutative diagram~\eqref{ref-5.13-90}; the above bi-gradings
naturally translate into bi-gradings on
$\Omega_{C^{\mathrm{aff},L}\widehat{\otimes}F/C^{\mathrm{aff},L}}$ and
$\widehat{C}_{C^{\mathrm{aff},L},\bullet}(C^{\mathrm{aff},L}\widehat\otimes F)$.
The component $\mathfrak S_{L,0}$ is a sum of terms $\mathfrak S_{L,0}^n$,
$n\geq 0$, of bi-degree $(n,-n)$.

We now prove that the morphisms $\mathfrak S_{L,0}$ and $\mathfrak
S_{L,0}^0$ coincide at the level of cohomology.  For this, we consider
on the double complexes
$\Omega_{C^{\mathrm{aff},L}\widehat{\otimes}F/C^{\mathrm{aff},L}}$ and
$\widehat{C}_{C^{\mathrm{aff},L},\bullet}(C^{\mathrm{aff},L}\widehat\otimes F)$ the
filtration with respect to the second degree: then, the corresponding spectral
sequences degenerate at their first terms, because of the results of
\S\ref{ref-4.3-65},\ \S\ref{ref-4.4-67}, and the resulting complexes
consist of single columns $(\Omega^L(R),0)$ and
$(C_\mathrm{poly}^L(R),\mathrm b_\mathrm H)$. 
Thus, the respective
second terms of the spectral sequences coincide with $\Omega^L(R)$ and
with $\mathrm H^\bullet(C_\mathrm{poly}^L(R),\mathrm b_\mathrm H)$.
Since both spectral sequences degenerate at their first term (i.e.\
the cohomology with respect to the first degree is concentrated in degree
$0$), $\mathfrak S_{L,0}$ and $\mathfrak S_{L,0}^0$ obviously coincide at the
level of cohomology, and this ends the proof.
\end{proof}

\subsection{Functoriality property of Theorem~\ref{ref-5.3-88}}\label{ref-5.6-91}
We consider two Lie algebroids $(L,R)$, $(M,S)$ as above. 
\begin{Def}\label{ref-5.4-92}
An {\bf algebraic morphism} from $(L,R)$ to $(M,S)$ consists of a pair $(\ell,\lambda)$, where $i)$ $\lambda$ is a $k$-algebra morphism from $R$ to $S$, and $ii)$ $\ell$ is a Lie algebra morphism from $L$ to $M$, enjoying the following compatibility properties with respect to the corresponding anchor maps:
\[
\lambda(l(r))=\ell(l)(\lambda(r)),\ \ell(rl)=\lambda(r)\ell(l),\ r\in R,\ l\in L.
\]
\end{Def}
The universal property of the universal enveloping algebra of a Lie
algebroid yields, for any algebraic morphism $\varphi=(\ell,\lambda)$ from
$(L,R)$ to $(M,S)$, a Hopf algebroid morphism $\varphi_D:\mathrm U_R(L)\to
\mathrm U_S(M)$.  Thus, $(\ell,\lambda)$ defines a morphism $\varphi_D$
of $B_\infty$-algebras from $D_\mathrm{poly}^L(R)$ to
$D_\mathrm{poly}^M(S)$: in particular, it restricts to a morphism of
Gerstenhaber algebras up to homotopy.

Further, the algebraic morphism $\varphi$ defines a morphism
$\varphi_T:T_\mathrm{poly}^L(R)\to T_\mathrm{poly}^M(S)$ by
extending (via the $S$-linear wedge product) the assignment
\[
\varphi_T:S\otimes_R L\rightarrow M: s\otimes_R l\mapsto s\ell(l)
\] 
Since $(\ell,\lambda)$ preserves the anchor map and Lie bracket, we have a morphism of Gerstenhaber algebras from $T_\mathrm{poly}^L(R)$ to $T_\mathrm{poly}^M(S)$.
\begin{Prop}\label{ref-5.5-93}
  We assume $(L,R)$, $(M,S)$ to be Lie algebroids over $R$ and $S$
  respectively, and $\varphi=(\ell,\lambda)$ to be an algebraic morphism between
  them as in Definition~\ref{ref-5.4-92}; we further assume that the
  morphism
\[
\varphi_T:S\otimes_R L\rightarrow M: s\otimes_R l\mapsto s\ell(l)
\]
is an isomorphism of $S$-modules.

The morphism $(\ell,\lambda)$ determines a morphism of DG-algebras
\[
\left(\Omega^L(R),\mathrm d_L\right)\overset{\varphi_\Omega}\to \left(\Omega^M(S),\mathrm d_M\right),
\]
which satisfies
\begin{equation}\label{ref-5.14-94}
  \varphi_\Omega(\gamma\cap \omega)=\varphi_T(\gamma)\cap \varphi_\Omega(\omega),\ \gamma\in T_\mathrm{poly}^L(R),\ \omega\in \Omega^L(R),
\end{equation}
and a morphism of DG-algebras 
\[
\varphi_J:J\!L\to J\!M,
\]
which satisfies
\begin{align}
\label{ref-5.15-95}\lambda(\alpha(E))&=\varphi_J(\alpha)(\varphi_D(E)),\ \alpha\in J\!L,\ E\in \mathrm U_R(L),\\ 
\label{ref-5.16-96}\varphi_J\!\left({}^1\nabla_l\alpha\right)&=
{}^1\nabla_{\ell(l)}\varphi_J(\alpha),\ \alpha\in J\!L,\ l\in L,\\
\label{ref-5.17-97}\varphi_J({}^2\nabla_l\alpha)&={}^2\nabla_{\ell(l)}\varphi_J(\alpha),\ \alpha\in J\!L,\ l\in L,
\end{align}
and which commutes with the algebra monomorphisms $\alpha_i$, $i=1,2$ (see \S\ref{ref-3.1.4-29}).
\end{Prop}
\begin{proof}
  Since $\varphi_T$ is an isomorphism of $S$-modules, we define
  $\varphi_\Omega$ on $L$-differential forms on $R$ via
\[
\varphi_\Omega(r)=\lambda(r),\qquad \varphi_\Omega(l^*)(s\ell(l))=s\lambda(l^*(l)),\ r\in R,\ s\in S,\ l\in L,\ l^*\in L^*,
\]
and we extend it to $\Omega^R(L)$ by $R$-linearity and by multiplicativity with respect to the wedge product.

To prove that $\varphi_\Omega$ intertwines $\mathrm d_L$ and $\mathrm d_M$, it suffices to verify the claim on $R$ and $L^*$.
In the first case, we have
\[
\varphi_\Omega(\mathrm d_L(r))(s\ell(l))=s\lambda(\mathrm d_L(r)(l))=s\lambda(l(r))=s\ell(l)(\lambda(r))=s\ell(l)(\varphi_\Omega(r))=\mathrm d_M(\varphi_\Omega(r))(s\ell(l)),
\]
for a general element $r$ of $R$, $s$ of $S$ and $l$ of $L$, while in the second case we have
\[
\begin{aligned}
\varphi_\Omega(\mathrm d_L l^*)\!(s_1\ell(l_1),s_2\ell(l_2))
&=s_1s_2\lambda(\mathrm d_L l^*(l_1,l_2))\\
&=s_1s_2\lambda(l_1(l^*(l_2)))-s_1s_2\lambda(l_2(l^*(l_1)))-s_1s_2\lambda(l^*([l_1,l_2]))\\
&=s_1 \ell(l_1)(s_2)\lambda(l^*(l_2))+s_1s_2\ell(l_1)(\lambda(l^*(l_2)))-\\
&\phantom{=}-s_2\ell(l_2)(s_1)\lambda(l^*(l_1))-s_1s_2\ell(l_2)(\lambda(l^*(l_1)))-\\
&\phantom{=}-s_1\ell(l_1)(s_2)\lambda(l^*(l_2))+s_2\ell(l_2)(s_1)\lambda(l^*(l_1))-s_1s_2\lambda(l^*([l_1,l_2]))\\
&=(s_1\ell(l_1))(\varphi_\Omega(l^*)(s_2\ell(l_2)))-(s_2\ell(l_2))(\varphi_\Omega(l^*)(s_1\ell(l_1)))-\\
&\phantom{=}-\varphi_\Omega(l^*)\!\left(s_1\ell(l_1)(s_2)\ell(l_2)-s_2\ell(l_2)(s_1)\ell(l_1)+s_1s_2\ell([l_1,l_2])\right)\\
&=\mathrm d_M(\varphi_\Omega(l^*))(s_1\ell(l_1),s_2\ell(l_2)).
\end{aligned}
\]
By compatibility with wedge products, it suffices to prove
\eqref{ref-5.14-94} for $\gamma$ in $R$ or in $L$, and for a general
$\omega$: we check exemplarily the claim for $\gamma$ in $L$, i.e.\
\[
\begin{aligned}
\varphi_\Omega(l\cap\omega)(s_1\ell(l_1),\dots,s_p\ell(l_p))&=s_1\cdots s_p\lambda((l\cap\omega)(l_1,\dots,l_p))\\
&=s_1\cdots s_p \lambda(\omega(l,l_1,\dots,l_p))\\&=\varphi_\Omega(\omega)(\ell(l),s_1\ell(l_1),\dots,s_p\ell(l_p))\\
&=(\varphi_T(l)\cap \varphi_\Omega(\omega))(s_1\ell(l_1),\dots,s_p\ell(l_p)).
\end{aligned}
\]
We now define the morphism $\varphi_J$ on $J\!L$: for a general element
$\alpha$ of $J\!L$, we set
\[
\varphi_J(\alpha)(s)=s\lambda(\alpha(1)),\
\varphi_J(\alpha)(s\ell(l_1)\cdots\ell(l_p))=s\lambda(\alpha(l_1\cdots
l_p)),\ s\in S,\ l_i\in L.
\]
It is sufficient to define $\varphi_J$ on such elements of $\mathrm
U_S(M)$, since, being $\varphi_T$ an isomorphism of $S$-modules, a general
element of $\mathrm U_S(M)$ is a sum of elements of the form
\begin{multline*}
(s_1\ell(l_1))\cdots(s_p\ell(l_p))=s_1\ell(l_1)s_2\ell(l_2)\cdots s_p\ell(l_p)\\
=s_1(\ell(l_1)(s_2))\ell(l_2)\cdots s_p\ell(l_p)+s_1 s_2\ell(l_1)\ell(l_2)\cdots s_p\ell(l_p)=\cdots,
\end{multline*}
where the product has to be understood in $\mathrm U_S(M)$.  

Since $\varphi_D$ is defined by extending $\lambda$ and $\ell$ in a way
compatible with the Lie algebroid structure of $\mathrm U_R(L)$,
\eqref{ref-5.15-95} follows immediately.

As for \eqref{ref-5.16-96}, it suffices to check the identity on $R$
and on elements of $\mathrm U_S(M)$ of the form
$s\ell(l_1)\cdots\ell(l_p)$.  In the first case, we have for $s\in S,\ l\in L$
\[
\begin{aligned}
  \left({}^1\nabla_{\ell(l)}\varphi_J(\alpha)\right)(s)&=\ell(l)\!\left(\varphi_J(\alpha)(s)\right)-\varphi_J(\alpha)\!\left(\ell(l)s\right)\\
  &=\ell(l)\!\left(s\lambda(\alpha(1))\right)-\varphi_J(\alpha)(\ell(l)s)\\
  &=\ell(l)(s)\lambda(\alpha(1))+s\ell(l)\lambda(\alpha(1))
-\varphi_J(\alpha)(\ell(l)(s))-\varphi_J(\alpha)(s\ell(l))\\
  &=s\ell(l)\lambda(\alpha(1))-s\lambda(\alpha(l))\\&=s\lambda(({}^1\nabla_l\alpha)(1))\\
  &=\varphi_J({}^1\nabla_l\alpha)(s).
\end{aligned}
\]
As for the second case, we have for $\alpha\in J\!L,\ l, l_i\in L,\ i=1,\dots,p,\ s\in S$ 
\[
\begin{aligned}
  \left({}^1\nabla_{\ell(l)}\varphi_J(\alpha)\right)(s\ell(l_1)\cdots\ell(l_p))&=\ell(l)\left(\varphi_J(\alpha)(s\ell(l_1)\cdots\ell(l_p))\right)-\varphi_J(\alpha)\!\left(\ell(l)s\ell(l_1)\cdots\ell(l_p)\right)\\
  &=\ell(l)\!\left(s\lambda(\alpha(l_1\cdots l_p))\right)-\varphi_J(\alpha)\!\left(\ell(l)s\ell(l_1)\cdots\ell(l_p)\right)\\
  &=\ell(l)(s)\lambda(\alpha(l_1\cdots l_p))+s\ell(l)\!\left(\lambda(\alpha(l_1\cdots l_p))\right)\\
  &\qquad-\varphi_J(\alpha)(\ell(l)(s)\ell(l_1)\cdots\ell(l_p))-\varphi_J(\alpha)(s\ell(l)\ell(l_1)\cdots\ell(l_p))\\
  &=s\lambda(l(\alpha(l_1\cdots l_p)))-s\lambda(\alpha)(l l_1\cdots l_p)\\
  &=\varphi_J({}^1\nabla_l\alpha)(s\ell(l_1)\cdots\ell(l_p)).
\end{aligned}
\]
The identity~\eqref{ref-5.17-97} as well as the compatibility with
$\alpha_i$, $i=1,2$ are verified by similar computations.
\end{proof}
Assume now that $\varphi=(\ell,\lambda):(L,R)\rightarrow (M,S)$ is as in the
previous lemma and that $\varphi_T:S\otimes_R L\rightarrow M$ is an 
isomorphism. As always we assume that $L$ (and hence $M$) is free of rank $d$.
Looking at associated graded objects we see that the extended map
\begin{equation}
\label{ref-5.18-98}
S_1\otimes_{R_1} J\!L\rightarrow J\!M:s\otimes \alpha\mapsto s\varphi_J(\alpha)
\end{equation}
is an isomorphism. Hence any $R_1$-linear differential operator
on $J\!L$ can be extended to an $S_1$-linear differential operator
on $J\!M$. We use this to define a
map
\[
\varphi_D:D_{R_1}(J\!L)\rightarrow D_{S_1}(J\!M)
\]
and a corresponding map of $B_\infty$-algebras
\begin{equation}
\label{ref-5.19-99}
\varphi_D:D_{\poly,R_1}(J\!L)\rightarrow D_{\poly,S_1}(J\!M)
\end{equation}
such that the following diagram is commutative
\begin{equation}
\label{ref-5.20-100}
\xymatrix{
D_{\poly}^L(R) \ar[r]^{\varphi_D}\ar@{_(->}[d] &D_{\poly}^M(S)\ar@{_(->}[d] \\
D_{\poly,R_1}(J\!L)\ar[r]_{\varphi_D}& D_{\poly,S_1}(J\!M)
}
\end{equation}
where the vertical monomorphisms have been defined in \eqref{ref-3.12-37}. 

An easy computation shows that $\varphi_D$ in \eqref{ref-3.12-37}
commutes with the action of the Grothendieck connection
$[{}^1\nabla,-]$. It follows by the discussion in \S\ref{ref-3.3-33} that if we take the invariants for
$[{}^1\nabla,-]$ of the lower line in \eqref{ref-5.20-100} we obtain the
upper line.

\medskip

We extend $\varphi_J$ to a map of graded vector spaces
\[
\varphi_C:\widehat{C}_{R,\bullet}(J\!L)\rightarrow \widehat{C}_{S,\bullet}(J\!M):
\alpha_1\otimes\cdots \otimes \alpha_n\mapsto
\varphi_J(\alpha_1)\otimes\cdots\otimes  \varphi_J(\alpha_n)
\]
which is again essentially just base extension over $S/R$. 
This map obviously commutes with the Grothendieck connection ${}^1\nabla$.
We obtain a map of pairs of graded vector spaces
\[
(\varphi_D,\varphi_C):(D_{\poly,R}(J\!L),\widehat{C}_{R,\bullet}(J\!L))
\rightarrow (D_{\poly,S}(J\!M),\widehat{C}_{R,\bullet}(J\!L))
\]
and as this map is just base extension over $S/R$ it is compatible
with all structures defined in \cite{CR2}, hence in particular
with the DG-Lie algebra and DG-Lie module structures and also with the
precalculi up to homotopy.

Taking invariants for ${}^1\nabla$ and using \eqref{ref-3.16-45}
we obtain a commutative diagram of precalculus structure up to homotopy
\[
\xymatrix{(D_\mathrm{poly}^L(R),\mathrm d_\mathrm H,[\ ,\ ],\cup)\ar@{~>}[d]\ar[r]^-{\varphi_D} & (D_\mathrm{poly}^M(S),\mathrm d_\mathrm H,[\ ,\ ],\cup)\ar@{~>}[d]\\
(C_\mathrm{poly}^L(R),\mathrm b_\mathrm H,\mathrm L,\cap)\ar[r]_-{\varphi_C} & (C_\mathrm{poly}^M(S),\mathrm b_\mathrm H,\mathrm L,\cap)}.
\]
One also obtains from Proposition~\ref{ref-5.5-93} 
a commutative diagram of precalculi.
\[
\xymatrix{(T_\mathrm{poly}^L(R),0,[\ ,\ ],\cup)\ar@{~>}[d]\ar[r]^-{\varphi_T} & (T_\mathrm{poly}^M(S),0,[\ ,\ ],\cup)\ar@{~>}[d]\\
(\Omega^L(R),0,\mathrm L,\cap)\ar[r]_-{\varphi_\Omega} & (\Omega^M(S),0,\mathrm L,\cap)}.
\]
Furthermore from \eqref{ref-5.18-98} and the universal
property of coordinate spaces (see \eqref{ref-4.1-55}) we obtain
an 
$R$-algebra morphism from $R^{\mathrm{coord},L}$ to
$S^{\mathrm{coord},M}$.  It extends further to a morphism of
DG-algebras from $C^{\mathrm{coord},L}$ to $C^{\mathrm{coord},M}$ thanks to~\eqref{ref-5.16-96} and the fact that $\varphi_\Omega$ is a morphism
of DG-algebras from $\Omega^L(R)$ to $\Omega^M(S)$.

Finally, the algebraic morphism $(\ell,\lambda)$ induces precalculi
morphisms (up to homotopy) between all corresponding Fedosov
resolutions, since the monomorphism $\alpha_2$ and the connection
${}^2\nabla$, which are needed in the construction of the Fedosov
resolutions of \S\ref{ref-4-53} (we refer to~\cite{CVdB} for more
details thereabout), have been proved to be preserved by
$(\ell,\lambda)$.

As a consequence of these arguments, we deduce the following theorem,
which expresses the functoriality properties of the commutative
diagram~\eqref{ref-5.12-89} of Theorem~\ref{ref-5.3-88}.
\begin{Thm}\label{ref-5.6-101}
  For a general algebraic morphism $\varphi=(\ell,\lambda)$ from
  $(L,R)$ to $(M,S)$ as in Definition~\ref{ref-5.4-92}, which induces
  an isomorphism $S\otimes_R L\cong M$ of $S$-modules, and such that
  $L$ is free of rank $d$, 
the
  $L_\infty$-quasi-isomorphisms $\mathfrak U_L$, $\mathfrak U_M$,
  $\mathfrak S_L$ and $\mathfrak S_M$ of DG-Lie algebras and DG-Lie
  modules fit into the commutative diagram
\begin{equation}\label{ref-5.21-102}
  \xymatrix{T_\mathrm{poly}^L(R) \ar@{^{(}->}[r]\ar[d]_{\varphi_T} & \mathfrak g_1^L \ar[r]^{\mathfrak U_L}\ar[d]_{\varphi_T}& \mathfrak g_2^L\ar[d]_{\varphi_D} & D_\mathrm{poly}^L(R)\ar@{_{(}->}[l]\ar[d]_{\varphi_D}\\
    T_\mathrm{poly}^M(S) \ar@{^{(}->}[r]\ar@{~>}[d]_{\mathrm L} & \mathfrak g_1^M \ar[r]^{\mathfrak U_M}\ar@{~>}[d]_{\mathrm L_1}& \mathfrak g_2^M\ar@{~>}[d]_{\mathrm L_2} & D_\mathrm{poly}^M(S)\ar@{_{(}->}[l]\ar@{~>}[d]_{\mathrm L}\\
    \Omega^M(S)\ar@{^{(}->}[r] & \mathfrak m_1^M & \mathfrak m_2^M\ar[l]_{\mathfrak S_M} & C_\mathrm{poly}^M(S)\ar@{_{(}->}[l]\\
    \Omega^L(R)\ar@{^{(}->}[r]\ar[u]^{\varphi_\Omega}& \mathfrak m_1^L \ar[u]^{\varphi_\Omega}& \mathfrak m_2^L\ar[l]_{\mathfrak S_L} \ar[u]^{\varphi_J}& C_\mathrm{poly}^L(R)\ar@{_{(}->}[l]\ar[u]^{\varphi_J}},
\end{equation}
where we have borrowed notation from Proposition~\ref{ref-5.5-93}; all such morphisms are compatible
with respect to the composition of algebraic morphisms between Lie algebroids.
\end{Thm}
Note that Theorem~\ref{ref-5.6-101} makes no reference to the
(homotopy) precalculus structures which we discussed above;  we will
need these below. 
\subsection{Proof of Theorem~\ref{ref-1.2-11}}\label{ref-5.7-103}
We now collect the results of \S\ref{ref-5.5-87}
and \S\ref{ref-5.6-91} to give the proof of Theorem~\ref{ref-1.2-11},
 {\em via} a well-suited gluing procedure.

 We consider a ringed site $(X,\mathcal O)$, and a sheaf of Lie
 algebroids $\mathcal L$, such that $\mathcal L$ is locally free of
 rank $d$ over $\mathcal O$. We replace $X$ by its full subcategory
of objects $U$
such that $\mathcal{L}(U)$ is free over $\mathcal{O}(U)$. This
does not change the category of sheaves. 

 All sheaves of DG-Lie algebras and DG-Lie
 modules in the commutative diagram~\eqref{ref-1.7-12} are obtained by
 sheafifying the corresponding presheaves of DG-Lie algebras and
 DG-Lie modules, i.e.\
\[
\begin{aligned}
U&\to T_\mathrm{poly}^{\mathcal L(U)}(\mathcal O(U)),\ &\ U&\to D_\mathrm{poly}^{\mathcal L(U)}(\mathcal O(U)),\\
U&\to \Omega^{\mathcal L(U)}(\mathcal O(U)),\ & \ U&\to C_\mathrm{poly}^{\mathcal L(U)}(\mathcal O(U)).
\end{aligned}
\]
Since $\mathcal L$ is locally free of order $d$ over $\mathcal O$, for
a morphism $V\rightarrow U$ in $X$, the
corresponding restriction morphism $(\mathcal O(U),\mathcal L(U))\to
(\mathcal O(V),\mathcal L(V))$, yields an isomorphism
\[
\mathcal O(V)\otimes_{\mathcal O(U)}\mathcal L(U)\cong \mathcal L(V).
\] 
Thus, any restriction morphism as above may be viewed as an algebraic morphism between Lie algebroids, satisfying the isomorphism property of Theorem~\ref{ref-5.6-101}.

If we then consider the DG-Lie algebras and DG-Lie modules
\[
U\to \mathfrak g_i^{\mathcal L(U)},\ U\to \mathfrak m_i^{\mathcal L(U)},\ i=1,2,
\]
Theorem~\ref{ref-5.3-88} produces, for any $U$ in $X$,
$L_\infty$-quasi-isomorphisms $\mathfrak U_{\mathcal L(U)}$ and
$\mathfrak S_{\mathcal L(U)}$ which
fit into a commutative diagram~\eqref{ref-5.12-89}. By Theorem \ref{ref-5.6-101}
these are actually morphisms of presheaves.

Sheafifying all presheaves and morphisms between presheaves concludes
the proof.

\section{The relationship between Atiyah classes and jet bundles}\label{ref-6-104}
In the present section we review some technical results
from~\cite[\S8]{CVdB}, to which we refer for more details.  We
need only the main notation and conventions for use in~\S\ref{ref-7-109}.

For a field $k$ of characteristic $0$, we consider a sheaf $\mathcal
L$ of Lie algebroids over a ringed site $(X,\mathcal O)$, which is
locally free of rank $d$ over $\mathcal O$.

We have a short exact sequence of $\mathcal{O}_1$-$\mathcal{O}_2$-bimodules
\[
0\rightarrow \mathcal{L}^\ast \rightarrow J^1\mathcal{L}\rightarrow
\mathcal{O}\rightarrow 0
\]
where $\mathcal O_i$, $i=1,2$, denotes a copy of $\mathcal O$ embedded
in $J\mathcal L$ {\em via} the monomorphism~$\alpha_i$ and where $J^1\mathcal{L}$
was introduced in \S\ref{ref-3.1.4-29}.

For a general $\mathcal O$-module $\mathcal E$, tensoring over $\mathcal{O}_2$
yields a short exact sequence
\[
\xymatrix{0\ar[r]& \mathcal L^*\otimes_{\mathcal O}\mathcal E \ar[r] & J^1\mathcal L\otimes_{\mathcal O_2}\mathcal E\ar[r] & \mathcal E \ar[r] & 0}
\]  
which we will call the $\mathcal L$-Atiyah sequence.  The
$\mathcal L$-Atiyah class $A_\mathcal L(\mathcal E)$ of $\mathcal E$
over~$\mathcal L$ is the extension class of this sequence in
$\mathrm{Ext}^1_{\mathcal O}\!\left(\mathcal E,\mathcal
  L^*\otimes_{\mathcal O}\mathcal E\right)$.  As explained in
\S\ref{ref-1.1-1}, if $\mathcal{E}$ is a vector bundle, the~$i$-th scalar
Atiyah class $a_{\mathcal L,i}(\mathcal E)$ of $\mathcal E$ is defined
as
\begin{equation}\label{ref-6.1-105}
a_{\mathcal L,i}(\mathcal E)=\mathrm{tr}(\bigwedge^i A_{\mathcal L}(\mathcal E))\in \mathrm H^i(X,\bigwedge^i \mathcal L^*),
\end{equation} 
Below we will only consider the case $\mathcal E=\mathcal L$. In that
case we simplify the notation~to
\[
A(\mathcal L)=A_{\mathcal L}(\mathcal L),\qquad a_i(\mathcal L)=a_{i,\mathcal L}(\mathcal L)
\]
Observe that the $a_i(\mathcal L)$ are cohomology classes. We now outline
how we may realize them as explicit cocycles.

By the very construction of $C^{\mathrm{coord},\mathcal L}$ and $C^{\mathrm{aff},\mathcal L}$, there are natural morphisms of DG-algebras
\begin{equation}
\label{ref-6.2-106}
\xymatrix{\Omega^{\mathcal L_2}(X)\ar@{^{(}->}[r]^-\theta & C^{\mathrm{aff},\mathcal L}\widehat\otimes_{\mathcal O_{1}}\Omega_{J\mathcal L/\mathcal O_{1}} \ar@{^(->}[r] & C^{\mathrm{coord},\mathcal L}\widehat\otimes_{\mathcal O_{1}}\Omega_{J\mathcal L/\mathcal O_{1}}\cong C^{\mathrm{coord},\mathcal L}\widehat\otimes\Omega_F}.
\end{equation}
The differentials on the first
three DG-algebras are the natural ones (see \S\ref{ref-4.1-54}). The
differential on the fourth DG-algebra is $\mathrm d+\mathrm L_\omega$ for a
certain MC element $\omega\in C^{\coord,\mathcal L}\widehat{\otimes}
\operatorname{Der}(F)$ and $\mathrm d$ the natural
differential. See again \S\ref{ref-4.1-54}.

The MC element $\omega$ can be expressed as
\[
\omega=\eta_\alpha\omega_{\alpha,i}\mathrm \partial_{x_i},\ i=1,\dots, d,
\]
where $\eta_\alpha$ is in $C^{\mathrm{coord},\mathcal L}$ and has degree $1$, $\omega_{\alpha,i}$ belongs to $F$ and $\partial_{x_i}=\partial/\partial x_i$. 

If we define $\Xi$ to be the matrix with entries
\begin{equation}
\label{ref-6.3-107}
\Xi_{ij}=\eta_\alpha\mathrm d_F(\partial_{x_j}\omega_{\alpha,i})\in  C^{\mathrm{coord},\mathcal L}\widehat\otimes\Omega_F,
\end{equation}
where $\mathrm d_F$ is the De Rham differential on $\Omega_F$, 
then on the
nose we have
\[
\mathrm{Tr}(\Xi^n)\in C^{\mathrm{coord},\mathcal L}\widehat\otimes\Omega_F
\]
Furthermore it is true that
\[
(\mathrm d+ \mathrm L_\omega)(\mathrm{Tr}(\Xi^n))=0
\]
It is shown in~\cite[\S8]{CVdB} that $\mathrm{Tr}(\Xi^n)$ is
actually the image of a (necessarily unique) element in
$C^{\mathrm{aff},\mathcal L}\widehat\otimes_{\mathcal
  O_{1}}\Omega_{J\mathcal L/\mathcal O_{1}}$. Abusing notation somewhat
we will still write this element as $\mathrm{Tr}(\Xi^n)$. It is still
a cocycle and in this way represents an element of
\[
\overline{\mathrm{Tr}(\Xi^n)}\in\Gamma(X,\mathrm H^{2n}(C^{\mathrm{aff},\mathcal L}\widehat\otimes_{\mathcal
  O_{1}}\Omega_{J\mathcal L/\mathcal O_{1}}))
\]
which maps naturally to the hypercohomology
\[
\mathbb{H}^{2n}(X,C^{\mathrm{aff},\mathcal L}\widehat\otimes_{\mathcal
  O_{1}}\Omega_{J\mathcal L/\mathcal O_{1}})
\]
Further, we observe that the injection $\Omega^\mathcal
L(X)\xrightarrow{\theta} C^{\mathrm{aff},\mathcal
  L}\widehat\otimes_{\mathcal O_{1}}\Omega_{J\mathcal L/\mathcal
  O_{1}}$ of DG-algebras is a  quasi-isomorphism, as discussed
in~\S\ref{ref-4.3-65}. Thus $\theta$ induces an isomorphism
\[
\bigoplus_{m,n} \mathrm H^m(X,\wedge^n \mathcal L^\ast)= \mathbb{H}^\bullet(X,\Omega^\mathcal
L(X))\xrightarrow{\mathbb{H}(\theta)}  \mathbb{H}^{\bullet}(X,C^{\mathrm{aff},\mathcal
  L}\widehat\otimes_{\mathcal O_{1}}\Omega_{J\mathcal L/\mathcal
  O_{1}})
\]
The following identity is \cite[eq.\ (8.8)]{CVdB}
\begin{equation}\label{ref-6.4-108}
a_{n}(\mathcal L)=\mathbb{H}(\theta)^{-1}\left(\overline{\mathrm{Tr}(\Xi^n)}\right),\ n\geq 1,
\end{equation}
which indeed expresses $a_n(\mathcal L)$ in terms of the 
explicit cocycle $\mathrm{Tr}(\Xi^n)$.
\section{Proof of Theorem~\ref{ref-1.1-9}}\label{ref-7-109}
The aim of this Section is to prove Theorem~\ref{ref-1.1-9} which
implies C{\u{a}}ld{\u{a}}raru's conjecture (Theorem~\ref{ref-1.3-14}) as has been outlined
in the introduction.

For this purpose, we first remind the reader of 
the main result of~\cite{CR2} about compatibility between cap
products. We then prove a ring-theoretical globalized version
of this result (compare to the proof of
Theorem~\ref{ref-5.3-88}). By
functoriality (see \S\ref{ref-5.6-91}), we obtain the sheaf-theoretical globalization.  Finally,
using results of~\cite{CR2},
we compute explicitly the isomorphism
appearing in the compatibility between cap products, which we
identify with the action of the homological
HKR-quasi-isomorphism followed by left multiplication by the square
root of the (modified) Todd class.

\subsection{A {\em memento} of compatibility between cup and cap products}\label{ref-7.1-110}
In this section, we present a {\em memento} of the main results
of~\cite{CR2,CVdB} concerning compatibility between cup and
cap products respectively.

First of all as before $F$ is the algebra of formal power series in
$d$ variables over the field $k$ which we assume to contain
$\mathbb{R}$ for now.  We recall the existence of (homotopy)
Gerstenhaber algebra structures on $T_\mathrm{poly}(F)$ and
$D_\mathrm{poly}(F)$, which together with $\Omega_F$ and
$\widehat{C}_\bullet(F)$ yield (homotopy) precalculi \cite{CR2}.

We recall also the $L_\infty$-quasi-isomorphisms $\mathcal U$ introduced in
Theorem~\ref{ref-5.1-85} and $\mathcal S$ introduced in
Theorem~\ref{ref-5.2-86}. We denote by $\mathcal
U_n$, $n\geq 1$, resp.\ $\mathcal S_n$, $n\geq 0$, the $n$-th Taylor
component of $\mathcal U$, resp.\ $\mathcal S$.

We further consider a commutative DG-algebra $(\mathfrak m,\mathrm
d_\mathfrak m)$. 
 The precalculus
structures on $(T_\mathrm{poly}(F),\Omega_F)$ and
$(D_\mathrm{poly}(F),\widehat{C}_\bullet(F))$, can be
extended by $\mathfrak m$-linearity to precalculi $(T_\mathrm{poly}^\mathfrak
m(F),\Omega_F^\mathfrak
m)
=(T_\mathrm{poly}(F)\widehat\otimes\mathfrak m,\Omega_F\widehat\otimes\mathfrak m)
$ and
$(D_\mathrm{poly}^\mathfrak m(F),\widehat C^\mathfrak m _\bullet
(F))=(D_{\mathrm{poly},\mathfrak
  m}(F\widehat\otimes \mathfrak m),\widehat{C}_{\bullet,\mathfrak m} (F\widehat\otimes \mathfrak m))
)$.

\medskip

\noindent
\textbf{Convention. }Below we will work with potentially infinite
series with coefficients in $\mathfrak{m}$. We make the standard
assumption that we are in a setting where all these series converge
and standard series manipulations are allowed. In our actual
application all series will be finite for degree reasons.

\medskip

A MC element $\gamma$ of $T_\mathrm{poly}^\mathfrak m(F)$ can be written as a sum
\[
\gamma=\gamma_{-1}+\gamma_0+\gamma_1+\gamma_2+\cdots,
\]
where $\gamma_i$ is an element of $T_\mathrm{poly}^\mathfrak m(F)$ of poly-vector degree $i$, $i\geq -1$, which satisfies the Maurer--Cartan equation
\[
\mathrm d_\mathfrak m \gamma+\frac{1}2[\gamma,\gamma]=0.
\]
We denote by $\mathcal U(\gamma)$ the image of a MC element $\gamma$
as above with respect to $\mathcal U$ (see \eqref{ref-5.7-78}).  This is
again a MC element.  Further, we set
\[
\begin{aligned}
\mathcal U_{\gamma,1}(\gamma_1)&=\sum_{n\geq 0}\frac{1}{n!}\mathcal U_{n+1}(\underset{n}{\underbrace{\gamma,\dots,\gamma}},\gamma_1),\ \gamma_1\in T_\mathrm{poly}^\mathfrak m(F),\\
\mathcal S_{\gamma,0}(c)&=\sum_{n\geq 0}\frac{1}{n!}\mathcal S_n(\underset{n}{\underbrace{\gamma,\dots,\gamma}};c),\ c\in \widehat{C}^\mathfrak m_\bullet(F).
\end{aligned}
\]
Since $\mathcal U$ and $\mathcal S_\gamma$ are
$L_\infty$-quasi-isomorphisms, $\mathcal U_{\gamma,1}$ and $\mathcal
S_{\gamma,0}$ are both quasi-isomorphisms of DG-vector spaces.
\begin{Thm}\label{ref-7.1-111}
  For a general commutative DG-algebra $(\mathfrak m,\mathrm
  d_\mathfrak m)$ as above, and for a general MC element $\gamma$ of
  $T_\mathrm{poly}^\mathfrak m(F)$, $\mathcal U_{\gamma,1}$ and
  $\mathcal S_{\gamma,0}$ descend to quasi-isomorphisms of (homotopy)
  precalculi, fitting into the commutative
  diagram
\[
\xymatrix{\left(T_\mathrm{poly}^\mathfrak m(F),\mathrm d_\mathfrak m+[\gamma,\bullet],[\ ,\ ],\cup\right)\ar@{~>}[d]\ar[r]^-{\mathcal U_{\gamma,1}} & \left(D_\mathrm{poly}^\mathfrak m(F),\mathrm d_\mathfrak m+\mathrm d_\mathrm H+[\mathcal U(\gamma),\bullet],[\ ,\ ],\cup\right)\ar@{~>}[d]\\
\left(\Omega^\mathfrak m_F,\mathrm d_\mathfrak m+\mathrm L_{\gamma},\mathrm L,\cap\right)  & \ar[l]_-{\mathcal S_{\gamma,0}}\left(\widehat{C}^\mathfrak m_\bullet(F),\mathrm d_\mathfrak m+\mathrm b_\mathrm H+\mathrm L_{\mathcal U(\gamma)},\mathrm L,\cap\right)},
\]
in the sense that $\mathcal U_{\gamma,1}$ and $\mathcal S_{\gamma,0}$ preserve Lie brackets, Lie actions, cup and cap products up to homotopy.
\end{Thm}
Kontsevich~\cite{K} has first stated and proved that $\mathcal
U_{\gamma,1}$ defines a quasi-isomorphism of Gerstenhaber algebras up
to homotopy from $T_\mathrm{poly}^\mathfrak m(F)$ to
$D_\mathrm{poly}^\mathfrak m(F)$ in the sense specified above.  We
observe that the identity $\mathcal
U_{\gamma,1}([\gamma_1,\gamma_2])=\left[\mathcal
  U_{\gamma,1}(\gamma_1),\mathcal U_{\gamma_1}(\gamma_2)\right]$ at
the level of cohomology, for $\gamma_i$ in $T_\mathrm{poly}^\mathfrak
m(F)$, $i=1,2$, holds true, because $\mathcal U$ is an
$L_\infty$-morphism: in particular, there is a homotopy operator
describing the compatibility with Lie brackets, expressible in terms of
the Taylor components of $\mathcal U$ twisted by the MC element
$\gamma$.  On the other hand, the identity $\mathcal
U_{\gamma,1}(\gamma_1\cup\gamma_2)=\mathcal
U_{\gamma,1}(\gamma_1)\cup\mathcal U_{\gamma,1}(\gamma_2)$ at the
level of cohomology comes from a more complicated identity up to
homotopy: in this situation, the homotopy operator is not expressible
in terms of the Taylor components of $\mathcal U$.  For an explicit
description of the homotopy operator, we refer to~\cite{MT,CVdB,CR2}.

The actual formulation of Theorem~\ref{ref-7.1-111} has been first proposed
as a conjecture in the particular case, where $\gamma$ is a (formal)
Poisson structure, by Shoikhet~\cite{Sh}: this conjecture has been first proved in~\cite{Sh1} only in degree $0$
and later in~\cite{CR1} for all degrees.  A more general result has been stated and
proved in~\cite{CR2}, to which we refer for more details. The
identity $\mathcal S_{\gamma,0}(\mathrm L_{\mathcal
  U_{\gamma,1}}(\gamma_1)(c))=\mathrm L_{\gamma_1}(\mathcal
S_{\gamma,0}(c))$ at the level of cohomology, for $\gamma_1$ in
$T_\mathrm{poly}^\mathfrak m(F)$, $c$ in $\widehat{C}^\mathfrak m_\bullet(F)$, is a
consequence of the fact that $\mathcal S_{L,\gamma}$ is an
$L_\infty$-morphism of $L_\infty$-modules (in particular, there is a
homotopy formula involving the Taylor components of $\mathcal U$ and
$\mathcal S$, twisted by $\gamma$).  The identity $\mathcal
S_{\gamma,0}(\mathcal U_{\gamma,1}(\gamma_1)\cap
c)=\gamma_1\cap\mathcal S_{\gamma,0}(c)$ at the level of cohomology
holds true in virtue of a homotopy formula, but the corresponding
homotopy operator does not involve the Taylor components of $\mathcal
U$ and $\mathcal S$: such an operator has been explicitly described
in~\cite{CR2}.

We briefly review in section \S\ref{ref-7.1.1-112} below the construction of the homotopy operator for the compatibility between cap products.

\subsubsection{The homotopy formula for the compatibility between cap products}\label{ref-7.1.1-112}
For later computations, we write down the explicit homotopy operator
for the compatibility between the $\cap$-actions: namely, for a MC
element $\gamma$ as in Theorem~\ref{ref-7.1-111}, for $\gamma_1$ a general
element of $T_\mathrm{poly}^\mathfrak m(F)$ and $c$ a general element
of $\widehat C^\mathfrak m_\bullet(F)$, we have the homotopy relation
\begin{equation}\label{ref-7.1-113}
\begin{aligned}
\mathcal S_{\gamma,0}(\mathcal U_{\gamma,1}(\gamma_1)\cap c)-\gamma_1\cap \mathcal S_{\gamma,0}(c)&=(\mathrm d_\mathfrak m+\mathrm L_\gamma)\mathcal H^\mathcal S_\gamma(\gamma_1,c)+\mathcal H^\mathcal S_\gamma(\mathrm d_\mathfrak m\gamma_1+[\gamma,\gamma_1],c)+\\
&\phantom{=}+(-1)^{|\gamma_1|}\mathcal H^\mathcal S_\gamma(\gamma_1,\mathrm d_\mathfrak m c+\mathrm b_\mathrm Hc+\mathrm L_{\mathcal U(\gamma)}c),  
\end{aligned}
\end{equation}
where
\begin{equation}\label{ref-7.2-114}
\mathcal H_\gamma^\mathcal S(\gamma_1,c)=\sum_{n\geq 0}\frac{1}{n!}\sum_{\Gamma\in \mathcal G_{n+1,m+1}^S}\overset{\circ}W_{D,\Gamma}\mathcal S_\Gamma(\gamma_1,\underset{n}{\underbrace{\gamma,\dots,\gamma}},c),
\end{equation}
$c$ being of Hochschild degree $-m$.

In~\eqref{ref-7.2-114}, the second sum is over ``$S$-admissible graphs''
of type $(n+1,m+1)$: these are directed graphs with $n+2$ vertices of
the first type and $m+1$ cyclically ordered vertices of the second type and with an
orientation of the outgoing edges from vertices of the first type, and
with a special vertex of the first type, labelled by $0$.  The
vertices of the second type can be only endpoints of edges, and
$S$-admissible graphs do not contain edges starting and ending at the
same vertex; finally, the vertex $0$ has only incoming edges.

To the vertex $1$ of the first type of an $S$-admissible graph $\Gamma$ is assigned the poly-vector field $\gamma_1$: the number of outgoing edges from $1$ equals the poly-vector degree of $\gamma_1$ plus $1$.
To any other vertex of the first type, except $0$, is assigned a copy of the MC element $\gamma$.
To the $i$-th vertex of the second type is assigned the $i+1$-th component of the Hochschild chain $c$. 
Pictorially, here is an $S$-admissible graph of type $(4,5)$, with corresponding coloring by poly-vector fields and Hochschild chains:
\bigskip
\begin{center}
\resizebox{0.36 \textwidth}{!}{\input{S_graph.pstex_t}}\\
\text{Figure 1 - An $S$-admissible graph of type $(4,5)$ } \\
\end{center}
\bigskip

The differential form $S_\Gamma(\gamma_1,\underset{n}{\underbrace{\gamma,\dots,\gamma}},c)$ is defined explicitly in~\cite{Sh,CR1,CR2}.

More important for our purposes is the integral weight $\overset{\circ}W_{D,\Gamma}$, for a general $S$-admissible graph of type $(n+1,m+1)$, 
\begin{equation}
\overset{\circ}W_{D,\Gamma}=\int_{\overset{\circ}{\mathcal Y}_{n+1,m+1}^+}\omega_{D,\Gamma}.
\end{equation}
First of all, $\overset{\circ}{\mathcal Y}_{n+1,m+1}^+$ denotes the codimension-$1$-submanifold (with corners) of the compactified configuration space $\mathcal D_{n+1,m+1}^+$ of $n+1$ points in the punctured unit disk $D^\times$ and $m+1$ cyclically oriented points in $S^1$, consisting of configurations of points, where the point labelled by $1$ moves on a smooth curve from the origin to the first point $\overline 1$ (with respect to the cyclic order) in $S^1$.
Graphically,
\bigskip
\begin{center}
\resizebox{0.30 \textwidth}{!}{\input{BNE_curconf.pstex_t}}\\
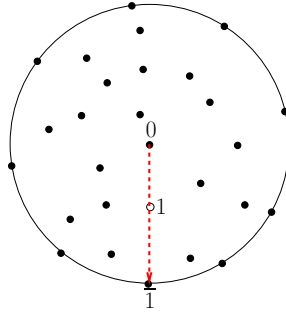
\text{Figure 2 - A general configuration of points in $\overset{\circ}{\mathcal Y}_{n+1,m+1}^+$} \\
\end{center}
\bigskip
In Figure 1, the dashed line represents the curve, along which the point $1$ (labelled as ``$\circ$'') moves.
The differential form $\omega_{D,\Gamma}$ associated to a graph in $\mathcal G^S_{n+1,m+1}$ is a product of smooth $1$-forms on $\mathcal D_{n+1,m+1}^+$: the basic ingredient is a slight modification of the exterior derivative of Kontsevich's angle function, see~\cite{K,CR2} for more details.

For the globalization procedure of the compatibility between cap
products, we need the following technical Lemma, which corresponds, in
the present framework, to 
Theorem~\ref{ref-5.2-86}, $ii)$.
\begin{Lem}\label{ref-7.2-115}
If $\Gamma$ is an $S$-admissible graph in $\mathcal G_{n+1,m+1}^S$, $n\geq 1$, and at least one of the poly-vector fields $\gamma_i$, $i\neq 1$, is linear on $F$, then
\[
\overset{\circ}W_{D,\Gamma}\mathcal S_\Gamma(\gamma_1,\gamma_2,\dots,\gamma_{n+1},c)=0.
\] 
\end{Lem}
\begin{proof}
  The first point of the first type in $\overset{\circ}{\mathcal
    Y}^+_{n+1,m+1}$, by the very construction of
  $\overset{\circ}{\mathcal Y}_{n+1,m+1}^+$, moves from the origin $0$
  to the first point in $S^1$ with respect to the cyclic order: to the
  former point is associated the poly-vector field $\gamma_1$.  Any
  other point associated to a vertex of the first type moves freely in
  the punctured unit disk $D^\times$.

  Without loss of generality we assume $\gamma_2$ to be an $\mathfrak
  m$-valued linear vector field: the valence (i.e.\ the number of outgoing edges) of the corresponding
  vertex of the first type is $1$, while the linearity of $\gamma_2$
  implies that there can be at most one incoming edge to the vertex
  corresponding to $\gamma_2$.  This follows from the construction of
  the differential form $\mathcal
  S_\Gamma(\gamma_1,\gamma_2,\dots,\gamma_{n+1},c)$.

Thus, we may safely restrict to $S$-admissible graphs $\Gamma$, such that the vertex $2$ has valence exactly $1$ and with at most one incoming edge.

If the vertex labelled by $2$ does not have incoming edges, the corresponding integral weight $\overset{\circ}{W}_{D,\Gamma}$ vanishes by dimensional reasons: in fact, we integrate a $1$-form (corresponding to the only outgoing edge from $2$) over a $2$-dimensional submanifold (with corners) of $D^\times$.

If the vertex labelled by $2$ has exactly one incoming and one
outgoing edge, we may apply~\cite[Lemma 6.1]{CR2}, to yield the
vanishing of the corresponding weight $\overset{\circ}W_{D,\Gamma}$.
\end{proof}

\subsection{The proof of Theorem~\ref{ref-1.1-9} in the ring case}\label{ref-7.2-116}
We will first assume that the ground field contains $\mathbb{R}$. 
At the end of the section we will show how to get rid of this restriction.

We consider a Lie algebroid $L$ over $R$ as in
Definition~\ref{ref-3.1-18} free of rank $d$
over $R$.  Then, we set $(\mathfrak m,\mathrm d_\mathfrak
m)=(C^{\mathrm{coord},L},\mathrm d)$, where $\mathrm d=\mathrm
d_{\Omega_{R^{\mathrm{coord},L}}}\widehat\otimes_{\Omega_{R_1}}1+1\widehat\otimes_{\Omega_{R_1}}\mathrm
d_{L_1}$ (see \S\ref{ref-4.1-54} for more details), and the
Maurer--Cartan form $\omega$ is an $\mathfrak m$-valued vector field
on $F$ obeying
\[
\mathrm d \omega+\frac{1}2[\omega,\omega]=0.
\]
By Theorem \ref{ref-5.1-85}, $ii)$ we have $\mathcal
U(\omega)=\omega$. Furthermore one checks that by degree reasons
$\mathcal U_\omega$ and $\mathcal S_\omega$ yield finite sums when
evaluated on specific elements. The same goes for the associated
homotopies. So the results of \S\ref{ref-7.1-110} apply.

Combining the arguments of the proof of Theorem~\ref{ref-5.3-88} with
Theorem~\ref{ref-7.1-111} we get the following commutative diagram of
precalculus structures up to homotopy
\begin{equation}\label{ref-7.4-117}
\xymatrix{\left(T_{\mathrm{poly},C^{\mathrm{aff},L}}(C^{\mathrm{aff},L}\widehat{\otimes}_{R_1}J\!L),{}^1\nabla^\mathrm{aff},[\ ,\ ],\cup\right)\ar@{~>}[d]\ar[r]^-{\mathfrak U_{L,1}} & \left(D_{\mathrm{poly},C^{\mathrm{aff},L}}(C^{\mathrm{aff},L}\widehat{\otimes}_{R_1}J\!L),{}^1\nabla^\mathrm{aff}+\mathrm d_\mathrm H,[\ ,\ ],\cup\right)\ar@{~>}[d]\\
\left(\Omega_{C^{\mathrm{aff},L}\widehat\otimes_{R_1}J\!L/C^{\mathrm{aff},L}},{}^1\nabla^{\mathrm{aff}},\mathrm L,\cap\right)  & \ar[l]_-{\mathfrak S_{L,0}}\left(\widehat{C}_{C^{\mathrm{aff},L},\bullet}(C^{\mathrm{aff},L}\widehat\otimes_{R_1}J\!L),{}^1\nabla^\mathrm{aff}+\mathrm b_\mathrm H,\mathrm L,\cap\right)}
\end{equation}
The fact that $\mathfrak U_{L,1}$ preserves the respective Lie
brackets up to homotopy is a consequence of the fact that $\mathfrak
U_L$ is an $L_\infty$-morphism; similarly, the fact that $\mathfrak
S_{L,0}$ preserves the Lie module structure up to homotopy is a
consequence of the fact that $\mathfrak S_L$ is an $L_\infty$-morphism
of $L_\infty$-modules.

On the other hand, $\mathfrak U_{L,1}$ is compatible with respect to
the products labelled by $\cup$ up to homotopy by the results
of~\cite[\S 10.1]{CVdB}.

As for the compatibility between the actions labelled by $\cap$ up to
homotopy, we first observe that the homotopy
formula~\eqref{ref-7.1-113} is well-defined in the case $(\mathfrak
m,\mathrm d_\mathfrak m)=(C^{\mathrm{coord},L},\mathrm d)$ and
$\gamma=\omega$, with the same notation as above: by the same
arguments as in the proof of Theorem~\ref{ref-5.3-88} it remains to
prove that the homotopy operator~\eqref{ref-7.2-114} descends to a
homotopy operator
\[
\mathfrak H_L^\mathcal S:T_{\mathrm{poly},C^{\mathrm{aff},L}}(C^{\mathrm{aff},L}\widehat{\otimes}_{R_1}J\!L)\otimes \widehat{C}_{C^{\mathrm{aff},L},\bullet}(C^{\mathrm{aff},L}\widehat\otimes_{R_1}J\!L)\to \Omega_{C^{\mathrm{aff},L}\widehat\otimes_{R_1}J\!L/C^{\mathrm{aff},L}}. 
\]
This holds true as a
consequence of Lemma~\ref{ref-7.2-115}
together with the verticality property of the Maurer--Cartan form
$\omega$, see~\S\ref{ref-4.1-54}.

If we now couple the commutative diagram~\eqref{ref-7.4-117} with the
results of \S\ref{ref-4.2-61},\ \S\ref{ref-4.3-65} and~\S\ref{ref-4.4-67}, and using
the same notation introduced at the end of the proof of
Theorem~\ref{ref-5.3-88} we get the
following commutative diagram of precalculi up to homotopy
\begin{equation}\label{ref-7.5-118}
\xymatrix{T_\mathrm{poly}^L(R) \ar@{^{(}->}[r]\ar[d] & \mathfrak g_1^L \ar[r]^{\mathfrak U_{L,1}}\ar[d]& \mathfrak g_2^L\ar[d] & D_\mathrm{poly}^L(R)\ar@{_{(}->}[l]\ar[d]\\
\Omega^L(R)\ar@{^{(}->}[r] & \mathfrak m_1^L & \mathfrak m_2^L\ar[l]_{\mathfrak S_{L,0}} & C_\mathrm{poly}^L(R)\ar@{_{(}->}[l]}.
\end{equation}
The quasi-isomorphisms $\mathfrak U_{L,1}$ and $\mathfrak S_{L,0}$ are
obtained from $\mathcal U_{L,\omega,1}$ and $\mathcal S_{L,\omega,0}$
respectively by means of the descent procedure: since $\omega$ is an
$\mathfrak m$-valued vector field in $T_\mathrm{poly}^\mathfrak
m(F)=\mathfrak g_1^L$, for $\mathfrak m=C^{\mathrm{coord},L}$, we can
use the results of~\cite[\S10.1]{CVdB}, and~\cite[\S6]{CR2}, to
evaluate explicitly $\mathcal U_{L,\omega,1}$ and $\mathcal
S_{L,\omega,0}$, namely
\begin{equation}\label{ref-7.6-119}
\begin{aligned}
\mathcal U_{L,\omega,1}=\mathrm{HKR}\circ\iota_{j(\omega)},\ \mathcal S_{L,\omega,0}=j(\omega)\wedge\mathrm{HKR},
\end{aligned}
\end{equation}
where
\begin{equation}
\label{ref-7.7-120}
j(\omega)=\det\sqrt{\frac{\Xi}{\exp\!\left(\frac{\Xi}2\right)-\exp\!\left(-\frac{\Xi}2\right)}},
\end{equation}
with $\Xi$ as defined in \eqref{ref-6.3-107}. To interpret \eqref{ref-7.7-120}
one should expand the right-hand side formally in terms of
$\operatorname{Tr}(\Xi^n)$ and then substitute the expression for $\Xi$
given in \eqref{ref-6.3-107}. This
yields an element of
$C^{\mathrm{coord},L}\widehat\otimes \Omega_F$ of degree $2n$. Thus
$j(\omega)$ is a sum of elements in
$C^{\mathrm{coord},L}\widehat\otimes \Omega_F$ of even total degree.

By the discussion in \S\ref{ref-6-104} the element
$\operatorname{Tr}(\Xi^n)\in C^{\mathrm{coord},L}\widehat\otimes
\Omega_F$ may be interpreted as an element in
$C^{\mathrm{aff},\mathcal L}\widehat\otimes_{\mathcal
  O_{1}}\Omega_{J\mathcal L/\mathcal O_{1}}$ {\em via} the inclusions
\eqref{ref-6.2-106}. Hence the same holds for $j(\omega)$. We keep
the same notation for this reinterpreted version of $j(\omega)$.

We thus get the following formul\ae:
\begin{equation}\label{ref-7.8-121}
  \mathfrak U_{L,1}=\mathrm{HKR}\circ \iota_{j(\omega)},\ \mathfrak S_{L,0}=j(\omega)\wedge \mathrm{HKR}.
\end{equation}

\subsection{Functoriality properties of the commutative diagram~\eqref{ref-7.5-118}}\label{ref-7.3-122}

The computations in the proof of Proposition~\ref{ref-5.5-93} imply
the following theorem, expressing the functoriality properties of the
commutative diagram~\eqref{ref-7.5-118}.
\begin{Thm}\label{ref-7.3-123}
  For a general algebraic morphism $(\ell,\lambda)$ from $(L,R)$ to
  $(M,S)$ as in Definition~\ref{ref-5.4-92}, which induces an isomorphism
  $S\otimes_R L\cong M$ of $S$-modules, and such that $L$ is free of
  rank $d$ over $R$ there exist quasi-isomorphisms $\mathfrak
  U_{L,1}$, $\mathfrak U_{M,1}$, $\mathfrak S_{L,0}$ and $\mathfrak
  S_{M,0}$, fitting into the commutative diagram of precalculi up to
  homotopy
\begin{equation}\label{ref-7.9-124}
  \xymatrix{T_\mathrm{poly}^L(R) \ar@{^{(}->}[r]\ar[d]_{\varphi_T} & \mathfrak g_1^L \ar[rrrr]^-{\mathfrak U_{L,1}=\mathrm{HKR}\circ 
\iota_{j(\omega_L)}}
\ar[d]_{\varphi_T}&&& &\mathfrak g_2^L\ar[d]_{\varphi_D} & D_\mathrm{poly}^L(R)\ar@{_{(}->}[l]\ar[d]_{\varphi_D}\\
      T_\mathrm{poly}^M(S) \ar@{^{(}->}[r]\ar@{~>}[d]_{\mathrm L} & \mathfrak g_1^M \ar[rrrr]^-{\mathfrak U_{M,1}=\mathrm{HKR}\circ 
\iota_{j(\omega_M)}}
\ar@{~>}[d]_{\mathrm L_1}&&& &\mathfrak g_2^M\ar@{~>}[d]_{\mathrm L_2} & D_\mathrm{poly}^M(S)\ar@{_{(}->}[l]\ar@{~>}[d]_{\mathrm L}\\
      \Omega^M(S)\ar@{^{(}->}[r] & \mathfrak m_1^M &&& &\mathfrak m_2^M\ar[llll]_-{\mathfrak S_{M,0}=
j(\omega_M)
\wedge \mathrm{HKR}} & C_\mathrm{poly}^M(S)\ar@{_{(}->}[l]\\
      \Omega^L(R)\ar@{^{(}->}[r]\ar[u]^{\varphi_\Omega}& \mathfrak m_1^L \ar[u]^{\varphi_\Omega}&&& & \mathfrak m_2^L\ar[llll]_-{\mathfrak S_{L,0}=
j(\omega_L)
\wedge \mathrm{HKR}} \ar[u]^{\varphi_J}& C_\mathrm{poly}^L(R)\ar@{_{(}->}[l]\ar[u]^{\varphi_J}},
\end{equation}
where we borrow notation from Proposition~\ref{ref-5.5-93}, and where $\omega_L$ and $\omega_M$, denote the Maurer--Cartan form on $C^{\mathrm{coord},L}$ and $C^{\mathrm{coord},M}$ respectively.  
The precalculus structures up to homotopy on $(\mathfrak{g}^\ast_i,\mathfrak{m}^\ast_i)$, $\ast=L,M$, $i=1,2$, are defined as in \S\ref{ref-5.5-87}.
{}{Moreover the implied homotopies are in a similar way functorial
for algebraic morphisms $(\ell,\lambda)$ from $(L,R)$ to
  $(M,S)$ satisfying $S\otimes_R L\cong M$.}
\end{Thm}
Almost all important objects appearing in Theorem~\ref{ref-7.3-123}
have already appeared in Theorem~\ref{ref-5.6-101}, hence the
functoriality properties extend  to the present situation. The commutativity
of the upper and lower squares involving $j(\omega)$ follows from the
compatibility of the inclusions \eqref{ref-6.2-106} with the base extension
$S/R$. 
{}{The functoriality properties of the implied homotopies are verified
in the same way.
See~\cite[Lemma 10.1.1]{CVdB} for results on $\mathfrak U_{\ast,1}$ and related homotopies; in virtue of Lemma~\ref{ref-7.2-115}, the hotomopy expressing the compatibility of $\mathfrak S_{\ast,0}$ with cap products descends correctly on $C^{\mathrm{aff},\ast}$, and the functoriality properties of such a homotopy follow along the same lines of the functoriality properties in Theorem~\ref{ref-5.6-101}, as the homotopy under consideration is expressed in terms of scalar combinations of poly-differential operators associated to certain graphs, as the $L_\infty$-quasi-isomorphisms of Kontsevich and Shoikhet.} 
 
\subsubsection{Arbitrary base fields}
\label{sec-arbitrary}
We now briefly indicate how we may replace $k$ by a general field of
characteristic zero. Our arguments depend on the existence of a number
of explicit homotopies. 
{}{These homotopies are constructed as scalar
linear combinations of poly-differential operators indexed by certain graphs, where the scalars depend only on the corresponding graphs.}
  For the arguments to
work the coefficients need to satisfy certain linear equations. These
equations have a solution over $\mathbb{R}$ (given that over this
field we have homotopies that work). Thus they have a solution over
any field of characteristic zero.

\medskip

{}{We will now be more specific.} We refer to \cite[\S10.4]{CVdB} for what concerns Lie brackets and cup
products; here we concentrate on the compatibility between cap products.  
{}{We embed $k$ in a field $K$ containing
$\mathbb{R}$.}
By virtue of \cite[\S
6]{CR2},  {}{$\mathcal{U}_{\gamma,1}$} and $\mathcal{S}_{\gamma,0}$ are
defined over $\mathbb{Q}$ and thus $k$ (while they are {\em a priori}
defined over  {}{$\mathbb{R}\subset K$}). Then observe that
equation \eqref{ref-7.1-113} is linear in the coefficients
$\overset{\circ}{W}_{D,\Gamma}$ of $\mathcal H_\gamma^{\mathcal{S}}$. Since we already
have a solution of these equations in 
{}{$\mathbb{R}\subset K$},
we get one in $k$ by applying any projection  {}{$K\to k$}.

\subsection{Proof of Theorem~\ref{ref-1.1-9} in the global case
}\label{ref-7.4-125}
Let $(X,\mathcal O)$ be a ringed site and $\mathcal L$ be a locally free sheaf of Lie algebroids over $\mathcal O$ of rank $d$. We denote by $D(X)$ the derived category of sheaves of $k$-vector spaces over $X$. According to the results of Section~\ref{ref-3-16}, transported to the
framework of sheaves of $k$-vector spaces, $(T_\mathrm{poly}^\mathcal
L(X),\Omega^\mathcal L(X))$ and $(D_\mathrm{poly}^\mathcal
L(X),C_\mathrm{poly}^\mathcal L(X))$ are precalculi up to
homotopy. Therefore, viewed as objects of $D(X)$ they are genuine
precalculi.

Additionally, the sheafification procedure can be applied to the
commutative diagram~\eqref{ref-7.9-124}, in virtue of the results
of \S\ref{ref-7.3-122} {}{(using the fact that the homotopies are functorial
as well)} : if we further consider the resulting
commutative diagram of sheaves of $k$-vector spaces in the derived
category $D(X)$, then using \eqref{ref-6.4-108}
we get the commutative diagram of precalculi
\[
\xymatrix{T_\mathrm{poly}^\mathcal L(X)\ar[dd]|{\mathrm{HKR}\circ
    \iota_{\widetilde{\mathrm{td}} (\mathcal L)^{1/2}}
  }\ar[rr]\ar@{~>}[dr] & & \mathfrak g_1^\mathcal L\ar[dd]|(.25){\mathrm{HKR}\circ\iota_{j(\omega)}}\ar@{~>}[drr] & & \\
  & \Omega^\mathcal L(X)\ar'[r][rrr] & & &\mathfrak m_1^\mathcal L\\
  D_\mathrm{poly}^\mathcal L(X)\ar@{~>}[dr]\ar[rr] & & \mathfrak g_2^\mathcal L\ar@{~>}[drr] & &\\
  & C_\mathrm{poly}^\mathcal
  L(X)\ar'[u][uu]|(0.25){\widetilde{\mathrm{td}}(\mathcal L)^{1/2}\wedge
    \mathrm{HKR}}\ar[rrr] & & &\mathfrak m_2^\mathcal
  L\ar[uu]|{j(\omega)\wedge
    \mathrm{HKR}}},
\]
where all horizontal and vertical arrows represent isomorphisms in the
derived category $D(X)$. Here $\widetilde{\mathrm{td}}(\mathcal{L})$
is the modified Todd class of $\mathcal L$ which is obtained
by replacing the function $q(x)$ in the definition of the Todd class
(see \eqref{ref-1.3-4}) by
\[
\widetilde{q}(x)=\frac{x}{e^{x/2}-e^{-x/2}}.
\]
Hence at this point we have proved Theorem~\ref{ref-1.1-9} provided
that we replace the Todd class by the modified one. To obtain the 
result for the ordinary Todd class we follow
the method of \cite[\S10.3]{CVdB}.
We have
\def\Td{\operatorname{td}}\def\Lscr{\mathcal{L}}\def\Tr{\operatorname{Tr}}
\begin{align*}
\widetilde{\Td}(\Lscr)&= \Td(\Lscr)\det(e^{-A(\Lscr)/2})\\
&=\Td(\Lscr)e^{-\Tr(A(\Lscr))/2}\\
&=\Td(\Lscr)e^{-a_1(\Lscr)/2}
\end{align*}
In other words it is sufficient to prove that $(\iota_{\textstyle e^{-a_1(\Lscr)/4}},
e^{a_1(\Lscr)/4}\wedge-)$
defines an automorphism of the precalculus $(T^{\mathcal L}_{\poly}(X),
\Omega^{\mathcal L}(X))$. 

Via the inclusions \eqref{ref-6.2-106} together with
\eqref{ref-6.4-108} we may as well prove that
$(\iota_{e^{-\Tr(\Xi)/4}}, e^{\Tr(\Xi)/4}\wedge-)$ defines an
automorphism of the precalculus $(C^{\coord}\widehat{\otimes}
T_{\poly}(F),C^{\coord}\widehat{\otimes} \Omega_F)$ or equivalently
that $(\iota_{\Tr(\Xi)},-\Tr(\Xi)\wedge-)$ act as derivations.  The
fact that $\iota_{\Tr(\Xi)}$ is a derivation with respect to the cup product
and Lie bracket has been checked in \cite[\S10.3]{CVdB}.  So it
remains to show compatibility with the cap product and Lie derivative.

As
$\Tr(\Xi)=\sum_{i,\alpha}\eta_\alpha \mathrm d_F(\partial_i\omega^i_\alpha)$
we first derive some identities for $\iota_{\mathrm d_F b}$ and $\mathrm
d_F b\wedge-$ with $b$ in $F$.

First we claim
\begin{equation}
\label{ref-7.10-126}
\mathrm d_F b\wedge (D\cap \sigma)=-\iota_{\mathrm d_F b}(D)\cap \sigma+(-1)^{|D|+1}D\cap (\mathrm d_F b\wedge \sigma)
\end{equation}
for $b\in F$, $D\in T_{\poly}(F)$, $\sigma\in \Omega_F$.  If $D=D_1\cup D_2$
and \eqref{ref-7.10-126} holds for $D_1$, $D_2$ then it holds for $D$ as well.  To
see this note
\begin{align*}
\mathrm d_F b\wedge ((D_1\cup D_2)\cap \sigma)&=\mathrm d_F b\wedge (D_1\cap (D_2\cap \sigma))\\
&=-\iota_{\mathrm d_F b}(D_1) \cap (D_2\cap \sigma)+(-1)^{|D_1|+1} D_1 \cap(\mathrm d_F b\wedge (D_2\cap \sigma))\\
&=-\iota_{\mathrm d_F b}(D_1) \cap (D_2\cap \sigma)-(-1)^{|D_1|+1}D_1\cap \iota_{\mathrm d_F b}(D_2)\cap
\sigma\\
&\qquad+(-1)^{|D_1|+|D_2|} D_1\cap D_2\cap (\mathrm d_F b\wedge \sigma)\\
&=-\iota_{\mathrm d_F b}(D_1\cup D_2)\cap \sigma+(-1)^{|D_1\cup D_2|+1}(D_1\cup D_2)\cap (\mathrm d_F b\wedge \sigma)
\end{align*}
So we only have to consider the case where $D$ is a function or a
vector field.  The case that $D$ is a function is trivial so assume
that $D$ is a vector field. In that case we find for the right-hand
side of \eqref{ref-7.10-126}
\begin{align*}
-\iota_{\mathrm d_F b}(D)\cap \sigma+(-1)^{|D|+1}D\cap (\mathrm d_F b\wedge \sigma)&=
Db\cap \sigma-Db\wedge \sigma+\mathrm d_F b\wedge(D\cap\sigma)\\
&=\mathrm d_F b\wedge (D\cap\sigma)
\end{align*}
which is equal to the left-hand side of \eqref{ref-7.10-126}.

For the Lie derivative we use
$
\mathrm L_D=[\mathrm d_F, D\cap-]
$.
It is clear that $\mathrm d_F$ and $\mathrm d_F b\wedge-$ commute. We
then compute using \eqref{ref-7.10-126}
\begin{align*}
  \mathrm d_F b\wedge \mathrm L_D\sigma&=\mathrm d_F b\wedge (\mathrm
  d_F(D\cap\sigma)-(-1)^{|D|+1}
  D\cap \mathrm d_F \sigma)\\
  &=-\mathrm d_F (\mathrm d_F b\wedge (D\cap \sigma))+(-1)^{|D|}\mathrm d_F b\wedge (D\cap \mathrm d_F \sigma)\\
  &=\mathrm d_F(\iota_{\mathrm d_F b}(D)\cap \sigma)+(-1)^{|D|}\mathrm
  d_F(D\cap
  (\mathrm d_F b\wedge \sigma))\\
  &\qquad +(-1)^{|D|+1} \iota_{\mathrm d_F b}(D)\cap d\sigma-D\cap (\mathrm d_F b\wedge\mathrm d_F \sigma)\\
  &=\mathrm L_{\iota_{\mathrm d_F b} D}(\sigma)+(-1)^{|D|}\mathrm L_D(\mathrm
  d_F b\wedge \sigma)
\end{align*}
If $\eta$ is an odd element in $C^{\coord}$ then $\iota_{\eta \mathrm
  d_F b} D=\eta \iota_{\mathrm d_F b}D$ and $\mathrm L_{\iota_{\eta \mathrm
    d_F b}D}(\sigma)= \mathrm L_{\eta \iota_{\mathrm d_F b}D}(\sigma)=-\eta
\mathrm L_{\iota_{\mathrm d_F b}}(\sigma)$. Using this we find
\[
\mathrm \Tr(\Xi)\wedge (D\cap \sigma)=-\iota_{\Tr(\Xi)}(D)\cap \sigma+D\cap (\Tr(\Xi)\wedge \sigma)
\]
and
\[
\Tr(\Xi)\wedge \mathrm L_D\sigma=-\mathrm L_{\iota_{\Tr(\Xi)}D}(\sigma)+\mathrm L_D(\Tr(\Xi)\wedge \sigma)
\]
We conclude that $(\imath_{\Tr(\Xi)},-\Tr(\Xi)\wedge-)$ does indeed define a
derivation of precalculi. 

\appendix
\section{Explicit formul\ae\ for the $B_\infty$-structure on
poly-differential operators}
\label{ref-A-127}
In this appendix and the next one we develop the precalculus structure on
$L$-chains over $L$-cochains up to homotopy.  The results in these
appendices are provided for background and are not essential for the
results in the body of the paper.

The graded vector space $V=D_\mathrm{poly}^L(R)$ is naturally a
$B_\infty$-algebra.  This means that the cofree
coassociative coalgebra (with counit) $\mathrm T(V)$ is canonically equipped
with the structure of a DG bialgebra.
The notion of $B_\infty$-algebra has been introduced in~\cite{Bau}; though, we make use here mainly of the $B_\infty$-algebra structure given by braces~\cite{GV,GJ}, to which we refer for more details, see also~\cite[Sections 1,2]{CR2}.

The corresponding associative product $\mathrm m$ on $\mathrm T(V)$ is uniquely
determined by its Taylor components $m_{p,q}:\mathrm T^p(V)\otimes \mathrm T^q(V)\rightarrow V$.
We have $m_{p,q}=0$ if $p\neq 1$ and
{\tiny
\begin{equation}\label{ref-A.1-128}
\begin{aligned}
&\mathrm m_{1,q}(D\otimes (D_1\otimes\cdots\otimes D_q))=D\{D_1,\dots,D_q\}=\\
&=\sum_{1\leq i_1\leq\cdots\leq i_q\leq |D|+\sum_{b=1}^{q-1}|D_b|+1}(-1)^{\sum_{k=1}^q |D_k|\left(i_k-1\right)}\\
&\left(1^{\otimes(i_1-1)}\otimes \Delta^{|D_1|}\otimes 1^{\otimes(i_2-i_1-|D_1|-1)}\otimes \Delta^{|D_2|}\otimes \cdots\otimes 1^{\otimes(i_q-i_{q-1}-|D_{q-1}|-1)}\otimes \Delta^{|D_q|}\otimes 1^{\otimes(|D|+\sum_{b=1}^{q-1}|D_b|-i_q)}\right)(D)\\
&\left(1^{\otimes(i_1-1)}\otimes D_1\otimes 1^{\otimes(i_2-i_1-|D_1|-1)}\otimes D_2\otimes \cdots\otimes 1^{\otimes(i_q-i_{q-1}-|D_{q-1}|-1)}\otimes D_q\otimes 1^{\otimes(|D|+\sum_{b=1}^{q-1}|D_b|-i_q)}\right),
\end{aligned}
\end{equation}
} 
 for
elements $D$, $D_i$, $i=1,\dots,q$, of $D_\mathrm{poly}^L(R)$, where
$|-|$ denotes the (shifted) degree of elements of
$D_\mathrm{poly}^L(R)$: accordingly, we have
$|D\{D_1,\dots,D_q\}|=|D|+\sum_{a=1}^q|D_a|$ and thus all brace operations
are of degree zero. In the sum
(\ref{ref-A.1-128}), we have $1\leq i_1$, $i_k+|D_k|+1\leq i_{k+1}$,
$k=1,\dots,q-1$, $i_q+|D_q|\leq |D|+\sum_{a=1}^q|D_a|+1$. 
 The sign
conventions are taken from \cite{CR2}. 
The brace operations~\eqref{ref-A.1-128}
satisfy an infinite family of quadratic identities (see
e.g.~\cite{CR2}), which are equivalent to the associativity of the
product $\mathrm m$.  

\medskip

We define the cup product by means of the brace operations, see also~\cite{GV, CR2}, {\em via} the assignment 
\begin{equation}\label{ref-A.2-129}
D_1\cup D_2=(-1)^{|D_1|+1}\mu\{D_1,D_2\},\ D_i\in D_\mathrm{poly}^L(R),\ i=1,2.
\end{equation}
It is obvious that the cup product has (shifted) degree $1$. 
An easy verification using Formula~\eqref{ref-A.1-128} shows that the previous definition of cup product coincides with the one given in Formula~\eqref{ref-3.11-36}.

We now have the following compatibilities
\begin{Lem}\label{ref-A.1-130}
  The degree $0$ operation~\eqref{ref-3.10-35} and the degree $1$
  operation~\eqref{ref-A.2-129} satisfy the following properties:
\begin{align}
\label{ref-A.3-131}&[D_1,D_2]=-(-1)^{|D_1||D_2|}[D_2,D_1],\\
\label{ref-A.4-132}&[D_1,[D_2,D_3]]=[[D_1,D_2],D_3]+(-1)^{|D_1||D_2|}[D_2,[D_1,D_3]],\\
\label{ref-A.5-133} &D_1\cup D_2=(-1)^{(|D_1|-1)(|D_2|-1)}D_2\cup D_1\pm \left(\mathrm d_\mathrm H(D_1\{D_2\})-(\mathrm d_\mathrm H D_1)\{D_2\}-(-1)^{|D_1|}D_1\{\mathrm d_\mathrm H D_2\}\right),\\
\label{ref-A.6-134} &D_1\cup (D_2\cup D_3)=(D_1\cup D_2)\cup D_3,
\end{align}
and 
\begin{equation}\label{ref-A.7-135}
\begin{aligned}
&[D_1,D_2\cup D_3]&=[D_1,D_2]\cup D_3+(-1)^{|D_1|(|D_2|-1)} D_2\cup [D_1,D_3]
+(-1)^{|D_1|}\biggl(\mathrm d_\mathrm H (D_1\{D_2,D_3\})-\\
& &\phantom{=}(\mathrm d_\mathrm H D_1)\{D_2,D_3\}-(-1)^{|D_1|}D_1\{\mathrm d_\mathrm H D_2,D_3\}-(-1)^{|D_1|+|D_2|}D_1\{D_2,\mathrm d_\mathrm H D_3\}\biggr),
\end{aligned}
\end{equation}
for general elements $D_i$ of $D_\mathrm{poly}^L(R)$, $i=1,2,3$, and where $\mathrm d_\mathrm H=[\mu,\bullet]$, $\mu=1\otimes_R 1$. 
\end{Lem}
\section{The precalculus structure on $L$-chains}\label{ref-B-136}
We need results from~\cite{CR2,TT} about algebraic structures on Hochschild (co)chains, which have to be adapted to the Lie algebroid framework.

According to~\cite{TT,CR2}, there are two distinct, non-compatible, left $B_\infty$-module structures
on the Hochschild chain complex of $A$, viewed as a $B_\infty$-algebra
with respect to the brace operations~\eqref{ref-A.1-128}.  Equivalently,
we view the two left $B_\infty$-module structures on the Hochschild chain
complex as the data of two left actions $\mathrm m_{L,i}$, $i=1,2$, on the left comodule cofreely cogenerated by the
Hochschild chain complex of $A$ over the coalgebra cofreely
cogenerated by the Hochschild cochain complex of $A$.

These results can be applied to the present situation with due
changes: $\widehat{C}_{R,\bullet}(J\!L)$ has two left $B_\infty$-module structures over the $B_\infty$-algebra
$D_\mathrm{poly}^L(R)$.

We borrow the main notation and sign conventions from~\cite{CR2}.  We
denote by $\mathrm m_{L,i}$, $i=1,2$ the two left $B_\infty$-module structures on $\widehat{C}_{R,\bullet}(J\!L)$: they are
uniquely determined by their Taylor components
\begin{equation}\label{ref-B.1-137}
{\tiny
\begin{aligned}
&\left(\mathrm m_{L,1}^{1,q,r}(P\otimes (Q_1\otimes\cdots\otimes Q_q)\otimes a\otimes (R_1\otimes\cdots\otimes R_r))\right)(D)=\\
&=\sum_{l=-|a|-\sum_{b=1}^q|Q_b|-q+1\ \mathrm{mod}\ (-|a|+1)}^{-|a|-|P|-\sum_{b=1}^q|Q_b|+r+1\ \mathrm{mod}\ (-|a|+1)}\sum_{l\leq j_1\leq\cdots\leq j_q\leq -|a|\atop 1\leq k_1\leq\cdots\leq k_r\leq |a|+|P|+l}(-1)^{l(-|a|-l+1)+\sum_{b=1}^q |Q_b|(j_b-l)+\sum_{c=1}^r|R_c|(k_c-l-1)}\\
&\sigma^{(-|a|-l+1)}(a)\left((\Delta^{(|P|+\sum_{b=1}^q|Q_b|+\sum_{c=1}^r|R_c|)}\otimes 1^{\otimes (-|a|-|P|-\sum_{b=1}^q|Q_b|-\sum_{c=1}^r|R_c|)})(D)\right.\\
&\left.\left(1^{(j_1-l)}\otimes \Delta^{|Q_1|}\otimes\cdots\otimes 1^{\otimes(j_q-j_{q-1}-|Q_{q-1}|-1)}\otimes \Delta^{|Q_q|}\otimes 1^{\otimes (-|a|-j_q-|Q_q|+k_1)} \otimes \Delta^{|R_1|}\otimes \cdots \otimes\right.\right.\\
&\left.\left. \otimes 1^{\otimes (k_r-k_{r-1}-|R_{r-1}|-1)}\otimes \Delta^{|R_r|}\otimes 1^{\otimes (|a|+|P|+\sum_{b=1}^q|Q_b|+\sum_{c=1}^{r-1}|R_c|+l-k_r-1)}\otimes 1^{\otimes (-|a|-|P|-\sum_{b=1}^q|Q_b|-\sum_{c=1}^r|R_c|)}\right)(P)\right.\\
&\left.\left(1^{(j_1-l)}\otimes Q_1\otimes\cdots\otimes 1^{\otimes(j_q-j_{q-1}-|Q_{q-1}|-1)}\otimes Q_q\otimes 1^{\otimes (-|a|-j_q-|Q_q|+k_1)} \otimes R_1\otimes \cdots \otimes\right.\right.\\
&\left.\left.\otimes 1^{\otimes (k_r-k_{r-1}-|R_{r-1}|-1)}\otimes R_r\otimes 1^{\otimes (|a|+|P|+\sum_{b=1}^q|Q_b|+\sum_{c=1}^{r-1}|R_c|+l-k_r-1)}\otimes 1^{\otimes (-|a|-|P|-\sum_{b=1}^q|Q_b|-\sum_{c=1}^r|R_c|)}\right)\right),
\end{aligned}
}
\end{equation}
where $\sigma$ is the operator on $\widehat{C}_{R,\bullet}(J\!L)$ defined via
\[
\sigma(a)(D_0\otimes\cdots D_{-|a|})=a(D_1\otimes\cdots D_{-|a|}\otimes D_0),\ D_i\in \mathrm U_R(L),\ i=0,\dots,-|a|,
\]
which obviously satisfies $\sigma^{(-|a|+1)}=\mathrm{id}$, and the indices in the summation satisfy $l\leq j_1$, $j_i+|Q_i|+1\leq j_{i+1}$, $i=1,\dots,q-1$, $j_q+|Q_q|\leq -|a|$, $k_i+|R_i|+1\leq k_{i+1}$, $i=1,\dots,r-1$, $k_r+|R_r|\leq |a|+|P|+\sum_{b=1}^q|Q_b|+\sum_{j=1}^c|R_c|+l-1$, and
\begin{equation}\label{ref-B.2-138}
{\tiny
\begin{aligned}
& \left(\mathrm m_{L,2}^{0,0,r}(a\otimes (R_1\otimes\cdots\otimes R_r))\right)(D)=\\
&=\sum_{1\leq i_1\leq \cdots\leq i_p\leq -|a|}(-1)^{\sum_{c=1}^r |R_c|\left(i_c-1\right)}\\
&\phantom{=}a\left(\left(1^{\otimes i_1}\otimes \Delta^{|R_1|}\otimes 1^{\otimes(i_2-i_1-|R_1|-1)}\otimes \Delta^{|R_2|}\otimes \cdots\otimes 1^{\otimes(i_r-i_{r-1}-|R_{r-1}|-1)}\otimes \Delta^{|R_r|}\otimes 1^{\otimes(|D|+\sum_{c=1}^{r-1}|R_c|-i_r)}\right)(D)\right.\\
&\phantom{=}\left.\left(1^{\otimes i_1}\otimes R_1\otimes 1^{\otimes(i_2-i_1-|R_1|-1)}\otimes R_2\otimes \cdots\otimes 1^{\otimes(i_r-i_{r-1}-|R_{r-1}|-1)}\otimes R_r\otimes 1^{\otimes(|D|+\sum_{c=1}^{r-1}|R_c|-i_r)}\right)\right),
\end{aligned}
}
\end{equation}
where the summation is over indices $i_1,\dots,i_r$, such that $1\leq i_1$, $i_k+|D_k|+1\leq i_{k+1}$, $k=1,\dots,p-1$, $i_p+|D_p|\leq -|a|$.
We observe that the components of $\mathrm m_{L}$, resp.\ $\mathrm m_{L}$, are non-trivial only if $p\leq 1$, with no restrictions on $q,r$, resp.\ only if $q=r=0$, with no restrictions on $p$.

It is not difficult but quite tedious to verify that both~\eqref{ref-B.1-137} and~\eqref{ref-B.2-138} have degree $0$ and satisfy an infinite family of quadratic relations involving braces.

The Taylor components of $\mathrm m_{L,i}$, $i=1,2$, permit to define a pairing of degree $0$ between $D_\mathrm{poly}^L(R)$ and $\widehat{C}_{R,\bullet}(J\!L)$ via
\begin{equation}\label{ref-B.3-139}
\mathrm L_Da=\mathrm m_{L,1}^{1,0,0}(D\otimes a)+(-1)^{|D|}\mathrm m_{L,2}^{0,0,1}(a\otimes D),\ D\in D_\mathrm{poly}^L(R),\ a\in \widehat{C}_{R,\bullet}(J\!L).
\end{equation}
Similarly, we may consider two distinct pairings between $D_\mathrm{poly}^L(R)$ and $\widehat{C}_{R,\bullet}(J\!L)$: for $\mu$ as above,
\begin{align}
\label{ref-B.4-140}D\cap a&=(-1)^{|D|}\mathrm m_{L,1}^{1,1,0}(\mu\otimes D\otimes a),\\
\label{ref-B.5-141}a\cap D&=(-1)^{|a|}\mathrm m_{L,1}^{1,0,1}(\mu\otimes a\otimes D),\ D\in\mathrm D_\mathrm{poly}^L(R),\ a\in \widehat{C}_{R,\bullet}(J\!L).
\end{align}
It follows from their very definition that both~\eqref{ref-B.4-140} and~\eqref{ref-B.5-141} have degree $1$.
\begin{Lem}\label{ref-B.1-142}
The pairing~\eqref{ref-B.3-139} of degree $0$ and the pairings~\eqref{ref-B.4-140} and~\eqref{ref-B.5-141} of degree $1$ satisfy the following properties:
\begin{align}\
\label{ref-B.6-143}&\mathrm L_{[D_1,D_2]}a=\mathrm L_{D_1}(\mathrm L_{D_2}a)-(-1)^{|D_1||D_2|}\mathrm L_{D_2}(\mathrm L_{D_1}a),\\ 
\label{ref-B.7-144}&D\cap a=(-1)^{(|D|-1)(|a|-1)}a\cap D\pm \left(\mathrm b_\mathrm H(\mathrm m_{L,1}^{1,0,0}(D\otimes a)-\mathrm m_{L,1}^{1,0,0}(\mathrm d_\mathrm H D\otimes a)-(-1)^{|D|}\mathrm m_{L,1}^{1,0,0}(D\otimes \mathrm b_\mathrm H a)\right),\\
\label{ref-B.8-145}&D_1\cap (D_2\cap a)=(D_1\cup D_2)\cap a,\\
\label{ref-B.9-146}&(a\cap D_1)\cap D_2=a\cap (D_1\cup D_2), 
\end{align}
\begin{equation}\label{ref-B.10-147}
\begin{aligned}
  \mathrm L_{D_1}(D_2\cap a)&=[D_1,D_2]\cap a+(-1)^{|D_1|(|D_2|-1)}D_2\cap\mathrm L_{D_1}a+(-1)^{|D_1|}\left(\mathrm b_\mathrm H(\mathrm m_{L,1}^{1,1,0}(D_1\otimes D_2\otimes a))\right.\\
  &\qquad\left.-\mathrm m_{L,1}^{1,1,0}(\mathrm d_\mathrm H
    D_1\otimes D_2\otimes a)-(-1)^{|D_1|}\mathrm
    m_{L,1}^{1,1,0}(D_1\otimes \mathrm d_\mathrm HD_2\otimes
    a)\right.\\
&\qquad \left.-(-1)^{|D_1|+|D_2|}\mathrm m_{L,1}^{1,1,0}(D_1\otimes D_2\otimes
    \mathrm b_\mathrm Ha)\right),
\end{aligned}
\end{equation}
\begin{equation}\label{ref-B.11-148}
\begin{aligned}
  \mathrm L_{D_1}(a\cap D_2)&=\mathrm L_{D_1}a\cap D_2+(-1)^{|D_1|(|a|-1)}a\cap[D_1,D_2]+(-1)^{|D_1|}\left(\mathrm b_\mathrm H(\mathrm m_{L,1}^{1,0,1}(D_1\otimes a\otimes D_2))\right.\\
  &\qquad\left.-\mathrm m_{L,1}^{1,0,1}(\mathrm d_\mathrm HD_1\otimes
    a\otimes D_2)-(-1)^{|D_1|}\mathrm m_{L,1}^{1,0,1}(D_1\otimes \mathrm
    b_\mathrm H a\otimes D_2)\right.\\
&\qquad\left.-(-1)^{|D_1|+|a|}\mathrm
    m_{L,1}^{1,0,1}(D_1\otimes a\otimes \mathrm d_\mathrm HD_2)\right),
\end{aligned}
\end{equation}
and finally
\begin{equation}\label{ref-B.12-149}
\begin{aligned}
  &\mathrm L_{D_1\cup D_2}a+(-1)^{(|D_1|-1)(|D_2|-1)}\mathrm L_{D_2\cup D_1}a\\
  &=\left(D_1\cap \mathrm L_{D_2}a+(-1)^{(|D_1|-1)(|D_2|+|a|-1)}\mathrm L_{D_2}a\cap D_1\right)\\
  &\qquad+(-1)^{|a|(|D_2|-1)}\left(\mathrm L_{D_1}a\cap D_2+(-1)^{(|D_1|+|a|-1)(|D_2|-1)}D_2\cap \mathrm L_{D_1}a\right)\\
  &\qquad+(-1)^{(|D_2|-1)}\left([D_1,D_2]\cap a+(-1)^{(|a|-1)(|D_1|+|D_2|-1)}a\cap [D_1,D_2]\right)\\
  &\qquad+(-1)^{|D_2|}\mathrm b_\mathrm H(\mathrm m_{L,2}^{0,0,2}(a\otimes R_1\otimes R_2))-(-1)^{|D_1|}\mathrm m_{L,2}^{0,0,2}(\mathrm b_\mathrm H a\otimes D_1\otimes D_2)\\
  &\qquad+(-1)^{|D_2|}\mathrm m_{L,2}^{0,0,2}(a\otimes \mathrm d_\mathrm H D_1\otimes D_2)+(-1)^{|D_1|+|D_2|}\mathrm m_{L,2}^{0,0,2}(a \otimes D_1\otimes \mathrm d_\mathrm H D_2)\\
&\qquad+(-1)^{|D_1|}\mathrm b_\mathrm H(\mathrm m_{L,2}^{0,0,2}(a\otimes D_2\otimes D_1))-(-1)^{|D_2|}\mathrm m_{L,2}^{0,0,2}(\mathrm b_\mathrm H a\otimes D_2\otimes D_1)\\
  &\qquad+(-1)^{|D_1|}\mathrm m_{L,2}^{0,0,2}(a\otimes \mathrm d_\mathrm
  H D_2\otimes D_1)+(-1)^{|D_1|+|D_2|}\mathrm m_{L,2}^{0,0,2}(a\otimes
  D_2\otimes \mathrm d_\mathrm H D_1),
\end{aligned}
\end{equation}
for a general element $a$ of $\widehat{C}_{R,\bullet}(J\!L)$ and general elements $D$, $D_i$, $i=1,2$, of $D_\mathrm{poly}^L(R)$, and where $\mathrm b_\mathrm H=\mathrm L_\mu$, for $\mu$ as before.
\end{Lem}
As for Lemma~\ref{ref-A.1-130}, the proof essentially makes use of the
brace identities, of the fact that $\mathrm m_{L,i}$, $i=1,2$, is a left action with respect to the brace operations,
and of the fact that $\mathrm m_{L,1}$ and $\mathrm m_{L,2}$ satisfy a
weak compatibility, as explained in more details in~\cite{CR2}.

Both actions $\mathrm m_{L,1}$ and $\mathrm m_{L,2}$ are compatible with the Grothendieck connection, i.e.\ 
\[
\begin{aligned}
{}^1\nabla_l\!\left(\mathrm m_{L,1}^{1,q,r}(D\otimes Q_1\otimes\cdots\otimes a\otimes R_1\otimes\cdots)\right)&=\mathrm m_{L,1}^{1,q,r}(D\otimes Q_1\otimes\cdots\otimes {}^1\nabla_l a\otimes R_1\otimes\cdots),\ q,r\geq 0,\\
{}^1\nabla_l\!\left(\mathrm m_{L,2}^{0,0,r}(D_1\otimes \cdots\otimes a)\right)&=\mathrm m_{L,2}^{p,0,0}(D_1\otimes \cdots\otimes {}^1\nabla_l a),\ p\geq 0, 
\end{aligned}
\]
for $D$, $D_i$ ($i=1,\dots,p$), $Q_j$ ($j=1,\dots,q$), $R_k$ ($k=1,\dots,r$) elements of $D_\mathrm{poly}^L(R)$, and $a$ of $\widehat{C}_{R,\bullet}(J\!L)$.
Both identities follow from the fact that ${}^1\nabla$ commutes with the operator $\sigma$ and from the fact that $\mathrm U_R(L)$ is a Hopf algebroid, in particular, the comultiplication is an algebra morphism.

Then, in virtue of Lemma~\ref{ref-B.1-142}, the
pairings~\eqref{ref-B.3-139},~\eqref{ref-B.4-140} and~\eqref{ref-B.5-141} are
compatible with the Grothendieck connection, implying in particular
that the Hochschild differential is also compatible therewith.  By the
very same arguments,
Formul\ae~\eqref{ref-B.7-144}, \eqref{ref-B.8-145}, \eqref{ref-B.9-146}, \eqref{ref-B.10-147}, \eqref{ref-B.11-148}
and~\eqref{ref-B.12-149} are compatible with the Grothendieck
connection, whence $(\mathrm{Ker}({}^1\nabla)\cap
\widehat{C}_{R,\bullet}(J\!L),\mathrm b_\mathrm H,\mathrm L,\cap)$, where $\cap$
denotes here both~\eqref{ref-B.4-140} and~\eqref{ref-B.5-141}, inherits a
structure of precalculus up to homotopy over the Gerstenhaber
algebra $(D_\mathrm{poly}^L(R),\mathrm d_\mathrm H,[\ ,\ ],\cup)$ up
to homotopy.

For the sake of completeness, we write down explicit formul\ae\ for
the Hochschild differential $\mathrm b_\mathrm H$ on the complex of
Hochschild $L$-chains on $R$ and for the pairing~\eqref{ref-B.5-141}
between $D_\mathrm{poly}^L(R)$ and $C_\mathrm{poly}^L(R)$;
in~\cite{CRVdB2}, we will deduce the same formul\ae\ in the framework
of homological algebra and derived functors.  Explicitly,
\[
\begin{aligned}
\mathrm b_\mathrm H(a)&=a\circ \mathrm d_\mathrm H,\\
a\cap D&=(-1)^{|a|}a(D\otimes_R\bullet),\ a\in C_\mathrm{poly}^L(R),\ D\in D_\mathrm{poly}^L(R).
\end{aligned}
\]
We observe that~\eqref{ref-B.6-143} implies that $\mathrm b_\mathrm
H$, the Hochschild differential on $L$-chains, is compatible with
respect to \eqref{ref-B.3-139}, and that~\eqref{ref-B.10-147}
and~\ref{ref-B.11-148}, in the special case $D_1=\mu$, imply that
$\mathrm b_\mathrm H$ satisfies Leibniz's rule with respect to~\eqref{ref-B.4-140}
and~\eqref{ref-B.5-141} respectively.

Thus, combining these arguments with Proposition~\ref{ref-3.6-43}, we have the following important
\begin{Thm}\label{ref-B.2-150}
  For a Lie algebroid $L$ over the ring $R$ as above, the twist
  of~\eqref{ref-B.3-139}, \eqref{ref-B.4-140}, \eqref{ref-B.5-141} and of the
  Hochschild differential $\mathrm b_\mathrm H$ with respect to the
  isomorphism~\eqref{ref-3.15-44} endow $C_\mathrm{poly}^L(R)$ with a structure of
  precalculus up to homotopy over the Gerstenhaber algebra
  $(D_\mathrm{poly}^L(R),\mathrm d_\mathrm H,[\ ,\ ],\cup)$ up to
  homotopy.
\end{Thm}

\end{document}

%% file: S_graph.pstex_t
\begin{picture}(0,0)%
\includegraphics{S_graph.pstex}%
\end{picture}%
\setlength{\unitlength}{3947sp}%
\begingroup\makeatletter\ifx\SetFigFont\undefined%
\gdef\SetFigFont#1#2#3#4#5{%
  \reset@font\fontsize{#1}{#2pt}%
  \fontfamily{#3}\fontseries{#4}\fontshape{#5}%
  \selectfont}%
\fi\endgroup%
\begin{picture}(6016,5791)(2536,-7201)
\put(4786,-5986){\makebox(0,0)[lb]{\smash{\SetFigFont{25}{30.0}{\rmdefault}{\mddefault}{\updefault}{\color[rgb]{0,0,0}$\gamma$}%
}}}
\put(5086,-2491){\makebox(0,0)[lb]{\smash{\SetFigFont{25}{30.0}{\rmdefault}{\mddefault}{\updefault}{\color[rgb]{0,0,0}$\gamma$}%
}}}
\put(7561,-3571){\makebox(0,0)[lb]{\smash{\SetFigFont{25}{30.0}{\rmdefault}{\mddefault}{\updefault}{\color[rgb]{0,0,0}$\gamma$}%
}}}
\put(6496,-5281){\makebox(0,0)[lb]{\smash{\SetFigFont{25}{30.0}{\rmdefault}{\mddefault}{\updefault}{\color[rgb]{0,0,0}$\gamma_1$}%
}}}
\put(5776,-7201){\makebox(0,0)[lb]{\smash{\SetFigFont{25}{30.0}{\rmdefault}{\mddefault}{\updefault}{\color[rgb]{0,0,0}$a_0$}%
}}}
\put(8176,-5881){\makebox(0,0)[lb]{\smash{\SetFigFont{25}{30.0}{\rmdefault}{\mddefault}{\updefault}{\color[rgb]{0,0,0}$a_1$}%
}}}
\put(3091,-2116){\makebox(0,0)[lb]{\smash{\SetFigFont{25}{30.0}{\rmdefault}{\mddefault}{\updefault}{\color[rgb]{0,0,0}$a_3$}%
}}}
\put(2536,-5686){\makebox(0,0)[lb]{\smash{\SetFigFont{25}{30.0}{\rmdefault}{\mddefault}{\updefault}{\color[rgb]{0,0,0}$a_4$}%
}}}
\put(8371,-2926){\makebox(0,0)[lb]{\smash{\SetFigFont{25}{30.0}{\rmdefault}{\mddefault}{\updefault}{\color[rgb]{0,0,0}$a_2$}%
}}}
\put(5791,-3991){\makebox(0,0)[lb]{\smash{\SetFigFont{25}{30.0}{\rmdefault}{\mddefault}{\updefault}{\color[rgb]{0,0,0}$0$}%
}}}
\end{picture}

%% file: BNE_curconf.pstex_t
\begin{picture}(0,0)%
\includegraphics{BNE_curconf.pstex}%
\end{picture}%
\setlength{\unitlength}{3947sp}%
\begingroup\makeatletter\ifx\SetFigFont\undefined%
\gdef\SetFigFont#1#2#3#4#5{%
  \reset@font\fontsize{#1}{#2pt}%
  \fontfamily{#3}\fontseries{#4}\fontshape{#5}%
  \selectfont}%
\fi\endgroup%
\begin{picture}(4114,4450)(876,-6946)
\put(2896,-4456){\makebox(0,0)[lb]{\smash{\SetFigFont{20}{24.0}{\rmdefault}{\mddefault}{\updefault}{\color[rgb]{0,0,0}$0$}%
}}}
\put(2881,-6946){\makebox(0,0)[lb]{\smash{\SetFigFont{20}{24.0}{\rmdefault}{\mddefault}{\updefault}{\color[rgb]{0,0,0}$\overline 1$}%
}}}
\put(3046,-5581){\makebox(0,0)[lb]{\smash{\SetFigFont{20}{24.0}{\rmdefault}{\mddefault}{\updefault}{\color[rgb]{0,0,0}$1$}%
}}}
\end{picture}